\documentclass[11pt]{smfart}
\usepackage{amssymb}
\usepackage{amsmath, amsthm, amsopn, amsfonts}
\usepackage[dvips]{graphicx}
\usepackage{graphpap}
\usepackage[english, francais]{babel}
\newtheorem{theoreme}{Theorem}
\renewcommand{\theequation}{\arabic{equation}}
\newtheorem{proposition}{Proposition}
\newtheorem{lemme}[proposition]{Lemma}

\newtheorem{remarque}[proposition]{Remark}

\numberwithin{equation}{section}
\numberwithin{proposition}{section}

\def\11{{\rm 1~\hspace{-1.4ex}l} }
\def\R{\mathbb R}
\def\C{\mathbb C}
\def\Z{\mathbb Z}
\def\N{\mathbb N}
\def\E{\mathbb E}

\def\O{\mathbb O}

\begin{document}
\selectlanguage{english}
\title[NLS and invariant measures]
{ Invariant measures for the Nonlinear Schr\"odinger equation on the disc }
\author{Nikolay Tzvetkov}
\address{D\'epartement de Math\'ematiques, Universit\'e Lille I, 59 655
Villeneuve d'Ascq Cedex, France}
\email{nikolay.tzvetkov@math.univ-lille1.fr}
\begin{abstract} 
We study Gibbs measures invariant under the flow of the NLS on the unit disc of $\R^2$.
For that purpose, we construct the dynamics on a phase space of
limited Sobolev regularity and a wighted Wiener measure invariant by the NLS
flow. The density of the measure is integrable with respect to the   Wiener
measure for sub cubic nonlinear interactions. 
The existence of the dynamics is obtained in Bourgain spaces of low
regularity. The key ingredient are bilinear Strichartz estimates for the free evolution.
The bilinear effect in our analysis results from simple properties of the Bessel
functions and estimates on series of Bessel functions.
\end{abstract}
\subjclass{ 35Q55, 35BXX, 37K05, 37L50, 81Q20 }
\keywords{nonlinear Schr\"odinger, eigenfunctions, dispersive equations,
invariant measures}
\maketitle
\selectlanguage{english}
%
\section{Introduction}
This work fits in the line of research initiated in \cite{BGT1} aiming to study
the possible extensions of the work of Bourgain on nonlinear
Schr\"odinger equation (NLS) posed on the flat torus to other compact manifolds.
We are concerned here with the long time behavior of solutions of the nonlinear
Schr\"odinger equation, posed on the unit disc of $\R^2$. 
Our aim is to construct the dynamics on a phase space of
limited Sobolev regularity and a wighted Wiener measure invariant by the NLS flow.
Consider the Nonlinear Schr\"odinger equation
\begin{equation}\label{1}
iu_t+\Delta u +F(u)=0, 
\end{equation}
where 
$
u(t)\,:\, \Theta\longrightarrow \C
$
is a function defined on the unit disc
$$
\Theta=\{(x_1,x_2)\in\R^2 \, : \, x_1^2+x_2^2<1\}\, .
$$
The nonlinear interaction in (\ref{1}) is induced by $F(z)$, $z\in\C$ which is a smooth (non linear) complex valued function. 
We also assume that $F(0)=0$ and $F=\bar{\partial}V$ with a real valued $V$ satisfying the gauge invariance
assumption
$$
V(e^{i\theta} z)=V(z),\quad \forall\, \theta\in \R,\,\, \forall\, z\in\C\, .
$$
In addition, we suppose that for some $\alpha>0$,
\begin{equation}\label{rast}
\big|
\partial^{k_1}\bar{\partial}^{k_2}V(z)
\big|\leq C_{k_1,k_2}\langle z\rangle^{2+\alpha-k_1-k_2}\,.
\end{equation}
The real number $\alpha$ involved in (\ref{rast}) corresponds to the ``degree'' of
the nonlinear interaction. 
A typical example for $F(u)$ is
$$F(u)=\pm\Big(1+|u|^{2}\Big)^{\alpha/2}u$$
or $|u|^{\alpha}u$ when $\alpha$ is an even integer. In this paper, we assume that the nonlinearity is sub-cubic which means that
\begin{equation}\label{subcubic}
\alpha<2\,.
\end{equation}
Assumption (\ref{subcubic}) on $\alpha$ will be assumed from now on in the rest of this paper.
Notice that we do not suppose the defocusing assumption which in the context
of (\ref{1}) would be of type $V\leq 0$. 
In the (easier) defocusing case, one
can expect to cover a larger set of possible
values of $\alpha$ (see the final remarks at the end of the paper). 
\\

It is important that the problem (\ref{1}) may, at last formally, be seen as the
Hamiltonian PDE 
$$
iu_{t}=\partial_{\bar{u}}H(u,\overline{u})
$$
in an infinite dimensional phase space, with Hamiltonian
\begin{equation}\label{energy}
H(u,\bar{u})=\int_{\Theta}|\nabla u|^{2}-\int_{\Theta}V(u)
\end{equation}
and canonical coordinates $(u,\bar{u})$.
\\

We are interested in the solutions of the initial boundary value problem associated
to (\ref{1}). This means that we study (\ref{1}) subject to an initial 
condition
\begin{equation}\label{2}
u(0,x_1,x_2)=u_0(x_1,x_2),
\end{equation}
where $u_0$ is a given function.
In this paper, we will only consider initial data of Sobolev regularity $<1/2$
and thus we will not need to specify the boundary conditions on $\R\times
\partial \Theta$, where
$$
\partial \Theta=\{(x_1,x_2)\in\R^2 \, : \, x_1^2+x_2^2=1\}
$$
is the border of $\Theta$ (see also Remark~\ref{granitza} below).
We will however use the Dirichlet eigenfunctions of $\Delta$ as basis of $L^2(\Theta)$
and this will be convenient for our well-posedness analysis of (\ref{1})-(\ref{2}).
\\

We will only consider radial solutions of (\ref{1}), i.e. solutions
depending only on $t$ and $x_1^2+x_2^2$.
Thus, we suppose that the data is radially symmetric, i.e.
\begin{equation}\label{tilda}
u_0(x_1,x_2)=\tilde{u}_0(r),
\end{equation}
where 
$$
x_1=r\cos\varphi,\quad x_2=r\sin\varphi,\qquad 0\leq r<1,\,\, \varphi\in
[0,2\pi]\, .
$$
Let $J_0$ be the Bessel function of order zero (see e.g. \cite{S}) and let 
$z_1,z_2,\dots$ be the zeros of $J_0$. We have that
$$
0<z_1<z_2<\dots z_n< \dots
$$
and the zeroes are simple.
We also have that $J_{0}(z_n r)$ are eigenfunctions of
the Dirichlet self adjoint realization of $-\Delta$, corresponding to eigenvalues $z_n^2$.
Moreover any $L^2(\Theta)$ radial function can be expanded with respect to
$J_{0}(z_n r)$.
Let us set
\begin{equation}\label{en}
e_n\equiv e_{n}(r)=\|J_0(z_n\cdot)\|_{L^2(\Theta)}^{-1}\, J_0(z_n r)\,
\end{equation}
and
$$
e_{n,s}=z_{n}^{-s}e_{n}\, .
$$
We can decompose the solutions of (\ref{1}) with data of type (\ref{tilda}) as 
$$
u(t)=\sum_{n\geq 1}c_n(t)\, e_{n,s}\, .
$$
The initial data is thus given by
$$
\tilde{u}_0=\sum_{n\geq 1}c_n(0)e_{n,s}\, ,
$$
i.e the initial data is uniquely determined from the sequence $(c_n(0))$, $n\in\N$.
Thus the equation (\ref{1}) 
can be written as
\begin{equation}\label{4}
iz_{n}^{-s}\dot{c_n}(t)-z_n^2\,z_{n}^{-s} c_n(t)+
\Pi_{n}\Big(
F\Big(
\sum_{m\geq 1}c_m(t)\, e_{m,s}\Big)
\Big)=0,
\end{equation}
$n\geq 1$, where $\Pi_n$ is the projection on the mode $e_n$.
For instance if $f\in
L^{1}(\Theta)$ (which will always be the case in this paper), we have
$$
\Pi_n(f)=\langle f, e_n\rangle =\int_{\Theta} f\, \overline{e_n}\,.
$$
Of course one can define the action of $\Pi_n$ on distributions but here we
will not need it.
Notice that if $f\in L^2(\Theta)$, 
$\Pi_n(f)$ is simply the $L^2(\Theta)$ scalar product of $f$ and $e_n$.
Formally, equation (\ref{4}) is in fact a Hamiltonian equation with with canonical
coordinates $(c,\overline{c})$ and Hamiltonian
$$
H(c,\overline{c})=\sum_{n\geq 1} z_n^{2-2s}\, |c_n|^2 - \int_{0}^1
V\Big(
\sum_{m\geq 1}c_m\, e_{m,s}(r)\Big)rdr\, ,
$$
where $c=(c_n)$, $n\in\N$. More precisely equation (\ref{4}) can be written as
$$
ic_t=J\frac{\delta H}{\delta \overline{c}},\quad i\overline{c}_{t}=-J\frac{\delta H}{\delta c}\, ,
$$
where $\delta$ denotes the variational derivative and $J={\rm diag}(z_n^{2s})_{n\geq 1}$ is the
map inducing the symplectic form in the coordinates $(c,\overline{c})$. The
only important consequence, for our analysis, of this discussion is that
$H(c,\overline{c})$ is , at least formally, conserved by the flow of  (\ref{4}). 
\\

Let us now describe the construction of Lebowitz-Rose-Speer (cf. \cite{LRS})
of a weighted Wiener measure which is at least formally invariant under under
flow of (\ref{1}). The rigorous justification of the invariance of the measure 
will require, among other things, a new well-posedness result for the initial
value problem (\ref{1})-(\ref{2}). Let us fix a real number $s$ such that
\begin{equation}\label{s}
0<s<\frac{\alpha}{\alpha+2}\,.
\end{equation}
This number $s$ will be fixed from now on in all the rest the paper.
Notice that thanks to the restriction  (\ref{subcubic}) on the degree of the
nonlinearity $\alpha$,
$$
\frac{\alpha}{\alpha+2}<\frac{1}{2}\,.
$$

For $\sigma\in [0,1/2[$, let us denote by $H^{\sigma}_{rad}(\Theta)$ the Sobolev space 
of radial functions on $\Theta$, i.e.
$
u\in H^{\sigma}_{rad}(\Theta)
$
if and only if
$$
u=\sum_{n\geq 1}c_n e_{n,s}, \quad c_n\in\C
$$
with
$$
\sum_{n\geq 1}z_{n}^{2(\sigma-s)} |c_n|^{2}<\infty\, .
$$
The Sobolev space $H^s_{rad}(\Theta)$ is naturally a complex Hilbert space with orthonormal
basis $e_{n,s}$. 
Denote by $(\cdot,\cdot)$ the scalar product in $H^s_{rad}(\Theta)$.
Our goal will be to construct a well defined (at least local in time) dynamics on $H^s_{rad}(\Theta)$
and to construct a bounded Borel measure on it, invariant by the flow of (\ref{1}).
\\

The free Hamiltonian is given by
$$
H_0(c,\overline{c})=\sum_{n\geq 1}z_n^{2-2s}\, |c_n|^2 .
$$
It turns out that a renormalization of the formal measure
$$
e^{-H_0(c,\overline{c})}d^2c=
\prod_{n\geq 1}e^{-z_n^{2-2s}|c_n|^2}d^2 c_n
$$
is a Wiener measure. 
More precisely, we can give a sense of the formal measure
\begin{equation*}
\frac{e^{-H_0(c,\overline{c})}d^2c}{\int e^{-H_0(c,\overline{c})}d^2c} 
=
\prod_{n\geq 1}\frac{e^{-z_n^{2-2s}|c_n|^2}d^2 c_n}
{\int_{\C} e^{-z_n^{2-2s}|c_n|^2}d^2 c_n} 
\end{equation*}
as a measure on the Hilbert space $H^s_{rad}(\Theta)$
(corresponding to a Gaussian distribution for each mode).

A set $U\subset H^s_{rad}(\Theta)$ is called cylindrical if there exists
$N\in\N$ and a Borel set $V\subset \C^{N}$ such that
\begin{equation}\label{cyl}
U=
\Big\{
u\in H^s_{rad}(\Theta)\, :\, 
\big( (u,e_{1,s}),\dots,(u,e_{N,s})\big)\in V
\Big\}.
\end{equation}
Let us denote by $\tilde{\mu}$ the measure, defined on the cylindrical sets $U$ determined by
(\ref{cyl}) as
\begin{multline}\label{cyl-mes}
\tilde{\mu}(U)=
\frac
{\int_{V}e^{-{\sum_{1\leq n\leq N} z_n^{2-2s}|c_n|^2}}d^{2}c_1\dots d^{2}c_{N}}
{\int_{\C^{N}}e^{-{\sum_{1\leq n\leq N} z_n^{2-2s}|c_n|^2}}d^{2}c_1\dots d^{2}c_{N}}
\\
=
\pi^{-N}\Big(\prod_{1\leq n\leq N}z_{n}^{2-2s}\Big)
\int_{V}e^{-{\sum_{1\leq n\leq N} z_n^{2-2s}|c_n|^2}}d^{2}c_1\dots d^{2}c_{N}\,.
\end{multline}
The cylindrical sets form an algebra in $H^s_{rad}(\Theta)$.
Moreover the minimal sigma algebra containing all cylindrical sets is the
Borel sigma algebra. Since (see (\ref{zeros}) below) $z_n\sim n$, we deduce
that the series
$$
\sum_{n\geq 1}z_{n}^{2s-2}
$$
converges. 
It implies that the linear map defined on $H^s_{rad}(\Theta)$ by
$$
e_{n,s}\longmapsto z_{n}^{2s-2}e_{n,s}
$$
is in the trace class.
Therefore (see e.g. \cite{Mourier,DY,Zh}) the measure $\tilde{\mu}$ is
countably additive on the cylindrical sets of $H^s_{rad}(\Theta)$. 
We then denote by $\mu$ the Borel probability measure on $H^s_{rad}(\Theta)$ which is the unique
extension (Caratheodory theorem) of $\tilde{\mu}$ to the Borel sigma algebra of $H^s_{rad}(\Theta)$.
For the sake of completeness, in Section~3 we present the proof of the
countable additivity $\tilde{\mu}$ on the algebra of the cylindrical sets of $H^s_{rad}(\Theta)$. 
As we will show in Proposition~\ref{high}, for $\sigma\in[s,1/2[$,
$\mu(H^{\sigma}_{rad}(\Theta))=1$ and thus, we may consider $\mu$ as a measure on the space
\begin{equation}\label{X}
{\mathcal X}=\bigcap_{s\leq\sigma<\frac{1}{2}}H^{\sigma}_{rad}(\Theta)\,.
\end{equation}
Thus one should not take the particular choice of $s$ that we made too seriously.
Notice that since $\sigma<1/2$ the boundary conditions are not of importance
in the definition of ${\mathcal X}$.
In addition, in (\ref{X}) the intersection may be assumed countable.
\\

One may hope that the expression
$
\exp(\int_{\Theta} V(u))d\mu(u)
$
which is a normalised version of the formal Gibbs measure 
$\exp(-H(u,\bar{u}))\,d^{2}u$
is a well defined measure. 
The expression $\exp(-H(u,\bar{u}))\,d^{2}u$ is formally
invariant by the flow thanks to the Hamiltonian
conservation. If we were in finite dimensions the invariance would follow from the
invariance of the Lebesgue measure by the flow (Liouville's theorem).  
There is however a problem with the integrability
of the above density with respect to $\mu$. We will solve this problem by
using the $L^2$ cut-off idea of Lebowitz-Rose-Speer \cite{LRS}. 
\\

This paper is devoted to the proof of the following statement.
\begin{theoreme}\label{glavna}
Fix $R>0$.
Let us denote by $\chi\, :\, \R\rightarrow\{0,1\}$ the characteristic function of the 
interval $[0,R]$. 
For $u\in {\mathcal X}$, we define the functional $f(u)$ by
$$
f(u)=\chi\big(\|u\|_{L^2(\Theta)}\big)
\exp\Big(\int_{\Theta} V(u)\Big)\,.
$$
Then for every $q\in [1,+\infty[$, 
\begin{equation}\label{parva}
f(u)\in L^{q}(d\mu(u))\,.
\end{equation}
Moreover, if we set $d\rho(u)=f(u)d\mu(u)$ then there exists a set $\Sigma$ of
full $\rho$ measure such that for every $u_0\in\Sigma$ the Cauchy problem (\ref{1})-(\ref{2}) 
has a unique (in a suitable functional framework) global in time solution.
Finally, if we denote by $\Phi(t)$, $t\in\R$ the flow of (\ref{1}) acting on
$\Sigma$ then the measure $\rho$ is invariant under the flow of (\ref{1}),
i.e. for every $\rho$ measurable set $A\subset\Sigma$, every $t\in\R$,
$
\rho(A)=\rho(\Phi(t)(A)).
$
\end{theoreme}
\begin{remarque}
The uniqueness statement of Theorem~\ref{glavna} can be precised as follows :
for every $T>0$ there exists a Banach space $X_{T}$ continuously embedded in
$C([-T,T];H^s_{rad}(\Theta))$ such that the solution of (\ref{1}) with data
$u_0\in\Sigma$ is unique in $X_{T}$. 
\end{remarque}
Notice that thanks to the growth assumption (\ref{rast}) and the Sobolev
embedding, the functional $f(u)$ is well-defined for $u\in {\mathcal X}$. 
\\

As a consequence of Theorem~\ref{glavna}, the Poincar\'e recurrence theorem implies that
almost surely on the support of $d\rho$ the solutions of (\ref{1}) are stable
according to Poisson (see \cite{Zh} and the references therein for more details).
\\

Similar results to  Theorem~\ref{glavna} in the case of the circle $S^1$ are
known thanks to the works \cite{Bo,Zh}. 
Gibbs type invariant measures for a Wicked ordered cubic defocusing NLS,
posed on the two dimensional rational torus are constructed in \cite{Bo2}.
\\

Invariant measures for defocusing NLS of type (\ref{1}) posed on an arbitrary
compact riemannian manifold are constructed in \cite{KS}. These measures are
not of Gibbs type (but still related to the conservation laws), 
and are living on functions in the Sobolev space $H^2$.
Let us notice that Dirac measures concentrated on a stationary (independent of
$t$) solutions of (\ref{1}) are clearly invariant. The measures constructed in \cite{KS} are not of
this trivial type since the defocusing nature of the problem excludes the existence of
stationary solutions.
\\

The proof of Theorem~\ref{glavna} follows the ideas developed by Zhidkov (see \cite{Zh} and the references by the same author therein) 
and Bourgain \cite{Bo}. 
The main difficulties we should overcome are to prove a new local well-posedness
results for (\ref{1}), posed on the unit disc as well as adapting some estimates on random
Fourier series to the case of functions on the unit disc of $\R^2$.
In the local well-posedness analysis, we need some bilinear Strichartz estimates. 
Starting from the work of Bourgain, estimates in this spirit were already used by many authors in the context of 
dispersive PDE's. In the analysis here, the crucial bilinear effect results from simple properties of the Bessel
functions and estimates on some series of Bessel functions.
Notice that the bilinear approach and the Bourgain spaces are needed to be
employed here since the well-posedness analysis of \cite{BGT1} based only on linear Strichartz inequalities and
Sobolev spaces requires the restriction $\sigma>1/2$ (thus missing ${\mathcal X}$) coming from the Sobolev
embedding $W^{\sigma,4}\subset L^{\infty}$, 
$\sigma>1/2$ in two dimensions.
\\

The rest of the paper is organized as follows. The next section is devoted to
some properties of the Bessel functions needed for our analysis of NLS
(\ref{1}). In Section~3, we collect some properties of Wiener 
type measures on Sobolev spaces of radial functions on the disc. Section~4 is devoted to
bilinear Strichartz type inequalities which are the basic analytical tool in
this paper. In Section~5, we introduce the Bourgain spaces of radial functions
on $\Theta$. The main nonlinear estimate are established in Section~6.
As a first consequence of these estimates, in Section~7 we prove some local well-posedness
results for NLS and its finite dimensional approximation. Next, in section~8,
we improve the result for the finite dimensional model. In Section~9, we
transfer the result of Section~8 to the NLS. The proof of
Theorem~\ref{glavna}, we will be completed in Section~10. The final section is
devoted to some straightforward extensions of Theorem~\ref{glavna} and open problems
that seem of interest to the author of the present paper.
\\

{\bf Notation.} 
Let us now introduce several notations that will be used in the paper. For two
positive real numbers $N_1$ and $N_2$, we denote by
$N_1\wedge N_2\equiv\min(N_1,N_2)$ the smaller one.
For $x\in \R$, we set $\langle x\rangle\equiv 1+|x|$.
We use the notations $\sim$ or $\approx$ for the equivalence of two
quantities, uniformly with respect to some parameters which will be clear in each appearance of these two symbols. 
Several positive constants uniform with respect to some parameters, which will be clear
in each appearance, will be denoted by $C$ or $c$. The parameter set will
always be a set of numbers or a set of functions.
\section{On the Bessel functions and their zeros}
In this section, we collect several facts on the zero order Bessel function
that will be used in the sequel.
These facts are essentially in the literature (see e.g. \cite{S,W}) but, in order to keep the paper
as self contained as possible, here we give the proofs. We will be interested on $J_{0}(x)$
for $x\geq 0$ and its zeros $z_n$ since $J_{0}(z_n r)$, $0\leq r<1$ form a
basis for the radial $L^2$ functions on the disc $\Theta$. The Bessel function
$J_0(x)$ is defined by
\begin{equation*}
J_{0}(x)=\sum_{j=0}^{\infty}\frac{(-1)^{j}}{(j!)^{2}}\Big(\frac{x}{2}\Big)^{2j}\,.
\end{equation*}
The function $J_0(x)$ solves the ordinary differential equation
\begin{equation*}
J''_{0}(x)+\frac{1}{x}J'_{0}(x)+J_{0}(x)=0\,.
\end{equation*}
The function $J_{0}(x)$ may be seen as the zero Fourier coefficient of the
function $\exp(ix\sin\theta)$, $\theta\in[-\pi,\pi]$ and thus
\begin{equation*}
J_{0}(x)=\frac{1}{2\pi}\int_{-\pi}^{\pi}e^{ix\sin\theta}d\theta\,.
\end{equation*}
Moreover, by the Lebesgue differentiation theorem,
\begin{equation*}
J'_{0}(x)=\frac{1}{2\pi}\int_{-\pi}^{\pi}(i\sin\theta)e^{ix\sin\theta}d\theta\,.
\end{equation*}
Recall that $e_n\,:\,\Theta\rightarrow \R$, defined by (\ref{en})
form an orthonormal basis of the $L^2$ radial functions on the disc $\Theta$.
Observe that $e_{n}(r)$ are real valued.
The next lemma provides $L^p(\Theta)$ bounds for $e_n$ in the regime $n\gg 1$.
\begin{lemme}\label{Lp}
Let $p\in [2,\infty]$. 
There exists $C$ such that for every $n\geq 1$,
\begin{equation}\label{eigen}
\|e_{n}\|_{L^{p}(\Theta)}\leq C\delta(n)\|e_{n}\|_{L^{2}(\Theta)}= C\delta(n)  ,\quad
\|e'_{n}\|_{L^{p}(\Theta)}\leq C\delta(n)\|e'_{n}\|_{L^{2}(\Theta)}\,,
\end{equation}
where
$$
\delta(n)=
\left\{
\begin{array}{ll}
1 & {\rm when }\quad 2\leq p<4,
\\
(\log(1+n))^{\frac{1}{4}} & {\rm when }\quad p=4,
\\
n^{-\frac{2}{p}+\frac{1}{2}} & {\rm when }\quad p>4\,.
\end{array}
\right.
$$
In particular for every $\varepsilon>0$ there exists $C_{\varepsilon}$ such
that for every $n_1,n_2\geq 1$,
\begin{equation}\label{util}
\|e_{n_1}e_{n_2}\|_{L^2(\Theta)}\leq
C_{\varepsilon}
(\min(n_1,n_2))^{\varepsilon}\|e_{n_1}\|_{L^2(\Theta)}\|e_{n_2}\|_{L^2(\Theta)}
=C_{\varepsilon}
(\min(n_1,n_2))^{\varepsilon}
\end{equation}
and
\begin{equation}\label{utilbis}
\|e_{n_1}e'_{n_2}\|_{L^2(\Theta)}\leq
C_{\varepsilon}
(\min(n_1,n_2))^{\varepsilon}\|e_{n_1}\|_{L^2(\Theta)}\|e'_{n_2}\|_{L^2(\Theta)}\,\, .
\end{equation}
Finally, there exist two positive constants $C_1$ and $C_2$ such that for
every $n\in\N$,
\begin{equation}\label{equiv}
C_1 n=C_1 n\|e_{n}\|_{L^{2}(\Theta)}\leq\|e'_{n}\|_{L^{2}(\Theta)}\leq C_2
n\|e_{n}\|_{L^{2}(\Theta)}=C_2 n.
\end{equation}
\end{lemme}
\begin{proof}
The proof is based on the asymptotics for $J_0(x)$ and $J'_{0}(x)$ for large
values of $x$. These asymptotics may be found by applying the stationary phase
formula to the integrals defining $J_0(x)$ and $J'_{0}(x)$. Indeed, in both
cases the phase $\sin\theta$ has two non-degenerate critical points $\pm
\frac{\pi}{2}$ on $[-\pi,\pi]$. Therefore, there exists $C>0$ and a function
$r_1(x)$ defined on $[1,+\infty[$ such that
$$
J_{0}(x)=\sqrt{\frac{2}{\pi}}\,\frac{\cos\big(x-\frac{\pi}{4}\big)}{\sqrt{x}}+r_1(x),\quad
|r_{1}(x)|\leq Cx^{-\frac{3}{2}}
$$
(the two critical points contribute with phases $\exp(i(\pm x\mp\pi/4))$).
Similarly, we have
$$
J'_{0}(x)=-\sqrt{\frac{2}{\pi}}\,\frac{\sin\big(x-\frac{\pi}{4}\big)}{\sqrt{x}}+\tilde{r}_1(x),\quad
|\tilde{r}_{1}(x)|\leq Cx^{-\frac{3}{2}}\,.
$$
A first consequence of the above representations of $J_0(x)$ and $J'_{0}(x)$
is that the $n$'th zero of $J_0(x)$ satisfy $z_n\sim n$. We can therefore
write that for $n\gg 1$
\begin{eqnarray*}
\|J_{0}(z_n\cdot)\|_{L^2(\Theta)}^{2} & = & 
\int_{0}^{1}|J_{0}(z_n r)|^{2}r\,dr
\\
& = & z_{n}^{-2}\int_{0}^{z_n}|J_{0}(\rho)|^{2}\rho d\rho
\\
& \geq & 
cn^{-2}\int_{0}^{cn}|J_{0}(\rho)|^{2}\rho d\rho
\\
& \geq &
Cn^{-2}\int_{1}^{cn}\Big(\frac{1+\cos(2\rho-\pi/2)}{2\rho}
-\frac{C}{\rho^2}\Big)\rho d\rho
\\
& \geq & 
cn^{-2}(cn-C\log(n))\geq cn^{-1}\,.
\end{eqnarray*}
Therefore
\begin{equation}\label{down}
\|J_{0}(z_n\cdot)\|_{L^2(\Theta)}\geq cn^{-1/2}.
\end{equation}
Similarly, we can show that
$$
\|J'_{0}(z_n\cdot)\|_{L^2(\Theta)}\geq cn^{-1/2}.
$$
On the other hand, using that $|J_{0}(x)|\leq Cx^{-1/2}$, $x\geq 1$,
and, $|J_{0}(x)|\leq C$, $x\leq 1$, we obtain that for $p\in[2,\infty[$,
\begin{eqnarray*}
\|J_{0}(z_n\cdot)\|_{L^p(\Theta)}^{p} & = &
\int_{0}^{1}|J_{0}(z_n r)|^{p}r\,dr
\\
& = &
z_{n}^{-2}\int_{0}^{z_n}|J_{0}(\rho)|^{p}\rho\,
d\rho
\\
& \leq & C n^{-2}
\Big(C+\int_{1}^{cn}
\rho^{-p/2}\,\rho \, d\rho\Big)
\end{eqnarray*}
which gives the bound (\ref{eigen}) 
for $e_n$ and $p<+\infty$
by distinguishing the three regimes for $p$
involved in the definition of $\delta(n)$.
The last estimate also implies that
$$
\|J_{0}(z_n\cdot)\|_{L^2(\Theta)}\leq Cn^{-1/2}
$$
and thus
$$
\|J_{0}(z_n\cdot)\|_{L^2(\Theta)}\sim n^{-1/2}.
$$
Estimate (\ref{eigen}) for $p=\infty$ and $e_n$ follows form the bound
$|J_0(x)|\leq C$ for all $x\geq 0$ and the inequality (\ref{down}).
This completes the proof of (\ref{eigen}) as far as $e_n$ is concerned. The
bound for $e'_{n}$ in (\ref{eigen}) can be established in a completely
analogous way, once we have the stationary phase approximation of
$J'_{0}(x)$.
We also have 
$$
\|J'_{0}(z_n\cdot)\|_{L^2(\Theta)}\leq Cn^{-1/2}
$$
and thus
$$
\|J'_{0}(z_n\cdot)\|_{L^2(\Theta)}\sim n^{-1/2}.
$$
Since 
$$
e'_{n}(r)=z_{n}
\|J_{0}(z_n\cdot)\|_{L^2(\Theta)}^{-1}
J'_{0}(z_n r)
$$
we get estimate (\ref{equiv}).
Finally, the assertion of (\ref{util}) results from (\ref{eigen}) and
H\"older inequality
$$
\|e_{n_1}e_{n_2}\|_{L^2(\Theta)}\leq 
\|e_{n_1}\|_{L^p(\Theta)}\|e_{n_2}\|_{L^q(\Theta)},\quad
\frac{1}{p}+\frac{1}{q}=\frac{1}{2},
$$
with $p,q$ close to $4$ and according to the order of $n_1$, $n_2$, the bigger
of $p,q$ is attached to the smaller of $n_1$, $n_2$. A similar argument yields
(\ref{utilbis}).
This completes the proof of Lemma~\ref{Lp}. 
\end{proof}
The next lemma provides a more precise asymptotics for the zeros $z_n$, $n\gg 1$.
\begin{lemme}\label{lemmazero}
For every $\kappa>0$ there exists $C>0$ such that the zeros of $J_0(x)$ can be written as
\begin{equation}\label{zeros}
z_{n}=\pi\big(n-\frac{1}{4}\big)+\frac{1}{8\pi\big(n-\frac{1}{4}\big)}+r(n),\quad
|r(n)|\leq Cn^{-(2-\kappa)}\,.
\end{equation}
\end{lemme}
\begin{remarque}
In fact, much better bounds on $r(n)$ may be proved. However, estimate
(\ref{zeros}) will be sufficient for our applications.
\end{remarque}
\begin{proof}[Proof of Proposition~\ref{lemmazero}]
Using the stationary phase formula at order $2$ in the integral representation of $J_0(x)$
gives the existence of a constant $C>0$ and a function $r_2(x)$ defined on $[1,+\infty[$ such that
$$
J_{0}(x)=
\sqrt{\frac{2}{\pi}}\,\frac{\cos\big(x-\frac{\pi}{4}\big)}{x^{1/2}}
+
\sqrt{\frac{2}{\pi}}\,\frac{\sin\big(x-\frac{\pi}{4}\big)}{8x^{3/2}}
+
r_2(x),\quad
|r_{2}(x)|\leq Cx^{-\frac{5}{2}}\,.
$$
Therefore, for $n\gg 1$, the zero $z_n$ solves the equation $F(z_n)=0$, where
$F(x)$ (with $x-\pi/4$ near the positive odd integer multiples of $\pi/2$) is a continuous
function of the form 
$$
F(x)=\frac{1}{\tan\big(x-\frac{\pi}{4}\big)}+\frac{1}{8x}+{\mathcal O}(n^{-2}).
$$
Here ${\mathcal O}(n^{-2})$ denotes a quantity $\leq Cn^{-2}$ with $C$
independent of $n$ and $x$.
For $\kappa>0$, we set
$$
z_{n}^{\pm}=\pi\big(n-\frac{1}{4}\big)+\frac{1}{8\pi\big(n-\frac{1}{4}\big)}\pm\frac{1}{n^{2-\kappa}}.
$$
Further, we set
$$
\varepsilon_{n}^{\pm}=\frac{1}{8\pi\big(n-\frac{1}{4}\big)}\pm\frac{1}{n^{2-\kappa}}.
$$
Notice that $\cos(z_{n}^{\pm}-\pi/4)=(-1)^{n}\sin\varepsilon_{n}^{\pm}$ and
$\sin(z_{n}^{\pm}-\pi/4)=(-1)^{n+1}\cos\varepsilon_{n}^{\pm}$.
Therefore, by expanding, we get
$$
F(z_{n}^{\pm})=-\tan(\varepsilon_{n}^{\pm})
+\frac{1}{8\pi\big(n-\frac{1}{4}\big)+8\varepsilon_{n}^{\pm}}
+{\mathcal O}(n^{-2})=
\mp\frac{1}{n^{2-\kappa}}+{\mathcal O}(n^{-2}).
$$
Therefore for $n\gg 1$ the zero $z_n$ lies between $z_n^{-}$ and $z_n^{+}$.
This completes the proof of Lemma~\ref{lemmazero}.
\end{proof}
\section{The measures $\mu$ and $\rho$ }
In this section, we prove (\ref{parva}) and we collect some properties of the
measures $\mu$ and $\rho$.
Let us first observe that the minimal sigma algebra containing the algebra of cylindrical sets
(\ref{cyl}) contains the closed balls of $H^{\sigma}_{rad}(\Theta)$,
$\sigma\in[s,1/2[$.
Indeed, if for $r>0$ and $v\in H^{\sigma}_{rad}(\Theta)$, we set
$$
B_{\sigma}(r,v)=\big(u\in H^s_{rad}(\Theta)\,:\, u\in
H^{\sigma}_{rad}(\Theta)\quad  {\rm and} \quad\|u-v\|_{H^{\sigma}(\Theta)}\leq
r\big)
$$
then
$$
B_{\sigma}(r,v)=\bigcap_{N\geq 1} U_{\sigma,N}(r,v),
$$
where the cylindrical sets $U_{\sigma,N}(r,v)$ are defined by
$$
U_{\sigma,N}(r,v)=\Big(u\in H^s_{rad}(\Theta)\,:\, \sum_{1\leq j\leq
  N}z_{n}^{\sigma-s}|(u-v,e_{j,s})|^{2}\leq r^2\Big).
$$
Since $H^{s}_{rad}(\Theta)$ is separable, we obtain that
the minimal sigma algebra containing all cylindrical sets is the
Borel sigma algebra.
\\

As mentioned in the introduction, for a sake of completeness, we give the
proof of the countable additivity of the measure
$\tilde{\mu}$.
\begin{proposition}\label{mutilde}
The measure $\tilde{\mu}$, defined on the algebra of cylindrical sets
(\ref{cyl}) by formula (\ref{cyl-mes}) is countably additive, i.e. for every sequence 
$U_n$, $n\in\N$ of cylindrical sets such that $U_{n+1}\subset U_{n}$ and
\begin{equation}\label{lidl}
\bigcap_{n\in\N} U_n=\emptyset,
\end{equation}
one has 
$$
\lim_{n\rightarrow\infty}\tilde{\mu}(U_n)=0\,.
$$
Thus $\tilde{\mu}$ has a unique extension that we denote by $\mu$ to the Borel sigma algebra
of $H^s_{rad}(\Theta)$ which is a Borel probability measure on $H^s_{rad}(\Theta)$.
\end{proposition}
\begin{proof}
Let $\sigma>0$ be such that $s+\sigma<1/2$. For $R\geq 1$, we consider the set
$$
K_{R}=\big\{u\in H^s_{rad}(\Theta)\,:\,\|u\|_{H^{s+\sigma}(\Theta)}\leq R \big\}.
$$
Thanks to the compactness of the embedding $H^{s+\sigma}_{rad}(\Theta)$ into
$H^s_{rad}(\Theta)$, we obtain that $K_{R}$ is a compact set of
$H^s_{rad}(\Theta)$.
Since $U_n$, $n\in\N$ are cylindrical sets, there exists a function
$r:\N\rightarrow\N$ such that for every $n$ the set $U_n$ can be seen as a subset of the
finite dimensional space $E_{r(n)}$ defined by
$
E_{r(n)}=
{\rm span}(e_{j,s})_{1\leq j\leq r(n)}.
$
More precisely, there exists a Borel set $\widetilde{U}_n$ of $ E_{r(n)}$ 
such that
$$
U_{n}=\Big(u\in H^s_{rad}(\Theta)\,:\, (u,e_{1,s})e_{1,s}+\cdots +(u,e_{r(n),s})e_{r(n),s}\in \widetilde{U}_n\Big).
$$
Consider the cylindrical sets $F_{r(n)}$ defined as
$$
F_{r(n)}\equiv 
\Big(u\in H^s_{rad}(\Theta)\,:\, (u,e_{1,s})e_{1,s}+\cdots + (u,e_{r(n),s})e_{r(n),s}\in K_{R}\Big).
$$
Then 
\begin{equation}\label{koleda}
\tilde{\mu}(F_{r(n)})\geq 1-CR^{-2},
\end{equation}
where $C$ is a constant independent of $R$ and what is more important,
independent of $n\in\N$. 
Set $m=r(n)$. In order to prove (\ref{koleda}), we observe that
\begin{equation*}
1-\tilde{\mu}(F_{r(n)})\leq I,
\end{equation*}
where $I$ is given by the integral
$$
I=\pi^{-m}\Big(\prod_{j=1}^{m}z_{j}^{2-2s}\Big)
\int_{V}e^{-{\sum_{1\leq j\leq m} z_j^{2-2s}|c_j|^2}}d^{2}c_1\dots d^{2}c_{m}\,,
$$
where $V$ is given by
$$
V=
\big\{
(c_1,\dots,c_m)\in\C^{m}\,:\, \sum_{j=1}^{m}z_{j}^{2\sigma}|c_j|^{2}\geq R^2
\big\}.
$$
Set $\theta\equiv s+\sigma<1/2$. By the change of the variable
$c_{j}\rightarrow z_{j}^{\sigma}c_{j}$, we obtain that
$$
I=\pi^{-m}\Big(\prod_{j=1}^{m}z_{j}^{2-2\theta}\Big)
\int_{W}e^{-{\sum_{1\leq j\leq m} z_j^{2-2\theta}|c_j|^2}}d^{2}c_1\dots d^{2}c_{m}\,,
$$
where $W$ is given by
$$
W=
\big\{
(c_1,\dots,c_m)\in\C^{m}\,:\, \sum_{j=1}^{m}|c_j|^{2}\geq R^2
\big\}.
$$
By introducing polar coordinates in each $c_j$ integration, we may estimate
\begin{eqnarray*}
R^2\, I& \leq & \pi^{-m}\Big(\prod_{j=1}^{m}z_{j}^{2-2\theta}\Big)
\int_{W}\Big(\sum_{j=1}^{m}|c_j|^{2}\Big)e^{-{\sum_{1\leq j\leq m} z_j^{2-2\theta}|c_j|^2}}d^{2}c_1\dots d^{2}c_{m}
\\
& \leq &
\pi^{-m}\Big(\prod_{j=1}^{m}z_{j}^{2-2\theta}\Big)
\int_{\C^m}\Big(\sum_{j=1}^{m}|c_j|^{2}\Big)e^{-{\sum_{1\leq j\leq m} z_j^{2-2\theta}|c_j|^2}}d^{2}c_1\dots d^{2}c_{m}
\\
& = &\sum_{j=1}^{m}z_{j}^{2\theta-2}\leq C
\end{eqnarray*}
where $C$ is a constant independent of $m=r(n)$.
This proves (\ref{koleda}).
\\

Let us fix $\varepsilon>0$. By the regularity of the Lebesgue measure, 
using that $U_{n+1}\subset U_{n}$ one can construct closed sets
$\widetilde{V}_n\subset E_{r(n)}$ such that
\begin{equation}\label{jap}
V_{n}=\Big(u\in H^s_{rad}(\Theta)\,:\, (u,e_{1,s})e_{1,s}+\cdots
+(u,e_{r(n),s})e_{r(n),s}\in \widetilde{V}_n\Big)
\end{equation}
satisfy
$$
V_{n}\subset U_{n},\quad V_{n+1}\subset V_{n},\quad
\tilde{\mu}(U_n\backslash V_n)<\varepsilon/2\,.
$$
Indeed, one first constructs closed sets 
$\widetilde{W}_n\subset E_{r(n)}$ such that
\begin{equation*}
W_{n}=\Big(u\in H^s_{rad}(\Theta)\,:\, (u,e_{1,s})e_{1,s}+\cdots
+(u,e_{r(n),s})e_{r(n),s}\in \widetilde{W}_n\Big)
\end{equation*}
satisfy
$$
W_{n}\subset U_{n},\quad
\tilde{\mu}(U_n\backslash W_n)<\varepsilon/2^{n+3}\,.
$$
Then, we set
$$
V_{n}\equiv \bigcap_{j=1}^{n}W_{j}
$$
and one easily verifies that $V_n$ satisfies (\ref{jap}).
\\

We have that $K_{R}\cap V_n$ is a compact set of $H^s_{rad}(\Theta)$ included
in $U_n$. Therefore (\ref{lidl}) yields
$$
\bigcap_{n\in\N} (K_{R}\cap V_n)=\emptyset\,.
$$
Hence, there exists $N\geq 1$ such that $K_{R}\cap V_n=\emptyset$ for
$n\geq N$. Moreover, $F_{r(n)}\cap V_n=\emptyset$ for $n\geq N$.
Indeed, if $u\in F_{r(n)}\cap V_n$ then by setting
$$
u_{n}\equiv (u,e_{1,s})e_{1,s}+\cdots + (u,e_{r(n),s})e_{r(n),s},
$$
we observe that $u_{n}\in K_{R}$ and $u_{n}\in V_{n}$ which is a
contradiction.
Thus $F_{r(n)}\cap V_n=\emptyset$.
Therefore, using (\ref{koleda}), we infer that
$$
1\geq\tilde{\mu}( F_{r(n)}\cup V_n)=
\tilde{\mu}(F_{r(n)})+\tilde{\mu}(V_n)
\geq 1-CR^{-2}+\tilde{\mu}(V_n).
$$
Hence $\tilde{\mu}(V_n)\leq CR^{-2}$ and thus for $R\gg 1$ (independently of $n$),
$$
\tilde{\mu}(U_n)\leq 
\tilde{\mu}(V_n)+\tilde{\mu}(U_n\backslash
V_n)<CR^{-2}+\varepsilon/2<\varepsilon.
$$
This completes the proof of Proposition~\ref{mutilde}.
\end{proof}
\begin{remarque}
One may show that for $s\geq 1/2$, the measure $\tilde{\mu}$ {\bf is not}
countably additive on the algebra of the cylindrical set on
$H^s_{rad}(\Theta)$ (see e.g. \cite{DY}).
\end{remarque}
If $u\in H^s_{rad}(\Theta)$ is given
by
$$
u=\sum_{n\geq 1}c_{n}e_{n,s}
$$
then we can consider a Littlewood-Paley decomposition of $u$ defined by
$$
u=\sum_{N-{\rm dyadic}}\Delta_{N}(u),
$$
where $N$ is running over the set of dyadic integers, i.e. the nonnegative
powers of $2$, and, the projector $\Delta_{N}$ is defined by 
$$
\Delta_{N}(u)=\sum_{n\,:\, N\leq \langle z_n\rangle < 2N}c_{n}e_{n,s}\,.
$$
We next state a bound on the $\mu$ measure of functions containing only high
frequencies in their Littlewood-Paley decompositions.
\begin{proposition}\label{high}
Let $\sigma\in [s,1/2[$. 
There exist $C>0$ and $c>0$ such that for every $N_0\geq 1$, every
$\lambda\geq 1$,
$$
\mu
\Big(
u\in H^s_{rad}(\Theta)\, :\, 
\big\|
\sum_{\stackrel{N\geq N_0}{ N-{\rm dyadic }}}
\Delta_{N}(u)\big\|_{H^{\sigma}(\Theta)}\leq\lambda\Big)
\geq 1- Ce^{-c\lambda^{2}N_0^{2(1-\sigma)}}\,.
$$
In particular
\begin{equation}\label{better}
\mu
\big(
u\in H^s_{rad}(\Theta)\, :\, \|u\|_{H^{\sigma}(\Theta)}\leq\lambda
\big)
\geq 1-Ce^{-c\lambda^2}
\end{equation}
and
$$
\mu(H^{\sigma}_{rad}(\Theta))=1\, .
$$
Therefore one can consider $\mu$ as a measure on the space ${\mathcal X}$
defined by (\ref{X}).
\end{proposition}
\begin{proof}[Proof of Proposition~\ref{high}] 
In view of (\ref{cyl-mes}), we obtain that the measure $\mu$ is the
distribution of the random series
\begin{equation}\label{5}
\varphi_{\omega}(r)=
\sum_{n\geq 1}\frac{g_n(\omega)}{z_n^{1-s}}e_{n,s}(r)
=
\sum_{n\geq 1}\frac{g_n(\omega)}{z_n}e_{n}(r)\, ,
\end{equation}
where $g_n(\omega)$ is a sequence of normalised (${\mathcal N}(0,1/\sqrt{2})$) independent identically distributed (i.i.d.)
complex Gaussian random variables, defined in a probability space 
$(\Omega,{\mathcal F},p)$.
More precisely, for $U$ a $\mu$-measurable set, we have
$$
\mu(U)=p(\omega\,:\,\varphi_{\omega}\in U).
$$
Consider a Littlewood-Paley decomposition of (\ref{5})
\begin{equation}\label{LP}
\varphi_{\omega}(r)=\sum_{N-{\rm dyadic}}\Delta_{N}\big(\varphi_{\omega}(r)\big)
\end{equation}
with
$$
\Delta_{N}\big(\varphi_{\omega}(r)\big)=\sum_{n\,:\, N\leq \langle z_n\rangle < 2N}z_{n}^{-1}g_{n}(\omega)e_{n}(r)\,.
$$
We need therefore to establish the bound
$$
p\Big(\omega\in\Omega \,:\,\big\|\sum_{\stackrel{N\geq N_0}{ N-{\rm dyadic }}}\Delta_{N}(\varphi_{\omega})\big\|_{H^{\sigma}(\Theta)}>\lambda\Big)
\leq Ce^{-c\lambda^{2}N_0^{2(1-\sigma)}}\,.
$$
Let us next prove an inequality for Gaussians.
\begin{lemme}\label{lem3}
Let $g_n(\omega)$ be a sequence of normalized i.i.d. complex Gaussian random
variables defined in a probability space 
$(\Omega,{\mathcal F},p)$.
Then there exist positive numbers $c_1,c_2$ such that for every 
finite set of indexes $\Lambda \subset \N $, every $\lambda >0$, 
$$
p\Big(\omega\in \Omega\, : \, \sum_{n\in \Lambda} |g_n(\omega)|^2 >\lambda\Big)\leq
e^{c_1|\Lambda|-c_2\lambda}
\, .
$$
\end{lemme}
\begin{proof}
For every $\zeta>0$,
\begin{equation}\label{thc1}
p\Big(\omega\in \Omega\, : \, \sum_{n\in \Lambda} |g_n(\omega)|^2
>\lambda\Big)=
p\Big(\omega\in \Omega\, : \, \prod_{n\in \Lambda} e^{\zeta|g_n(\omega)|^2}
>e^{\zeta\lambda}\Big)\,.
\end{equation}
For $\zeta<1$, using the Tchebishev inequality and the independence of
$g_n(\omega)$, we deduce that (\ref{thc1}) is bounded by
$$
e^{-\zeta\lambda}\prod_{n\in\Lambda}
\E\big(
e^{\zeta|g_n(\omega)|^2}
\big)=e^{-\zeta\lambda}\, z^{|\Lambda|},
$$
where the positive number $z$ is given by
$$
z=\pi^{-1}
\Big(
\int_{-\infty}^{\infty}
e^{-(1-\zeta)x^2}dx
\Big)^2=\frac{1}{1-\zeta}>1\,.
$$
This completes the proof of Lemma~\ref{lem3}.
\end{proof}
Let us now turn to the proof of Proposition~\ref{high}. 
For $N_0\geq 1$, we set
$$
A_{N_0}=\Big(\omega\in\Omega \,:\,\big\|\sum_{\stackrel{N\geq N_0}{ N-{\rm dyadic }}}\Delta_{N}(\varphi_{\omega})\big\|_{H^{\sigma}(\Theta)}>
\lambda\Big)\,.
$$
Let $\theta$ be a real number such that 
\begin{equation}\label{theta}
0<2\theta<1-2\sigma\,.
\end{equation}
For $\kappa>0$ and $N\geq N_0$, $N$ being a dyadic integer, we set
$$
A_{N,\kappa}=
\Big(\omega\in\Omega
\,:\,\big\|\Delta_{N}(\varphi_{\omega})\big\|_{H^{\sigma}(\Theta)}>\lambda\kappa (N^{-\theta}+(N^{-1}N_0)^{1-\sigma})\Big)\,.
$$
Now, we observe that there exists $\kappa$ sufficiently small depending on $\sigma$
but independent of $N_0$ and $N$ such that
\begin{equation}\label{theta2}
A_{N_0}\subset\,
\bigcup_{\stackrel{N\geq N_0}{ N-{\rm dyadic }}} \,
A_{N,\kappa}\,.
\end{equation}
The restriction
$$
\big\|\Delta_{N}(\varphi_{\omega})\big\|_{H^{\sigma}(\Theta)}>\lambda\kappa
(N^{-\theta}+(N^{-1}N_0)^{1-\sigma})
$$
implies that
$$
\sum_{n\,:\, N\leq \langle z_n\rangle < 2N}z_{n}^{2\sigma}z_{n}^{-2}|g_n(\omega)|^{2}
>
\lambda^2\kappa^2
(N^{-\theta}+(N^{-1}N_0)^{1-\sigma})^2
$$
and therefore, in view of (\ref{zeros}),
$$
\sum_{n\,:\, N\leq \langle z_n\rangle< 2N}
|g_n(\omega)|^{2}
>
c\lambda^2\kappa^2
N^{2-2\sigma}
(N^{-2\theta}+(N^{-1}N_0)^{2-2\sigma}) \,.
$$
Once again invoking (\ref{zeros}), we infer that
$$
\#\{
n\,:\, N\leq \langle z_n\rangle  < 2N
\}\leq CN
$$
and therefore Lemma~\ref{lem3} yields the bound
$$
p(A_{N,\kappa})\leq
e^{c_1N-c_2\lambda^2\kappa^2N_{0}^{2-2\sigma}-c_2\lambda^2\kappa^2N^{2-2\sigma-2\theta}}\,.
$$
The assumption (\ref{theta}) implies that $1<2-2\sigma-2\theta$ and thus
$$
p(A_{N,\kappa})\leq Ce^{-c\lambda^2 N_0^{2-2\sigma}}\, e^{-cN^{2-2\sigma-2\theta}}\,.
$$
Using (\ref{theta2}), a summation over $N$ yields
$$
p(A_{N_0})\leq C e^{-c\lambda^2 N_{0}^{2-2\sigma}}
$$
which completes the proof of Proposition~\ref{high}.
\end{proof}
\begin{remarque}
One can use the method of proof of Proposition~\ref{high} to improve
(\ref{koleda}) to exponential bounds.
\end{remarque}
Let us now turn to the proof of (\ref{parva}) of Theorem~\ref{glavna}.
It is a consequence of the following statement.
\begin{proposition}\label{high-bis}
Let $g_n(\omega)$ be a sequence of normalised i.i.d.
complex Gaussian random variables, defined in a probability space 
$(\Omega,{\mathcal F},p)$.
Let $\chi\, :\, \R\rightarrow\{0,1\}$ be the characteristic function of the 
interval $[0,R]$, $R>0$. 
Define the random series $\varphi_{\omega}$ by
$$
\varphi_{\omega}(r)=\sum_{n\geq 1}\frac{g_n(\omega)}{z_n}e_{n}(r)\, .
$$
Then for every $q>0$,
$$
\E
\Big(
\chi(\|\phi_{\omega}\|_{L^2(\Theta)})\,
\exp(q\int_{\Theta}|V(\phi_{\omega}))|
\Big)
<\infty\,.
$$
\end{proposition} 
\begin{proof}
Thanks to (\ref{rast})
$$
\int_{\Theta}|V(\phi_{\omega})|
\leq C+C\|\phi_{\omega}\|_{L^{\alpha+2}(\Theta)}^{\alpha+2}\,.
$$
Therefore, we have to show that
$$
\E\Big(\chi(\|\phi_{\omega}\|_{L^2(\Theta)})\,\exp(Cq\|\phi_{\omega}\|_{L^{\alpha+2}(\Theta)}^{\alpha+2})\Big)<\infty\,.
$$
Observe that it suffices to show that
$$
\int_{1}^{\infty}f(\lambda)d\lambda<\infty,
$$
where
$$
f(\lambda)=
p\Big(
\omega\in\Omega\,:\,
\|\phi_{\omega}\|_{L^{\alpha+2}(\Theta)}\geq\Big(\frac{\log\lambda}{Cq}\Big)^{\frac{1}{\alpha+2}},\quad
\|\phi_{\omega}\|_{L^{2}(\Theta)}\leq R
\Big).
$$
Set
$$
\gamma:=\Big(\frac{\log\lambda}{Cq}\Big)^{\frac{1}{\alpha+2}}\,.
$$
Let us now fix the real number $\sigma$ according to the Sobolev embedding
restriction
$$
\sigma=2\Big(\frac{1}{2}-\frac{1}{\alpha+2}\Big)\,.
$$
Notice that thanks to (\ref{subcubic}) and (\ref{s}), $1/2>\sigma\geq s$
(of course the important point here is that $\sigma<1/2$).
The Sobolev embedding $H^{\sigma}(\Theta)\subset L^{\alpha+2}(\Theta)$ thus
yields the bound
\begin{equation}\label{sobolev}
\|\phi_{\omega}\|_{L^{\alpha+2}(\Theta)}\leq C_{sob}\|\phi_{\omega}\|_{H^{\sigma}(\Theta)}\,.
\end{equation}
Therefore
$$
f(\lambda)\leq p\Big(\omega\in\Omega\,:\,\|\phi_{\omega}\|_{H^{\sigma}(\Theta)}\geq \gamma/C_{sob},\quad
\|\phi_{\omega}\|_{L^{2}(\Theta)}\leq R \Big).
$$
Consider again the Littlewood-Paley decomposition (\ref{LP}). In the spirit of
the Br\'ezis-Gallouet argument, we set
$$
N_0=\kappa \gamma^{\frac{1}{\sigma}},
$$
where $\kappa>0$ is a small number to be fixed later. Then
$$
\Big(\omega\in\Omega\,:\,\|\phi_{\omega}\|_{H^{\sigma}(\Theta)}\geq \gamma/C_{sob},\quad
\|\phi_{\omega}\|_{L^{2}(\Theta)}\leq R \Big)\subset A_1\cup A_2
$$
with
$$
A_1=
\Big(\omega\in\Omega\,:\,
\Big\|
\sum_{\stackrel{N\leq N_0}{ N-{\rm dyadic }}}\Delta_{N}(\varphi_{\omega})
\Big\|_{H^{\sigma}(\Theta)}\geq \gamma/(4C_{sob}),\quad
\|\phi_{\omega}\|_{L^{2}(\Theta)}\leq R \Big)
$$
and
$$
A_2=
\Big(\omega\in\Omega\,:\,
\Big\|
\sum_{\stackrel{N> N_0}{ N-{\rm dyadic }}}\Delta_{N}(\varphi_{\omega})
\Big\|_{H^{\sigma}(\Theta)}\geq \gamma/(4C_{sob})\Big)\,.
$$
On the other hand
$$
\Big\|
\sum_{\stackrel{N\leq N_0}{ N-{\rm dyadic }}}\Delta_{N}(\varphi_{\omega})
\Big\|_{H^{\sigma}(\Theta)}\leq
CN_{0}^{\sigma}\|\varphi_{\omega}\|_{L^2(\Theta)}\leq CR\kappa^{\sigma}\gamma\,.
$$
Hence for $\kappa\ll 1$, the set $A_1$ is empty. 
This fixes the parameter $\kappa$.
On the other hand, thanks to Proposition~\ref{high},
$$
p(A_2)\leq Ce^{-c\gamma^{2}N_0^{2(1-\sigma)}}\,.
$$
Therefore
$$
f(\lambda)\leq Ce^{-c[\log\lambda]^{\frac{2}{\alpha+2}}N_0^{2(1-\sigma)}}\,.
$$
Coming back to the definitions of $\sigma$ and $N_0$, we get 
$$
[\log\lambda]^{\frac{2}{\alpha+2}}N_0^{2(1-\sigma)}=c[\log\lambda]^{\frac{2}{\alpha}}\,.
$$
The assumption $\alpha<2$ implies $2/\alpha>1$ and therefore $f(\lambda)$ is
integrable on $[1,+\infty[$.
This completes the proof of Proposition~\ref{high-bis}.
\end{proof}
We now state a corollary of  Proposition~\ref{high} and
Proposition~\ref{high-bis}.
\begin{proposition}\label{high-tris}
Let $\sigma\in[s,1/2[$. Then there exist $C>0$ and $c>0$ such that 
$$
\rho\big(u\in H^s_{rad}(\Theta)\, :\, \|u\|_{H^{\sigma}(\Theta)}>\lambda
\big)
\leq Ce^{-c\lambda^2}\,.
$$
\end{proposition}
\begin{proof}
Set
$$
A_{\lambda}=\big(u\in H^s_{rad}(\Theta)\, :\, \|u\|_{H^{\sigma}(\Theta)}>\lambda
\big).
$$
Then using  Proposition~\ref{high} and
Proposition~\ref{high-bis}, we can write
\begin{eqnarray*}
\rho(A_{\lambda})=\int_{A_{\lambda}}d\rho
& = &
\int_{A_{\lambda}}f(u)d\mu(u)
\\
& \leq &
\Big(\int_{A_{\lambda}}f^{2}(u)d\mu(u)\Big)^{1/2}
\Big(\int_{A_{\lambda}}d\mu(u)\Big)^{1/2}
\\
& \leq & C\big(\mu(A_{\lambda})\big)^{\frac{1}{2}}\leq Ce^{-c\lambda^2}\,.
\end{eqnarray*}
This completes the proof of Proposition~\ref{high-tris}.
\end{proof}
Next, we define the finite dimensional sup-spaces of $H^{s}_{rad}(\Theta)$,
$$
E_{N}={\rm span}\,(e_{1,s},\dots,e_{N,s})\,.
$$
We equip $E_N$ with the measures $\mu_N$ and $\rho_N$ which are the natural
restrictions to $E_N$ of $\mu$ and $\rho$ respectively. More precisely for a
Borel set $V\subset \C^N$, we set
\begin{equation}\label{EN}
\tilde{V}=\left\{
c_1e_{1,s}+\dots +c_{N}e_{N,s},\quad (c_1,\dots,c_{N})\in V
\right\}.
\end{equation}
We define the measures $\mu_N$ and $\rho_{N}$ on the sigma algebra of sets of
type (\ref{EN}) by
$$
\mu_{N}(\tilde{V})=\pi^{-N}\Big(\prod_{1\leq n\leq N}z_{n}^{2-2s}\Big)
\int_{V}e^{-{\sum_{1\leq n\leq N} z_n^{2-2s}|c_n|^2}}d^{2}c_1\dots
d^{2}c_{N}
$$
and
\begin{multline*}
\rho_{N}(\tilde{V})=
\pi^{-N}\Big(\prod_{1\leq n\leq N}z_{n}^{2-2s}\Big)\times
\\
\times\int_{V}
f(c_1e_{1,s}+\dots+c_{N}e_{N,s})\,e^{-{\sum_{1\leq n\leq N} z_n^{2-2s}|c_n|^2}}d^{2}c_1\dots
d^{2}c_{N}\,.
\end{multline*}
It is now clear that to every cylindrical set of $H^s_{rad}(\Theta)$ we may
naturally associate a $\mu_N$ and $\rho_{N}$ measurable set on $E_N$, provided $N$ being
sufficiently large. 
For $u\in H^{s}_{rad}(\Theta)$, we define the projector $S_N$,
$$
S_N\,:\, H^{s}_{rad}(\Theta)\longrightarrow E_N
$$
via the formula
\begin{equation}\label{formula}
S_{N}(u)=\sum_{ n= 1}^{N}(u,e_{n,s}) e_{n,s}\, .
\end{equation}
We have the following statement.
\begin{lemme}\label{lim2}
Let $U$ be an open set of $H^{\sigma}_{rad}(\Theta)$, $\sigma\in[s,1/2[$. Then
\begin{equation}\label{jean1}
\rho(U)\leq
\liminf_{N\rightarrow\infty }\rho_{N}(U\cap E_N).
\end{equation}
Moreover, if $F$ is a closed set of $H^{\sigma}_{rad}(\Theta)$, $\sigma\in[s,1/2[$ then
\begin{equation}\label{jean2}
\limsup_{N\rightarrow\infty }\rho_{N}(F\cap E_N)\leq\rho(F)\,.
\end{equation}
\end{lemme}
\begin{proof}
Define the sets
$$
U_{N}\equiv
\big\{u\in H^s_{rad}(\Theta)\,:\, S_{N}(u)\in U\big\}.
$$
Observe that $U\cap E_{N}$ is $\rho_{N}$ measurable and
$$
\rho_{N}(U\cap E_N)=\rho(U_{N}).
$$
We have the inclusion
\begin{equation}\label{exam1}
U\subset \liminf_{N}(U_{N}),
\end{equation}
where
$$
\liminf_{N}(U_{N})\equiv \bigcup_{N\geq 1}\bigcap_{N_1\geq N}U_{N_1}\,.
$$
Indeed, we have that for every $u\in H^{\sigma}_{rad}(\Theta)$,
\begin{equation}\label{exam2}
\lim_{N\rightarrow\infty }\|u-S_{N}(u)\|_{H^{\sigma}(\Theta)}=0\,.
\end{equation}
Therefore, using that $U$ is an open set, we conclude that for every $u\in U$
there exists $N_{0}\geq 1$ such that for $N\geq N_0$ one has $u\in
U_{N}$. Hence we have (\ref{exam1}).
If $A$ is a $\rho$-measurable set, we denote by $\chi_{A}$ the characteristic
function of $A$. Define the non negative functions $f_{N}$ by
$$
f_{N}(u)\equiv\chi_{U_{N}}(u)f(u)\,.
$$
Notice that thanks to (\ref{exam1}),
$$
\liminf_{N\rightarrow\infty}\chi_{U_{N}}\geq\chi_{U}\,.
$$
Next, we set
$$
F(u)\equiv \liminf_{N\rightarrow\infty} f_{N}(u).
$$
Thus
$$
F(u)\geq \chi_{U}(u)f(u)\,.
$$
Since
$$
\rho_{N}(U\cap E_{N})=\int_{H^s_{rad}(\Theta)}f_{N}(u)d\mu(u),
$$
using the Fatou lemma, we get
$$
\liminf_{N\rightarrow\infty}\rho_{N}(U\cap E_{N})\geq
\int_{H^s_{rad}(\Theta)}F(u)d\mu(u)
\geq \int_{U}f(u)d\mu(u)=\rho(U)\,.
$$
Next, we define the sets
$$
F_{N}\equiv
\big\{u\in H^s_{rad}(\Theta)\,:\, S_{N}(u)\in F\big\}.
$$
Thus
$$
\rho_{N}(F\cap E_N)=\rho(F_{N}).
$$
We have that
\begin{equation}\label{exam1bis}
\limsup_{N}(F_{N})\subset F,
\end{equation}
where
$$
\limsup_{N}(F_{N})\equiv \bigcap_{N\geq 1}\bigcup_{N_1\geq N}F_{N_1}\,.
$$
Indeed, suppose that $u\in \limsup_{N}(F_{N})$. Thus there exists a sequence
$(N_k)_{k\in\N}$ tending to infinity such that $u\in F_{N_k}$ which means that for every $k$ one
has $S_{N_k}(u)\in F$. Since $F$ is closed, coming back to (\ref{exam2}), we obtain that $u\in F$ and therefore
we get (\ref{exam1bis}).
If we set
$$
G(u)\equiv \limsup_{N\rightarrow\infty} \tilde{f}_{N}(u),
$$
where
$$
\tilde{f}_{N}(u)\equiv\chi_{F_{N}}(u)f(u)\,.
$$
then we have
$$
G(u)\leq \chi_{F}(u)f(u)
$$
and the Fatou lemma gives,
\begin{eqnarray*}
\limsup_{N\rightarrow\infty}\rho_{N}(F\cap E_{N})
& = &
\limsup_{N\rightarrow\infty}
\int_{H^s_{rad}(\Theta)}
\tilde{f}_{N}(u)d\mu(u)
\\
& \leq &
\int_{H^s_{rad}(\Theta)}G(u)d\mu(u)
\\
& \leq & \int_{F}f(u)d\mu(u)
\\
& = &
\rho(F)\,.
\end{eqnarray*}
This completes the proof of Lemma~\ref{lim2}.
\end{proof}
The next lemma shows that every $\rho$ measurable set can be approximated by
subsets of compact sets in $H^{s}_{rad}(\Theta)$.
\begin{lemme}\label{lim3}
Lets $\sigma\in]s,1/2[$ and denote by $K_{R}$, $R>0$ the ball
$$
K_{R}\equiv \{u\in H^{s}_{rad}(\Theta)\,:\, \|u\|_{H^{\sigma}(\Theta)}\leq R\}.
$$
Then, for every $\rho$ measurable set $A$,
$$
\rho(A)=\lim_{R\rightarrow\infty}\rho(A\cap K_{R})\, .
$$
\end{lemme}
\begin{proof}
Since $\mu(H^{s}_{rad}(\Theta))=\mu(H^{\sigma}_{rad}(\Theta))$ and since
$f(u)\in L^{q}(d\mu(u))$, $1\leq q <\infty$, we obtain that
$\rho(H^{s}_{rad}(\Theta))=\rho(H^{\sigma}_{rad}(\Theta))$.
Therefore, using Proposition~\ref{high-tris}, we can write
\begin{eqnarray*}
0\leq \rho(A)-\rho(A\cap K_{R})& = &\rho(A\cup K_{R})-\rho(K_{R})
\\
& \leq &\rho(H^{s}_{rad}(\Theta))-\rho(K_{R})
\\
& = & \rho(H^{\sigma}_{rad}(\Theta))-\rho(K_{R})
\\
& \leq &
Ce^{-CR^2}
\end{eqnarray*}
which completes the proof of Lemma~\ref{lim3}.
\end{proof}

\section{Bilinear Strichartz estimates}
We now state a localized Strichartz type bilinear estimate associated to the
linear Schr\"odinger group on the unit disc.
\begin{proposition}\label{thm3}
For every $\varepsilon>0$, there exists $\beta<1/2$, there exists $C>0$ such that for
every $N_1,N_2\geq 1$, every $L_1,L_2\geq 1$, every $u_1$, $u_2$ two functions on
$\R\times\Theta$ of the form
$$
u_{j}(t,r)=
\sum_{N_j\leq \langle z_n\rangle < 2N_j}\,c_j(n,t)\, e_{n}(r),\quad j=1,2
$$
where the Fourier transform of $c_j(n,t)$ with respect to $t$ satisfies 
$$
{\rm supp}\, \widehat{c_j}(n,\tau)\subset 
\{
\tau\in\R\,:\, L_{j}\leq \langle\tau+z_n^2\rangle\leq 2L_j
\},\quad j=1,2
$$
one has the bound
$$
\|u_1 u_2\|_{L^2(\R\times \Theta)}\leq C(N_1\wedge
N_2)^{\varepsilon}(L_1L_2)^{\beta}
\|u_1\|_{L^2(\R\times \Theta)}
\|u_2\|_{L^2(\R\times \Theta)}\,.
$$
\end{proposition}
\begin{proof}
Let us first notice that for $j=1,2$,
\begin{equation}\label{l2}
\|u_j\|_{L^2(\R\times \Theta)}^{2}=c\,\sum_{ N_j\leq \langle z_n\rangle < 2N_j}\,\int_{-\infty}^{\infty}|\widehat{c_j}(n,\tau)|^{2}d\tau\,.
\end{equation}
Denote $v(t,r)=u_1(t,r)u_2(t,r)$. Our purpose is thus to estimate
$\|v\|_{L^2(\R\times\Theta)}$.
Equivalently, we need to bound 
$\|\widehat{v}(\tau,\cdot)\|_{L^2(\R_{\tau}\times\Theta)}$.
Write
$$
\|\widehat{v}(\tau,r)\|_{L^2(\R\times\Theta)}^{2}
=
c\int_{0}^{1}\int_{-\infty}^{\infty}
\Big|
\int_{-\infty}^{\infty}\widehat{u_1}(\tau_1,r)\widehat{u_2}(\tau-\tau_1,r)d\tau_1
\Big|^{2}r d\tau dr.
$$
On the other hand $\widehat{u_1}(\tau_1,r)\widehat{u_2}(\tau-\tau_1,r)$ is
equal to
$$
\sum_{N_1\leq \langle z_{n_{1}}\rangle< 2N_1}\, \sum_{ N_2\leq \langle z_{n_{2}}\rangle< 2N_2}\,
\widehat{c_1}(n_1,\tau_1)\widehat{c_2}(n_2,\tau-\tau_1)e_{n_1}(r)e_{n_2}(r)\,.
$$
Therefore, by invoking (\ref{util}), we can write
\begin{multline*}
\Big\|
\int_{-\infty}^{\infty}\widehat{u_1}(\tau_1,r)\widehat{u_2}(\tau-\tau_1,r)d\tau_1
\Big\|_{L^2(\Theta)}
\leq
\int_{-\infty}^{\infty}
\big\|
\widehat{u_1}(\tau_1,r)\widehat{u_2}(\tau-\tau_1,r)
\big\|_{L^2(\Theta)}d\tau_1
\\
\leq \int_{-\infty}^{\infty}
\sum_{ N_1\leq \langle z_{n_{1}}\rangle < 2N_1}\, \sum_{ N_2\leq \langle
  z_{n_{2}}\rangle < 2N_2}\,
|\widehat{c_1}(n_1,\tau_1)||\widehat{c_2}(n_2,\tau-\tau_1)|\|e_{n_1}e_{n_2}\|_{L^2(\Theta)}
d\tau_1
\\
\leq C_{\varepsilon}(N_1\wedge N_2)^{\varepsilon}
\int_{-\infty}^{\infty}
\sum_{ N_1\leq \langle z_{n_{1}}\rangle < 2N_1}\, \sum_{ N_2\leq
  \langle z_{n_{2}}\rangle< 2N_2}\,
|\widehat{c_1}(n_1,\tau_1)||\widehat{c_2}(n_2,\tau-\tau_1)|
d\tau_{1}\,.
\end{multline*}
Our aim is estimate the $L^2(\R_{\tau})$ norm of the last expression. 
For this purpose, we will of course make use of the support properties of
$\widehat{c_j}(n,\tau)$. Using the Cauchy-Schwarz inequality in
$(\tau_1,n_1,n_2)$ gives the bound
\begin{multline*}
\int_{-\infty}^{\infty}
\sum_{ N_1\leq \langle z_{n_{1}}\rangle < 2N_1}\, \sum_{ N_2\leq
 \langle z_{n_{2}}\rangle< 2N_2}\,
|\widehat{c_1}(n_1,\tau_1)||\widehat{c_2}(n_2,\tau-\tau_1)|
d\tau_{1}
\leq
\\
\leq 
\Big(
\int_{-\infty}^{\infty}
\sum_{ N_1\leq \langle z_{n_{1}}\rangle < 2N_1}\, 
\sum_{ N_2\leq\langle  z_{n_{2}}\rangle < 2N_2}\,
|\widehat{c_1}(n_1,\tau_1)|^{2}|\widehat{c_2}(n_2,\tau-\tau_1)|^{2}
d\tau_{1}
\Big)^{\frac{1}{2}}
|\Lambda(\tau,L_1,L_2,N_1,N_2)|^{\frac{1}{2}},
\end{multline*}
where $\Lambda(\tau,L_1,L_2,N_1,N_2)$ is the following set of $\R\times\N\times\N$,
\begin{multline}\label{Lambda1}
\Lambda(\tau,L_1,L_2,N_1,N_2)=
\{
(\tau_1,n_1,n_2)\in\R\times\N\times\N\,:\, 
L_1\leq \langle\tau_1+z_{n_1}^2\rangle\leq 2L_1,
\\
L_2\leq \langle\tau-\tau_1+z_{n_2}^2\rangle\leq 2L_2,\,\,
\langle z_{n_1}\rangle\in [N_1,2N_1], \, \langle z_{n_2}\rangle \in[N_2,2N_2]
\}.
\end{multline}
The next lemma contains the main combinatorial ingredient of our analysis.
\begin{lemme}\label{combonatorics}
For every $\varepsilon>0$ there exists $C_{\varepsilon}>0$ such that
for every $\tau\in\R$, every $N_1,N_2\geq 1$, every $L_1,L_2\geq 1$,
$$
|\Lambda(\tau,L_1,L_2,N_1,N_2)|\leq C_{\varepsilon}(N_1\wedge
N_2)^{\varepsilon}(L_1 L_2).
$$
\end{lemme}
\begin{proof}
Consider the set $\widetilde{\Lambda}(\tau,L_1,L_2,N_1,N_2)$ of $\N\times\N$, defined by
\begin{multline*}
\widetilde{\Lambda}(\tau,L_1,L_2,N_1,N_2)=
\{
(n_1,n_2)\in\N\times\N\,:\, 
\langle\tau+z_{n_1}^2 +z_{n_2}^2\rangle\leq 2(L_1+L_2),
\\
\langle z_{n_1}\rangle\in [N_1,2N_1], \, \langle z_{n_2}\rangle \in[N_2,2N_2]
\}.
\end{multline*}
Let $(\tau_1,n_1,n_2)\in \Lambda(\tau,L_1,L_2,N_1,N_2)$. Then the triangle
inequality yields
$$
\langle\tau_1+z_{n_1}^2+z_{n_2}^2\rangle
\leq
\langle\tau_1+z_{n_1}^2\rangle+\langle\tau-\tau_1+z_{n_2}^2\rangle\leq
2(L_1+L_2).
$$
Therefore $(\tau_1,n_1,n_2)\in \Lambda(\tau,L_1,L_2,N_1,N_2)$ implies that
$(n_1,n_2)\in \widetilde{\Lambda}(\tau,L_1,L_2,N_1,N_2)$. On the other hand
for a fixed $(n_1,n_2)\in \widetilde{\Lambda}(\tau,L_1,L_2,N_1,N_2)$ the
Lebesgue measure of the possible $\tau_1$ such that $(\tau_1,n_1,n_2)\in
\Lambda(\tau,L_1,L_2,N_1,N_2)$ is bounded by $C(L_1\wedge L_2)$. Therefore
\begin{equation}\label{back}
|\Lambda(\tau,L_1,L_2,N_1,N_2)|\leq C(L_1\wedge
L_2)|\widetilde{\Lambda}(\tau,L_1,L_2,N_1,N_2)|.
\end{equation}
We next estimate $|\widetilde{\Lambda}(\tau,L_1,L_2,N_1,N_2)|$.
Observe that $z_{n_1}^2 +z_{n_2}^2$ ranges in an interval of size $\leq
C(L_1+L_2)$. Hence, thanks to (\ref{zeros}) the expression 
$(4n_1-1)^{2}+(4n_2-1)^{2}$ also ranges in
an interval of size $C(L_1+L_2)$, where the constant $C$ is independent of
$\tau$, $L_1,L_2$, $N_1,N_2$. 
Indeed, suppose that for some $A\in\R$,
\begin{equation}\label{irina1}
A\leq z_{n_1}^2 +z_{n_2}^2\leq A+C(L_1+L_2).
\end{equation}
In (\ref{irina1}), $A$ is the parameter we have no control on.
Using (\ref{zeros}), we obtain that (\ref{irina1}) implies
\begin{equation}\label{irina2}
\frac{16A}{\pi^2}\leq
(4n_1-1)^{2}+(4n_2-1)^{2}+R(n_1,n_2)\leq
\frac{16A+16C(L_1+L_2)}{\pi^2}\,,
\end{equation}
where, thanks to the estimate on the remainder in (\ref{zeros}), the function
$R(n_1,n_2)$ in (\ref{irina2}) satisfies
$$
|R(n_1,n_2)|\leq \widetilde{C}\,.
$$
Therefore, (\ref{irina2}) implies
$$
\frac{16A}{\pi^2}-\widetilde{C}\leq
(4n_1-1)^{2}+(4n_2-1)^{2}\leq
\frac{16A+16C(L_1+L_2)}{\pi^2}
+\widetilde{C}\,.
$$
Thus $(4n_1-1)^{2}+(4n_2-1)^{2}$ ranges in an interval of size
$$
\frac{16C(L_1+L_2)}{\pi^2}+2\widetilde{C}
\leq \Big(\frac{16C}{\pi^2}+\widetilde{C}\Big)(L_1+L_2)
$$
exactly as we claimed.
Denote the interval where $(4n_1-1)^{2}+(4n_2-1)^{2}$ can range by $\Delta$.
Another appeal to (\ref{zeros}) yields that the restrictions $\langle z_{n_1}\rangle\in
[N_1,2N_1]$
and $\langle z_{n_2}\rangle \in[N_2,2N_2]$ imply 
the bounds
$$
0\leq n_1\leq cN_1,\quad 
0\leq n_2\leq cN_2 \,.
$$
Let $l$ be and integer in the interval $\Delta$. 
Then we have the bound
\begin{multline}\label{malko}
\#\Big((n_1,n_2)\in \N\times\N\,:\, l= (4n_1-1)^{2}+(4n_2-1)^{2}, 
\\
0\leq n_1\leq cN_1,\quad 0\leq n_2\leq cN_2
\Big)
\leq C_{\varepsilon}(N_1\wedge N_2)^{\varepsilon}
\end{multline}
Indeed, if $l\leq 2c^{2}(N_1\wedge N_2)^{2005}$ then 
the left hand-side of (\ref{malko}) is bounded by $C_{\varepsilon}(\min(N_1,N_2))^{\varepsilon}$
by the standard bound (see e.g. \cite{Gr,Nat}) on the number of divisors in the ring of Gaussian
integers $\Z[i]$.
Let us next suppose that 
$$
l\geq 2c^{2}(N_1\wedge N_2)^{2005}+1.
$$
By symmetry, we can suppose that $N_2\geq N_1$. We have that $4n_2-1\in I$, where
the interval $I$ is defined by 
$$
I=\big[\sqrt{l-(4cN_1-1)^{2}},\sqrt{l}\big].
$$
But the size of $I$ is bounded by
$$
\frac{(4cN_{1}-1)^{2}}{\sqrt{l}}
\leq
\frac{CN_{1}^{2}}{\sqrt{c^{2}N_{1}^{2005}+1}}\leq C.
$$
Therefore the size of $I$ is bounded by a constant uniform in $N_1$, $N_2$ and
$l$.
Thus in the case $l\geq 2c_{3}^{2}(N_1\wedge N_2)^{2005}+1$, we can get even
better then (\ref{malko}), namely we have a bound by a uniform constant.
This completes the proof of (\ref{malko}).
Using (\ref{malko}) we infer that
$$
|\widetilde{\Lambda}(\tau,L_1,L_2,N_1,N_2)|\leq
C_{\varepsilon}|\Delta|(N_1\wedge N_2)^{\varepsilon}=
C_{\varepsilon}(L_1+L_2)(N_1\wedge N_2)^{\varepsilon}\,.
$$
Observe that $(L_1\wedge L_2)(L_1+L_2)\leq 2(L_1L_2)$. Therefore coming back
to (\ref{back}) completes the proof of Lemma~\ref{combonatorics}.
\end{proof}
Thanks to Lemma~\ref{combonatorics}, we may write
\begin{multline*}
\Big\|
\int_{-\infty}^{\infty}\widehat{u_1}(\tau_1,r)\widehat{u_2}(\tau-\tau_1,r)d\tau_1
\Big\|_{L^2(\Theta)}\leq
C_{\varepsilon}(N_1\wedge
N_2)^{\varepsilon}(L_1 L_2)^{\frac{1}{2}}
\\
\Big(
\int_{-\infty}^{\infty}
\sum_{ N_1\leq \langle z_{n_{1}}\rangle< 2N_1}\, 
\sum_{ N_2\ \leq\langle z_{n_{2}}\rangle < 2N_2}\,
|\widehat{c_1}(n_1,\tau_1)|^{2}|\widehat{c_2}(n_2,\tau-\tau_1)|^{2}
d\tau_{1}
\Big)^{\frac{1}{2}}\,.
\end{multline*}
Squaring the above inequality, integration over $\tau\in\R$ and using (\ref{l2}) gives the bound
\begin{equation}\label{glaven}
\|u_1 u_2\|_{L^2(\R\times \Theta)}\leq C_{\varepsilon}(N_1\wedge
N_2)^{\varepsilon}(L_1 L_2)^{\frac{1}{2}}
\|u_1\|_{L^2(\R\times \Theta)}
\|u_2\|_{L^2(\R\times \Theta)}\,.
\end{equation}
We however claimed that the power of $L_1L_2$ can be smaller than $1/2$. In
order to obtain this better bound with respect to the $L_1,L_2$ localization,
we will get an inequality which is better than (\ref{glaven}) as far as 
$(L_1 L_2)$ is concerned but which is very weak with respect to the $N_1,N_2$
localization. 
Using the formula for the inverse of the Fourier transform, the support
properties of the Fourier transform of $c_{j}(n,t)$, and the Cauchy-Schwarz
inequality, we obtain the bound
$$
|c_{j}(n,t)|^{2}\leq
CL_{j}\int_{-\infty}^{\infty}|\widehat{c_j}(n,\tau)|^{2}d\tau,\quad j=1,2.
$$
Hence, we infer that
\begin{eqnarray*}
\|u_{j}(t,\cdot)\|_{L^2(\Theta)}^{2} & = &
\sum_{N_j\leq \langle z_n\rangle < 2N_j}|c_j(n,t)|^{2}
\\
& \leq &
CL_{j}\sum_{N_j\leq \langle z_n\rangle < 2N_j}
\int_{-\infty}^{\infty}|\widehat{c_j}(n,\tau)|^{2}d\tau
\\
& = &
CL_{j}\|u_{j}\|_{L^2(\R\times\Theta)}^{2}.
\end{eqnarray*}
Therefore
\begin{equation}\label{m1}
\|u_{j}\|_{L^{\infty}(\R;L^2(\Theta))}
\leq CL_{j}^{\frac{1}{2}}\|u_{j}\|_{L^2(\R\times\Theta)}.
\end{equation}
Interpolation (it is in fact simply the H\"older inequality) with the equality
$$
\|u_{j}\|_{L^{2}(\R;L^2(\Theta))}=\|u_{j}\|_{L^2(\R\times\Theta)}
$$
gives the bound
\begin{equation}\label{m2}
\|u_{j}\|_{L^{4}(\R;L^2(\Theta))}
\leq CL_{j}^{\frac{1}{4}}\|u_{j}\|_{L^2(\R\times\Theta)},\quad j=1,2.
\end{equation}
Recall that (\ref{zeros}) implies that 
$$
\#(n\in\N\,:\, N\leq \langle z_n\rangle< 2N)\leq CN.
$$
Therefore, using (\ref{eigen}) and the Cauchy-Schwarz inequality, we get
\begin{eqnarray*}
|u_{j}(t,r)|
& \leq  & CN_{j}^{\frac{1}{2}}
\sum_{N_j\leq \langle z_n\rangle < 2N_j}|c_{j}(n,t)|
\\
&\leq &
CN_{j}^{\frac{1}{2}}N_{j}^{\frac{1}{2}}
\Big(\sum_{N_j\leq \langle z_n\rangle < 2N_j}|c_{j}(n,t)|^{2}\Big)^{\frac{1}{2}}
\\
& \leq &
CN_{j}L_{j}^{\frac{1}{2}}\|u_{j}\|_{L^2(\R\times\Theta)}\,.
\end{eqnarray*}
Thus
\begin{equation}\label{m3}
\|u_{j}\|_{L^{\infty}(\R;L^{\infty}(\Theta))}
\leq CN_{j}L_{j}^{\frac{1}{2}}\|u_{j}\|_{L^2(\R\times\Theta)}.
\end{equation}
Next, we can write 
\begin{eqnarray*}
\|u_{j}(t,\cdot)\|_{L^{\infty}(\Theta)}^{2}
& \leq  & CN_{j}
\Big(\sum_{N_j\leq \langle z_n\rangle < 2N_j}|c_{j}(n,t)|\Big)^{2}
\\
&\leq &
CN_{j}^{2}\sum_{N_j\leq \langle z_n\rangle < 2N_j}|c_{j}(n,t)|^{2}
\\
& = &
CN_{j}^{2}\|u_{j}(t,\cdot)\|_{L^2(\Theta)}^{2}.
\end{eqnarray*}
Integration of the last inequality over $t\in\R$ gives 
\begin{equation}\label{m4}
\|u_{j}\|_{L^{2}(\R;L^{\infty}(\Theta))}
\leq CN_{j}\|u_{j}\|_{L^2(\R\times\Theta)}.
\end{equation}
Interpolation between (\ref{m3}) and (\ref{m4}) now gives
\begin{equation}\label{m5}
\|u_{j}\|_{L^{4}(\R;L^{\infty}(\Theta))}
\leq CL_{j}^{\frac{1}{4}}N_{j}\|u_{j}\|_{L^2(\R\times\Theta)}.
\end{equation}
Suppose that $N_{1}\leq N_2$. Then
using (\ref{m2}), (\ref{m5}) and the H\"older inequality, we obtain
\begin{eqnarray*}
\|u_1 u_2\|_{L^2(\R\times\Theta)} & \leq & 
\|u_{1}\|_{L^{4}(\R;L^{\infty}(\Theta))}
\|u_{2}\|_{L^{4}(\R;L^{2}(\Theta))}
\\
& \leq &
C(L_1L_2)^{\frac{1}{4}}N_{1}
\|u_{1}\|_{L^2(\R\times\Theta)}\|u_{2}\|_{L^2(\R\times\Theta)}\,.
\end{eqnarray*}
Therefore, we arrive at
\begin{equation}\label{neglaven}
\|u_1 u_2\|_{L^2(\R\times\Theta)}\leq C(L_1L_2)^{\frac{1}{4}}(N_1\wedge N_2)
\|u_{1}\|_{L^2(\R\times\Theta)}\|u_{2}\|_{L^2(\R\times\Theta)}\,.
\end{equation}
Interpolation between (\ref{glaven}) and (\ref{neglaven}) completes the proof of Proposition~\ref{thm3}.
\end{proof}
We will also need the following variant of Proposition~\ref{thm3}.
\begin{proposition}\label{thm3triis}
For every $\varepsilon>0$, there exists $\beta<1/2$, there exists $C>0$ such that for
every $N_1,N_2\geq 1$, every $L_1,L_2\geq 1$, every $u_1$, $u_2$ two functions on
$\R\times\Theta$ of the form
$$
u_{1}(t,r)=
\sum_{N_1\leq \langle z_n\rangle < 2N_1}\,c_1(n,t)\, e_{n}(r)
$$
and
$$
u_{2}(t,r)=
\sum_{N_2\leq \langle z_n\rangle < 2N_2}\,c_2(n,t)\, e'_{n}(r)
$$
where the Fourier transform of $c_j(n,t)$ with respect to $t$ satisfies 
$$
{\rm supp}\, \widehat{c_j}(n,\tau)\subset \{\tau\in\R\,:\, L_{j}\leq \langle\tau+z_n^2\rangle\leq 2L_j\},\quad j=1,2
$$
one has the bound
$$
\|u_1 u_2\|_{L^2(\R\times \Theta)}\leq C(N_1\wedge N_2)^{\varepsilon}(L_1L_2)^{\beta}
\|u_1\|_{L^2(\R\times \Theta)}\|u_2\|_{L^2(\R\times \Theta)}\,.
$$
\end{proposition}
\begin{proof}
Recall that the function $e_{n}$ satisfies the equation
$$
re''_{n}(r)+e'_{n}(r)=-z_{n}^{2}r\,e_{n}(r).
$$
Therefore, using that for $m\neq n$, $e_{m}$ and $e_{n}$ are orthogonal in
$L^2(\Theta)$ and vanishing at $r=1$, an integration by parts gives for $m\neq n$,
\begin{eqnarray*}
\int_{0}^{1}e'_{m}(r)e'_{n}(r)rdr & = &
-\int_{0}^{1}e_{m}(r)(e'_{n}(r)r)'dr+
\Big[e_{m}(r)e'_{n}(r)r\Big]^{1}_{0}
\\
& = &
-\int_{0}^{1}e_{m}(r)(e'_{n}(r)+re''_{n}(r))dr
\\
& = &
z_{n}^{2}\int_{0}^{1}e_{m}(r)e_{n}(r)rdr=0\,.
\end{eqnarray*}
Therefore, we obtain that
\begin{equation}\label{l2bis}
\|u_2\|_{L^2(\R\times \Theta)}^{2}=c\,
\sum_{ N_2\leq \langle z_n\rangle < 2N_2}\,
\|e'_{n}\|_{L^2(\Theta)}^{2}
\int_{-\infty}^{\infty}|\widehat{c_2}(n,\tau)|^{2}d\tau
\end{equation}
and from now on the proof of Proposition~\ref{thm3triis} follows the lines of the proof of Proposition~\ref{thm3}.
Indeed, using (\ref{utilbis}), we can write
\begin{multline*}
\Big\|
\int_{-\infty}^{\infty}\widehat{u_1}(\tau_1,r)\widehat{u_2}(\tau-\tau_1,r)d\tau_1
\Big\|_{L^2(\Theta)}
\leq
\int_{-\infty}^{\infty}
\big\|
\widehat{u_1}(\tau_1,r)\widehat{u_2}(\tau-\tau_1,r)
\big\|_{L^2(\Theta)}d\tau_1
\\
\leq \int_{-\infty}^{\infty}
\sum_{ N_1\leq \langle z_{n_{1}}\rangle < 2N_1}\, \sum_{ N_2\leq \langle
  z_{n_{2}}\rangle < 2N_2}\,
|\widehat{c_1}(n_1,\tau_1)||\widehat{c_2}(n_2,\tau-\tau_1)|\|e_{n_1}e'_{n_2}\|_{L^2(\Theta)}
d\tau_1
\\
\leq C_{\varepsilon}(N_1\wedge N_2)^{\varepsilon}
\int_{-\infty}^{\infty}
\sum_{ N_1\leq \langle z_{n_{1}}\rangle < 2N_1}\, \sum_{ N_2\leq
  \langle z_{n_{2}}\rangle< 2N_2}\,
|\widehat{c_1}(n_1,\tau_1)|
\|e'_{n_2}\|_{L^2(\Theta)}
|\widehat{c_2}(n_2,\tau-\tau_1)|
d\tau_{1}\,.
\end{multline*}
Again,
our goal is to estimate the $L^2(\R_{\tau})$ norm of the last expression. 
Using the Cauchy-Schwarz inequality in $(\tau_1,n_1,n_2)$ yields
\begin{multline*}
\int_{-\infty}^{\infty}
\sum_{ N_1\leq \langle z_{n_{1}}\rangle < 2N_1}\, \sum_{ N_2\leq
 \langle z_{n_{2}}\rangle< 2N_2}\,
|\widehat{c_1}(n_1,\tau_1)|\,\|e'_{n_2}\|_{L^2(\Theta)}\,|\widehat{c_2}(n_2,\tau-\tau_1)|
d\tau_{1}
\leq
\\
\leq 
\Big(
\int_{-\infty}^{\infty}
\sum_{ N_1\leq \langle z_{n_{1}}\rangle < 2N_1}\, 
\sum_{ N_2\leq\langle  z_{n_{2}}\rangle < 2N_2}\,
|\widehat{c_1}(n_1,\tau_1)|^{2}
\|e'_{n_2}\|_{L^2(\Theta)}^{2}
|\widehat{c_2}(n_2,\tau-\tau_1)|^{2}
d\tau_{1}
\Big)^{\frac{1}{2}}
\\
\times
|\Lambda(\tau,L_1,L_2,N_1,N_2)|^{\frac{1}{2}},
\end{multline*}
where $\Lambda(\tau,L_1,L_2,N_1,N_2)$ is defined by (\ref{Lambda1}).
A use of Lemma~\ref{combonatorics} now gives
\begin{multline*}
\Big\|
\int_{-\infty}^{\infty}\widehat{u_1}(\tau_1,r)\widehat{u_2}(\tau-\tau_1,r)d\tau_1
\Big\|_{L^2(\Theta)}\leq
C_{\varepsilon}(N_1\wedge
N_2)^{\varepsilon}(L_1 L_2)^{\frac{1}{2}}
\\
\Big(
\int_{-\infty}^{\infty}
\sum_{ N_1\leq \langle z_{n_{1}}\rangle< 2N_1}\, 
\sum_{ N_2\ \leq\langle z_{n_{2}}\rangle < 2N_2}\,
|\widehat{c_1}(n_1,\tau_1)|^{2}|
\|e'_{n_2}\|_{L^2(\Theta)}^{2}
\widehat{c_2}(n_2,\tau-\tau_1)|^{2}
d\tau_{1}
\Big)^{\frac{1}{2}}\,.
\end{multline*}
and therefore
\begin{equation*}
\|u_1 u_2\|_{L^2(\R\times \Theta)}\leq C_{\varepsilon}(N_1\wedge
N_2)^{\varepsilon}(L_1 L_2)^{\frac{1}{2}}
\|u_1\|_{L^2(\R\times \Theta)}
\|u_2\|_{L^2(\R\times \Theta)}\,.
\end{equation*}
Next, using the localisation of the Fourier transforms of $c_{1}(n,t)$, as in
the proof of Proposition~\ref{thm3}, we get the bound
\begin{equation*}
\|u_{1}\|_{L^{4}(\R;L^2(\Theta))}
\leq CL_{1}^{\frac{1}{4}}\|u_{1}\|_{L^2(\R\times\Theta)}.
\end{equation*}
Next, we estimate $u_2$ as follows
\begin{eqnarray*}
\|u_{2}(t,\cdot)\|_{L^2(\Theta)}^{2} & = & 
\sum_{ N_2\leq\langle  z_{n}\rangle < 2N_2}
|c_{2}(n,t)|^{2}\|e'_{n}\|_{L^2(\Theta)}^{2}
\\
&\leq &
C \sum_{ N_2\leq\langle  z_{n}\rangle < 2N_2} 
\Big(
L_{2}\int_{-\infty}^{\infty}|\widehat{c_2}(n,\tau)|^{2}d\tau
\Big)\|e'_{n}\|_{L^2(\Theta)}^{2}
\\
& = &
CL_{2}\|u_2\|_{L^2(\R\times\Theta)}^{2}\,.
\end{eqnarray*}
Therefore
\begin{equation*}
\|u_{2}\|_{L^{\infty}(\R;L^2(\Theta))}
\leq CL_{2}^{\frac{1}{2}}\|u_{2}\|_{L^2(\R\times\Theta)}.
\end{equation*}
Interpolating with the equality
\begin{equation*}
\|u_{2}\|_{L^{2}(\R;L^2(\Theta))}
= \|u_{2}\|_{L^2(\R\times\Theta)}.
\end{equation*}
gives
\begin{equation*}
\|u_{2}\|_{L^{4}(\R;L^2(\Theta))}
\leq CL_{2}^{\frac{1}{4}}\|u_{2}\|_{L^2(\R\times\Theta)}.
\end{equation*}
Thus
\begin{equation*}
\|u_{j}\|_{L^{4}(\R;L^2(\Theta))}
\leq CL_{j}^{\frac{1}{4}}\|u_{j}\|_{L^2(\R\times\Theta)},\quad j=1,2.
\end{equation*}
Next, using (\ref{eigen}), we get the bound 
\begin{equation}\label{irinai}
\|u_{j}\|_{L^{\infty}(\R;L^{\infty}(\Theta))}
\leq CN_{j}
L_{j}^{\frac{1}{2}}\|u_{j}\|_{L^2(\R\times\Theta)},\quad j=1,2.
\end{equation}
Indeed, for $j=1$ such an inequality is already proved in Proposition~\ref{thm3}.
For $j=2$, we can write by using (\ref{eigen}),
\begin{eqnarray*}
|u_{2}(t,r)|& \leq & CN_{2}^{\frac{1}{2}}\sum_{ N_2\leq\langle  z_{n}\rangle < 2N_2}|c_{2}(n,t)|\|e'_{n}\|_{L^2(\Theta)}
\\
& \leq & CN_{2}\Big( \sum_{ N_2\leq\langle  z_{n}\rangle < 2N_2}|c_{2}(n,t)|^{2}\|e'_{n}\|_{L^2(\Theta)}^{2}\Big)
\\
& \leq & CN_{2}L_{2}^{\frac{1}{2}}\|u_2\|_{L^2(\R\times\Theta)}
\end{eqnarray*}
and thus (\ref{irinai}) for $j=2$.
Moreover,
\begin{equation}\label{irinaii}
\|u_{j}\|_{L^{2}(\R;L^{\infty}(\Theta))}
\leq CN_{j}\|u_{j}\|_{L^2(\R\times\Theta)},\quad j=1,2.
\end{equation}
Indeed, for $j=1$ it is already proved in Proposition~\ref{thm3}.
For $j=2$, by invoking once again (\ref{eigen}), and the Cauchy-Schwartz
inequality, we obtain
\begin{eqnarray*}
\|u_{2}(t,\cdot)\|_{L^{\infty}(\Theta)}^{2} & \leq &
CN_2
\Big(
\sum_{ N_2\leq\langle  z_{n}\rangle < 2N_2}|c_{2}(n,t)|\|e'_{n}\|_{L^2(\Theta)}
\Big)^{2}
\\
& \leq &
CN_{2}^{2}
\sum_{ N_2\leq\langle  z_{n}\rangle < 2N_2}|c_{2}(n,t)|^{2}\|e'_{n}\|_{L^2(\Theta)}^{2}\,.
\end{eqnarray*}
Integration of the last inequality over $t$ gives (\ref{irinaii}) for $j=2$.
An interpolation gives
\begin{equation*}
\|u_{j}\|_{L^{4}(\R;L^{\infty}(\Theta))}
\leq
CL_{j}^{\frac{1}{4}}N_{j}\|u_{j}\|_{L^2(\R\times\Theta)},\quad j=1,2.
\end{equation*}
Then the H\"older inequality gives
\begin{equation*}
\|u_1 u_2\|_{L^2(\R\times\Theta)}\leq C(L_1L_2)^{\frac{1}{4}}(N_1\wedge
N_2)
\|u_{1}\|_{L^2(\R\times\Theta)}\|u_{2}\|_{L^2(\R\times\Theta)}\,.
\end{equation*}
A final interpolation completes the proof of Proposition~\ref{thm3triis}.
\end{proof}
\section{Bourgain spaces}
We denote by $L^2_{rad}(\Theta)$ the $L^2$ radial functions on the unit disc.
We endow $L^2_{rad}(\Theta)$ with the natural Hilbert space structure.
Similarly, we denote by $L^2_{rad}(\R\times \Theta)$ the $L^2$ functions on $\R\times
\Theta$, radial with respect to the second argument.
For $\sigma<1/2$, the norm in $H^{\sigma}_{rad}(\Theta)$ of a radial
function
$$
v=\sum_{n\geq 1}c_n e_n
$$
can be expressed as
$$
\|v\|_{H^{\sigma}_{rad}(\Theta)}^{2}=
\sum_{n\geq 1}z_n^{2\sigma}\, |c_n|^{2}\, .
$$
In this paper, we will only consider spaces of Sobolev regularity $<1/2$ and
thus there is no need to specify the boundary conditions on $\partial\Theta$
(in our context it simply means $r=1$).
More precisely the choice of the Dirichlet eigenfunctions $e_n$ as basis of
$L^2_{rad}(\Theta)$ is not of importance in the definition of
$H^{\sigma}_{rad}(\Theta)$, $\sigma<1/2$.
\\

Next, we define the Bourgain spaces $X^{\sigma,b}_{rad}(\R\times\Theta)$ of
functions on $\R\times\Theta$ which are radial with respect to the second
argument. These spaces are equipped with the norm
$$
\|u\|_{X^{\sigma,b}_{rad}(\R\times\Theta)}^{2}=\sum_{n\geq 1}
z_{n}^{2\sigma}
\|\langle
\tau+z_n^2\rangle^{b}\widehat{c_n}(\tau)\|_{L^{2}(\R_{\tau})}^{2}
\,,
$$
where
$$
u(t)=\sum_{n\geq 1}c_n(t) e_n\,.
$$
Notice that 
\begin{equation}\label{libre}
\|u\|_{X^{\sigma,b}_{rad}(\R\times\Theta)}=\|\exp(-it\Delta)(u(t))\|_{H^{b}(\R;H^{\sigma}_{rad}(\Theta))}\,.
\end{equation}
Indeed, using that
$$
\exp(-it\Delta)(u(t))=\sum_{n\geq 1}\exp(itz_{n}^{2})c_{n}(t)e_{n}
$$
and since
$$
\widehat{\exp(itz_{n}^{2})c_{n}(t)}(\tau)=\widehat{c_{n}}(\tau-z_{n}^{2}),
$$
we arrive at
\begin{eqnarray*}
\|\exp(-it\Delta)(u(t))\|_{H^{b}(\R;H^{\sigma}_{rad}(\Theta))}^{2} & = &
\sum_{n\geq 1}z_{n}^{2\sigma}\|\exp(itz_{n}^{2})c_{n}(t)\|_{H^b(\R)}^{2}
\\
& = &
\sum_{n\geq 1}z_{n}^{2\sigma}\|\langle\tau\rangle^{b}\widehat{c_n}(\tau-z_n^2)\|_{L^2(\R)}^{2}
\\
& = &
\|u\|_{X^{\sigma,b}_{rad}(\R\times\Theta)}^{2}\,.
\end{eqnarray*}
This proves (\ref{libre}).
Clearly $X^{\sigma,b}_{rad}(\R\times\Theta)$ have a Hilbert space structure
and for $0\leq \sigma<1/2$ we
can see $X^{-\sigma,-b}_{rad}(\R\times\Theta)$ as its dual via the
$L^2(\R\times\Theta)$ pairing.
A one dimensional Sobolev embedding (for functions with values in
$H^{\sigma}_{rad}(\Theta)$) yields the estimate
\begin{equation}\label{Sobolev}
\|u\|_{L^{\infty}(\R\,;\,H^{\sigma}_{rad}(\Theta))}\leq C_{b}\|u\|_{X^{\sigma,b}_{rad}(\R\times \Theta)},\quad b>\frac{1}{2}.
\end{equation}
Next, for $T>0$, we define the restriction spaces
$X^{\sigma,b}_{rad}([-T,T]\times\Theta)$, equipped with the natural norm
$$
\|u\|_{X^{\sigma,b}_{rad}([-T,T]\times\Theta)}=
\inf\{
\|w\|_{X^{\sigma,b}_{rad}(\R\times\Theta)},\quad w\in
X^{\sigma,b}_{rad}(\R\times\Theta)\quad {\rm with}\quad w|_{]-T,T[}=u 
\}.
$$
Therefore (\ref{Sobolev}) yields
\begin{equation*}
\|u\|_{L^{\infty}([-T,T]\,;\,H^{\sigma}_{rad}(\Theta))}\leq C_{b}\|u\|_{X^{\sigma,b}_{rad}([-T,T]\times \Theta)},\quad b>\frac{1}{2}
\end{equation*}
which implies that for $b>1/2$ the space $X^{\sigma,b}_{rad}([-T,T]\times \Theta)$ 
is continuously embedded in
$C([-T,T]\,;\,H^{\sigma}_{rad}(\Theta))$.
Similarly, for $I\subset\R$ an interval, we can define the the restriction spaces
$X^{\sigma,b}_{rad}(I\times\Theta)$, equipped with the natural norm.
\\

Following \cite{BGTens}, our next purpose is to express the norm in
$X^{\sigma,b}_{rad}(\R\times\Theta)$ in terms of some basic localisation
operators. 
Recall that for $u=\sum_{n\geq 1}c_n e_n$, the projector $\Delta_{N}$ is defined by
$$
\Delta_{N}(u)=\sum_{n\,:\, N\leq\langle z_n\rangle< 2N}\, c_n e_n\,.
$$
For $N\geq 2$ a dyadic integer, we define the projector $\tilde{S}_N$ by
$$
\tilde{S}_{N}= 
\sum_{\stackrel{N_1\leq N/2}{ N_1-{\rm dyadic }}}\Delta_{N_1}\, .
$$
For a notational convenience, we assume that $\tilde{S}_{1}$ is zero.
Notice that $\tilde{S}_{N}$ is essentially equivalent to $S_{N}$, where the projector $S_{N}$ is defined in (\ref{formula}).
For $N,L$ positive integers, we define $\Delta_{N,L}$ by
\begin{equation}\label{jpc}
\Delta_{N,L}(u)=\frac{1}{2\pi}
\sum_{n\,:\, N\leq\langle z_n\rangle< 2N}
\Big(\int_{L\leq \langle\tau+z_n^2\rangle\leq 2L}
\widehat{c_n}(\tau)e^{it\tau}d\tau\Big)e_{n},
\end{equation}
where 
$$
u(t)=\sum_{n\geq 1}c_n(t)e_n\,.
$$
Then for $u\in X^{\sigma,b}_{rad}(\R\times\Theta)$ (with the natural
interpretation of the $\tau$ integration in (\ref{jpc}) if $b<0$), we can
write the identity
$$
u=\sum_{L,N-{\rm dyadic }}\Delta_{N,L}(u)
$$
in $X^{\sigma,b}_{rad}(\R\times\Theta)$.
Next, we have that there exists a constant $C_{\sigma,b}>1$ which depends
continuously on $\sigma$ and $b$ such that
\begin{eqnarray*}
C_{\sigma,b}^{-1}
L^{b}N^{\sigma} \|\Delta_{N,L}(u)\|_{L^2(\R\times \Theta)}
& \leq &
\|\Delta_{N,L}(u)\|_{X^{\sigma,b}_{rad}(\R\times \Theta)}
\\
& \leq &
C_{\sigma,b}
L^{b}N^{\sigma} \|\Delta_{N,L}(u)\|_{L^2(\R\times \Theta)}
\end{eqnarray*}
and 
\begin{eqnarray}
\nonumber
C_{\sigma,b}^{-1}
\|u\|_{X^{\sigma,b}_{rad}(\R\times\Theta)}^{2} & \leq &
\sum_{L,N-{\rm dyadic }}
L^{2b}N^{2\sigma} \|\Delta_{N,L}(u)\|_{L^2(\R\times\Theta)}^{2}
\\\label{harct2}
& \leq & C_{\sigma,b}\|u\|_{X^{\sigma,b}_{rad}(\R\times\Theta)}^{2}\,.
\end{eqnarray}
Moreover there exists $C_{b}>1$, a continuous function of $b$ such that
\begin{eqnarray}
\nonumber
C_{b}^{-1}\|\Delta_{N}(u)\|_{X^{0,b}_{rad}(\R\times\Theta)}^{2}
& \leq &
\sum_{L-{\rm dyadic }}L^{2b}\|\Delta_{N,L}(u)\|_{L^2(\R\times\Theta)}^{2}
\\\label{mus1}
& \leq & C_{b}
\|\Delta_{N}(u)\|_{X^{0,b}_{rad}(\R\times\Theta)}^{2}
\end{eqnarray}
and there exists $C_{\sigma}>1$, a continuous function of $\sigma$ such that
\begin{eqnarray}
\nonumber
C_{\sigma}^{-1}\|u\|_{X^{\sigma,b}_{rad}(\R\times\Theta)}^{2}
& \leq &
\sum_{N-{\rm dyadic }}N^{2\sigma}\|\Delta_{N}(u)\|_{X^{0,b}_{rad}(\R\times\Theta)}^{2}
\\\label{mus2}
& \leq & C_{\sigma}\|u\|_{X^{\sigma,b}_{rad}(\R\times\Theta)}^{2}\,.
\end{eqnarray}
Proposition~\ref{thm3} now has a natural formulation in terms of the basic localization
projectors.
\begin{proposition}\label{thm3bis}
For every $\varepsilon>0$, there exist $\beta<1/2$ and  $C>0$ such that
for every $N_1,N_2, L_1,L_2\geq 1$, every $u_1,u_2\in 
L^{2}_{rad}(\R\times \Theta)$,
\begin{multline*}
\|\Delta_{N_1,L_1}(u_1)\Delta_{N_2,L_2}(u_2)\|_{L^{2}(\R\times\Theta)} 
\leq 
\\
\leq
C(L_1 L_2)^{\beta}
\min(N_1,N_2)^{\varepsilon}\|\Delta_{N_1,L_1}(u_1)\|_{L^{2}(\R\times\Theta)}
\|\Delta_{N_2,L_2}(u_2)\|_{L^{2}(\R\times\Theta)}\,.
\end{multline*}
\end{proposition}
\begin{proof}
It suffices to observe that $\Delta_{N_j,L_j}(u)$, $j=1,2$ satisfy the
localisation properties needed to apply Proposition~\ref{thm3}.
\end{proof}
\section{Nonlinear estimates}
The next statement contains the main analytical ingredient in the proof of
Theorem~\ref{glavna}.
\begin{proposition}\label{main}
Let $0<\sigma_1\leq\sigma<1/2$.
Then there exist two positive numbers $b,b'$ such that $b+b'<1$, $b'<1/2<b$, there
exists $C>0$ such that for every $u,v\in X^{\sigma,b}_{rad}(\R\times\Theta)$,
\begin{equation}\label{main1}
\|F(u)\|_{X^{\sigma,-b'}_{rad}(\R\times\Theta)}\leq
C\Big(1+\|u\|^{2}_{X^{\sigma_{1},b}_{rad}(\R\times\Theta)}\Big)
\|u\|_{X^{\sigma,b}_{rad}(\R\times\Theta)}
\end{equation}
and
\begin{multline}\label{main2} 
\|F(u)-F(v)\|_{X^{\sigma,-b'}_{rad}(\R\times\Theta)}
\leq 
\\
C\Big(1+\|u\|^{2}_{X^{\sigma,b}_{rad}(\R\times\Theta)}+\|v\|^{2}_{X^{\sigma,b}_{rad}(\R\times\Theta)}\Big)
\|u-v\|_{X^{\sigma,b}_{rad}(\R\times\Theta)}\,.
\end{multline}
\end{proposition}
\begin{proof}
Using the gauge invariance of the nonlinearity $F(u)$, we observe that
$F(u)-(\partial F)(0)u$ is vanishing at order $3$ at $u=0$. It therefore
suffices to prove that
\begin{equation}\label{van}
\|F(u)\|_{X^{\sigma,-b'}_{rad}(\R\times\Theta)}\leq
C\|u\|^{2}_{X^{\sigma_{1},b}_{rad}(\R\times\Theta)}
\|u\|_{X^{\sigma,b}_{rad}(\R\times\Theta)},
\end{equation}
under the additional assumption that $F(u)$ is vanishing at order $3$ at
$u=0$. Indeed, by writing
$$
\|F(u)-(\partial F)(0)u\|_{X^{\sigma,-b'}_{rad}(\R\times\Theta)}\geq
\|F(u)\|_{X^{\sigma,-b'}_{rad}(\R\times\Theta)}-C\|u\|_{X^{\sigma,-b'}_{rad}(\R\times\Theta)},
$$
we deduce that the claimed estimate (\ref{main1}) follows from (\ref{van})
applied to $F(u)-(\partial F)(0)u$. By duality, in order to prove (\ref{van}),
it suffices to establish the bound
\begin{equation}\label{vanbis}
\Big|\int_{\R\times\Theta}F(u)\bar{v}\Big|
\leq C\|v\|_{X^{-\sigma,b'}_{rad}(\R\times\Theta)}
\|u\|^{2}_{X^{\sigma_{1},b}_{rad}(\R\times\Theta)}
\|u\|_{X^{\sigma,b}_{rad}(\R\times\Theta)}.
\end{equation}
Next, we have the decompositions
$$
v=\sum_{N_0-{\rm dyadic }}
\Delta_{N_0}(v)
$$
and (recall that $F$ is smooth),
$$
F(u)=\sum_{N_{1}-{\rm dyadic }}
\Big(F(\tilde{S}_{2N_1}(u))-F(\tilde{S}_{N_1}(u))\Big)
$$
with the convention that $\tilde{S}_{1}(u)=0$. Since
$\Delta_{N}=\tilde{S}_{2N}-\tilde{S}_{N}$, we can therefore write
\begin{multline*}
F(u) =  \sum_{N_{1}-{\rm dyadic }}\Delta_{N_1}(u)G_{1}(\Delta_{N_1}(u),\tilde{S}_{N_1}(u))+
\\
\sum_{N_{1}-{\rm dyadic }}
\overline{\Delta_{N_1}(u)}G_{2}(\Delta_{N_1}(u),\tilde{S}_{N_1}(u))
\equiv  F_{1}(u)+F_{2}(u),
\end{multline*}
where $G_{1}(z_1,z_2)$ and $G_{2}(z_1,z_2)$ are smooth functions with a
control on their growth at infinity coming from (\ref{rast}).
We are going only to show that
\begin{equation*}
\Big|\int_{\R\times\Theta}F_{1}(u)\bar{v}\Big|
\leq C\|v\|_{X^{-\sigma,b'}_{rad}(\R\times\Theta)}
\|u\|^{2}_{X^{\sigma_{1},b}_{rad}(\R\times\Theta)}
\|u\|_{X^{\sigma,b}_{rad}(\R\times\Theta)}
\end{equation*}
since the argument for
$$
\Big|\int_{\R\times\Theta}F_{2}(u)\bar{v}\Big|
$$
is completely analogous.
We can write
$$
F_{1}(u)=\sum_{N_1-{\rm dyadic }}
\Delta_{N_1}(u)G_{1}(\Delta_{N_1}(u),\tilde{S}_{N_1}(u)).
$$
Next, we set
$$
I=\Big|\int_{\R\times\Theta}F_{1}(u)\bar{v}\Big|
$$
and
$$
I(N_0,N_1)=
\Big|\int_{\R\times\Theta}
\Delta_{N_1}(u)
\overline{\Delta_{N_0}(v)}
G_{1}(\Delta_{N_1}(u),\tilde{S}_{N_1}(u))
\Big|.
$$
Then $I\leq I_1+I_2$, where
$$
I_{1}=
\sum_{\stackrel{N_0\leq N_1}{N_0,N_1-{\rm dyadic }}}I(N_0,N_1),\quad
I_{2}=
\sum_{\stackrel{N_0\geq N_1}{N_0,N_1-{\rm dyadic }}}I(N_0,N_1).
$$
We estimate first $I_1$. 
Similarly to the above expansion for $F$, 
using the vanishing property at the origin of $F$,
we now decompose $G_{1}(\Delta_{N_1}(u),\tilde{S}_{N_1}(u))$ as follows,
$$
\sum_{N_2-{\rm dyadic }}
\Big(G_{1}(\tilde{S}_{2N_2}\Delta_{N_1}(u),\tilde{S}_{2N_2}\tilde{S}_{N_1}(u))-G_{1}(\tilde{S}_{N_2}\Delta_{N_1}(u),\tilde{S}_{N_2}\tilde{S}_{N_1}(u))\Big).
$$
Therefore, using that $\Delta_{N_1}\Delta_{N_2}=\Delta_{N_1}$, if $N_1=N_2$
and zero elsewhere, we obtain that
\begin{multline*}
G_{1}(\Delta_{N_1}(u),\tilde{S}_{N_1}(u))=
\sum_{\stackrel{N_2\leq N_1}{N_2-{\rm dyadic }}}
\Delta_{N_2}(u)G_{11}^{N_2}(\Delta_{N_2}(u),\tilde{S}_{N_2}(u))
+
\\
\sum_{\stackrel{N_2\leq N_1}{N_2-{\rm dyadic }}}
\overline{\Delta_{N_2}(u)}G_{12}^{N_2}(\Delta_{N_2}(u),\tilde{S}_{N_2}(u)).
\end{multline*}
Finally, we expand for $j=1,2$,
\begin{multline*}
G_{1j}^{N_2}(\Delta_{N_2}(u),\tilde{S}_{N_2}(u))=
\sum_{\stackrel{N_3\leq N_2}{N_3-{\rm dyadic }}}
\Delta_{N_3}(u)G_{1j1}^{N_3}(\Delta_{N_3}(u),\tilde{S}_{N_3}(u))
+
\\
\sum_{\stackrel{N_3\leq N_2}{N_3-{\rm dyadic }}}
\overline{\Delta_{N_3}(u)}G_{1j2}^{N_3}(\Delta_{N_3}(u),\tilde{S}_{N_3}(u)),
\end{multline*}
where, thanks to the growth assumption on the nonlinearity $F(u)$, we obtain
that the functions $G_{1 j_1 j_2}^{N_3}(z_1,z_2)$, $j_1,j_2=1,2$ satisfy
$$
|G_{1 j_1 j_2}^{N_3}(z_1,z_2)|\leq C.
$$
We therefore have the bound
$$
I_1\leq C
\sum_{\stackrel{N_0\leq N_1}{N_0,N_1-{\rm dyadic }}}
\sum_{\stackrel{N_1\geq N_2\geq N_3}{N_2,N_3-{\rm dyadic }}}
\int_{\R\times\Theta}
|\Delta_{N_0}(v)\Delta_{N_1}(u)\Delta_{N_2}(u)\Delta_{N_3}(u)|
$$
and moreover using the equality
$$
\Delta_{N}=\sum_{L-{\rm dyadic}}\Delta_{N,L},
$$
we arrive at
\begin{multline*}
I_1\leq C
\sum_{L_0,L_1,L_2,L_3-{\rm dyadic }}
\\
\sum_{\stackrel{N_1\geq N_2\geq N_3, N_1\geq N_0}{N_0,N_1,N_2,N_3-{\rm dyadic }}}
\int_{\R\times\Theta}
|\Delta_{N_0,L_0}(v)\Delta_{N_1,L_1}(u)\Delta_{N_2,L_2}(u)\Delta_{N_3,L_3}(u)|.
\end{multline*}
Using Proposition~\ref{thm3bis} and the Cauchy-Schwarz inequality, we have
that for every $\varepsilon>0$ there exist $\beta<1/2$ and $C_{\varepsilon}$
such that
\begin{multline*}
\int_{\R\times\Theta}
|\Delta_{N_0,L_0}(v)\Delta_{N_1,L_1}(u)\Delta_{N_2,L_2}(u)\Delta_{N_3,L_3}(u)|\leq
\\
\leq
\|\Delta_{N_0,L_0}(v)\Delta_{N_2,L_2}(u)\|_{L^2(\R\times\Theta)}
\|\Delta_{N_1,L_1}(u)\Delta_{N_3,L_3}(u)\|_{L^2(\R\times\Theta)}
\leq
\\
\leq C_{\varepsilon}(N_2 N_3)^{\varepsilon}(L_0L_1L_2L_3)^{\beta}
\|\Delta_{N_0,L_0}(v)\|_{L^2(\R\times\Theta)}
\prod_{j=1}^{3}\|\Delta_{N_j,L_j}(u)\|_{L^2(\R\times\Theta)}.
\end{multline*}
Therefore, if we set
\begin{multline}\label{Q}
Q\equiv Q(N_0,N_1,N_2,N_3,L_0,L_1,L_2,L_3)= 
CN_{0}^{-\sigma}N_1^{\sigma}(N_2 N_3)^{\sigma_{1}}L_{0}^{b'}(L_1 L_2 L_3)^{b}
\\
\times
\|\Delta_{N_0,L_0}(v)\|_{L^2(\R\times\Theta)}
\prod_{j=1}^{3}\|\Delta_{N_j,L_j}(u)\|_{L^2(\R\times\Theta)},
\end{multline}
we can write
$$
I_1\leq \sum_{L_0,L_1,L_2,L_3-{\rm dyadic }}\,\sum_{\stackrel{N_1\geq N_2\geq N_3, N_1\geq N_0}{N_0,N_1,N_2,N_3-{\rm dyadic }}}
L_{0}^{\beta-b'}(L_1L_2L_3)^{\beta-b}\Big(\frac{N_0}{N_1}\Big)^{\sigma}(N_2N_3)^{\varepsilon-\sigma_{1}}Q
$$
Let us take $\varepsilon>0$ such that $\varepsilon-\sigma_{1}<0$. 
This fixes $\beta$.
Then we choose
$b'$ such that $\beta<b'<1/2$. We finally choose $b>1/2$ such that
$b+b'<1$. With this choice of the parameters, using (\ref{mus1}) and
after summing geometric series in $L_0$, $L_1$, $L_2$, $L_3$, $N_2$, $N_3$, we can write that
$$
I_{1}\leq
C\|u\|^{2}_{X^{\sigma_{1},b}_{rad}(\R\times\Theta)}
\sum_{\stackrel{ N_0\leq N_1}{N_0,N_1-{\rm dyadic }}}
\Big(\frac{N_0}{N_1}\Big)^{\sigma}c(N_0)d(N_1),
$$
where
\begin{equation}\label{cd}
c(N_0)=N_{0}^{-\sigma}\|\Delta_{N_0}(v)\|_{X^{0,b'}_{rad}(\R\times\Theta)},\quad
d(N_1)=N_{1}^{\sigma}\|\Delta_{N_1}(u)\|_{X^{0,b}_{rad}(\R\times\Theta)}\,.
\end{equation}
We now make appeal to the following lemma which is a discreet variant of the
Schur test.
\begin{lemme}\label{precedent}
For every $\sigma>0$ there exists $C>0$ such that for every couple of
functions $c_{j}(N)$, $j=1,2$, defined on the set of the dyadic integers such
that
$$
\|c_j\|\equiv\sum_{N-{\rm dyadic }}|c_{j}(N)|^{2}<\infty,\quad j=1,2
$$
one has
\begin{equation}\label{ineg}
\Big|\sum_{\stackrel{ N_0\leq N_1}{N_0,N_1-{\rm dyadic }}}\Big(\frac{N_0}{N_1}\Big)^{\sigma}c_{1}(N_0)c_{2}(N_1)\Big|
\leq C\|c_1\|\|c_2\|.
\end{equation}
\end{lemme}
\begin{proof}
Write $N_{1}=2^{j}N_{0}$ with $j\geq 0$, $j\in\Z$. Thus
the left hand-side of (\ref{ineg}) can be rewritten as
$$
\Big|
\sum_{j=0}^{\infty}\sum_{N_0-{\rm dyadic }}2^{-j\sigma}c_{1}(N_0)c_{2}(2^{j}N_0)
\Big|
$$
which by the Cauchy-Schwartz inequality in $N_0$ is bounded by
$C\|c_1\|\|c_2\|$ with
$$
C=\sum_{j=0}^{\infty}2^{-j\sigma}\,.
$$
This completes the proof of Lemma~\ref{precedent}.
\end{proof}
Next using (\ref{mus2}) and Lemma~\ref{precedent}, we deduce that
$$
I_1\leq
C\|v\|_{X^{-\sigma,b'}_{rad}(\R\times\Theta)}
\|u\|^{2}_{X^{\sigma_{1},b}_{rad}(\R\times\Theta)}
\|u\|_{X^{\sigma,b}_{rad}(\R\times\Theta)}.
$$
This ends the analysis for $I_1$.
We next turn to the estimate for $I_2$.
The basic idea is that after an integration by parts, the structure of $I_2$
becomes very close to the structure of $I_1$, by simply exchanging the roles of
$N_0$ and $N_1$. In this context the Proposition~\ref{thm3triis}
gives the relevant bound. We start by some preliminary observations.
For $u\in L^{2}_{rad}(\R\times\Theta)$ we can write 
\begin{equation*}
\Delta_{N,L}(u)=\sum_{N\leq \langle z_{n}\rangle < 2N}\,c(n,t)\, e_{n}(r),
\end{equation*}
where
$$
{\rm supp}\, \widehat{c}(n,\tau)\subset 
\{
\tau\in\R\,:\, L\leq \langle\tau+z_{n}^2\rangle\leq 2L
\}
$$
and
$$
\|\Delta_{N,L}(u)\|_{L^2(\R\times\Theta)}^{2}
=c\sum_{N\leq \langle z_{n}\rangle < 2N}
\int_{-\infty}^{\infty}|\widehat{c}(n,\tau)|^{2}d\tau.
$$
Moreover
$$
\partial_{r}\Big(\Delta_{N,L}(u)\Big)=
\sum_{N\leq \langle z_{n}\rangle < 2N}\,c(n,t)\, e'_{n}(r).
$$
Recall that for $m\neq n$, $e'_{m}$ and $e'_{n}$ are orthogonal in $L^2(\Theta)$. 
Moreover, thanks to (\ref{equiv}),
$$
\|e'_{n}\|_{L^2(\Theta)}\approx n\|e_{n}\|_{L^2(\Theta)}
$$
and thus using that
$$
\big\|\partial_{r}\big(\Delta_{N,L}(u)\big)\big\|_{L^2(\R\times\Theta)}^{2}=
c\sum_{N\leq \langle z_{n}\rangle < 2N}\,
\|e'_{n}\|_{L^2(\Theta)}^{2}
\int_{-\infty}^{\infty}
|\widehat{c}(n,\tau)|^{2}d\tau
$$
we arrive at the crucial relation
\begin{equation}\label{jm}
\big\|\partial_{r}\big(\Delta_{N,L}(u)\big)\big\|_{L^2(\R\times\Theta)}\approx 
N\big\|\Delta_{N,L}(u)\big\|_{L^2(\R\times\Theta)}\,.
\end{equation}
Let us observe that
$$
e_{n}(r)=-\frac{1}{z_n^2}\frac{1}{r}\partial_{r}(r\partial_{r}e_{n}(r)).
$$
Since $\Delta_{N_1}(u)G_{1}(\Delta_{N_1}(u),\tilde{S}_{N_1}(u))$ is vanishing on the
boundary of $\Theta$, an integration by parts yields
\begin{multline*}
\int_{0}^{1}e_{n}(r)\Delta_{N_1}(u)G_{1}(\Delta_{N_1}(u),\tilde{S}_{N_1}(u))rdr
=
\\
=
\frac{1}{z_n^2}\int_{0}^{1}e'_{n}(r)\partial_{r}\Big(\Delta_{N_1}(u)G_{1}(\Delta_{N_1}(u),\tilde{S}_{N_1}(u))\Big)rdr\,.
\end{multline*}
Write
\begin{equation}\label{write1}
\Delta_{N_0,L_0}(v)=\sum_{N_0\leq \langle z_{n_0}\rangle < 2N_0}\,c(n_0,t)\, e_{n_0}(r),
\end{equation}
where
$$
{\rm supp}\, \widehat{c}(n_0,\tau)\subset 
\{
\tau\in\R\,:\, L_{0}\leq \langle\tau+z_{n_0}^2\rangle\leq 2L_0
\}.
$$
Then, for $n\in\N$ such that $\langle z_n\rangle\in[N_0,2N_0[$, we set
$$
\widetilde{c}(n,t)=\frac{c(n,t)}{z_n^2}\,,
$$
where $c(n,t)$ are the coefficients involved in (\ref{write1}). 
Define $\widetilde{\Delta}_{N_0,L_0}$ as
\begin{equation*}
\widetilde{\Delta}_{N_0,L_0}(v)=\sum_{N_0\leq \langle z_{n_0}\rangle \leq 2N_0}\,\widetilde{c}(n_0,t)\, e'_{n_0}(r).
\end{equation*}
Clearly $\widetilde{\Delta}_{N_0,L_0}(v)$ is an object which fits in the scope
of applicability of Proposition~\ref{thm3triis} and
\begin{equation}\label{jm2}
\|\widetilde{\Delta}_{N_0,L_0}(v)\|_{L^2(\R\times\Theta)}
\approx
N_{0}^{-1}
\|\Delta_{N_0,L_0}(v)\|_{L^2(\R\times\Theta)}.
\end{equation}
Recall that $e_{n}(r)$ are real valued. In view of the above discussion,
we need to control the expression 
$$
E=
\sum_{\stackrel{N_0\geq N_1}{L_0,N_0,N_1-{\rm dyadic }}}
\Big|\int_{\R\times\Theta}\overline{\widetilde{\Delta}_{N_0,L_0}(v)}
\partial_{r}\Big(\Delta_{N_1}(u)G_{1}(\Delta_{N_1}(u),\tilde{S}_{N_1}(u))\Big)\Big|.
$$
Now, we can write $E\leq E_1+E_2$, where
$$
E_1=
\sum_{\stackrel{N_0\geq N_1}{L_0,N_0,N_1-{\rm dyadic }}}
\Big|\int_{\R\times\Theta}\overline{\widetilde{\Delta}_{N_0,L_0}(v)}
\partial_{r}\Big(\Delta_{N_1}(u)\Big)G_{1}(\Delta_{N_1}(u),\tilde{S}_{N_1}(u))\Big|
$$
and
$$
E_2=
\sum_{\stackrel{N_0\geq N_1}{L_0,N_0,N_1-{\rm dyadic }}}
\Big|\int_{\R\times\Theta}\overline{\widetilde{\Delta}_{N_0,L_0}(v)}
\Delta_{N_1}(u)\partial_{r}\Big(G_{1}(\Delta_{N_1}(u),\tilde{S}_{N_1}(u))\Big)\Big|.
$$
By expanding $G_1(z_1,z_2)$ and using the growth and vanishing assumptions on the nonlinear
interaction $F$, we can write
\begin{multline*}
E_1\leq C
\sum_{L_0,L_1,L_2,L_3-{\rm dyadic }}
\\
\sum_{\stackrel{N_0\geq N_1\geq N_2\geq N_3}{N_0,N_1,N_2,N_3-{\rm dyadic }}}
\int_{\R\times\Theta}
|\widetilde{\Delta}_{N_0,L_0}(v)
\partial_{r}\big(\Delta_{N_1,L_1}(u)\big)\Delta_{N_2,L_2}(u)\Delta_{N_3,L_3}(u)|.
\end{multline*}
Using Proposition~\ref{thm3bis},  Proposition~\ref{thm3triis}, the
Cauchy-Schwarz inequality, (\ref{jm}) and (\ref{jm2}), we have
that for every $\varepsilon>0$ there exist $\beta<1/2$ and $C_{\varepsilon}$
such that
\begin{multline*}
\int_{\R\times\Theta}
|\widetilde{\Delta}_{N_0,L_0}(v)
\partial_{r}\big(\Delta_{N_1,L_1}(u)\big)\Delta_{N_2,L_2}(u)\Delta_{N_3,L_3}(u)|\leq
\\
\leq\|\widetilde{\Delta}_{N_0,L_0}(v)\Delta_{N_2,L_2}(u)\|_{L^2(\R\times\Theta)}
\|
\partial_{r}\big(\Delta_{N_1,L_1}(u)\big)
\Delta_{N_3,L_3}(u)\|_{L^2(\R\times\Theta)}
\leq
\\
\leq C_{\varepsilon}(N_2 N_3)^{\varepsilon}(L_0L_1L_2L_3)^{\beta}
\|
\widetilde{\Delta}_{N_0,L_0}(v)
\|_{L^2(\R\times\Theta)}
\\
\times
\|\partial_{r}\big(\Delta_{N_1,L_1}(u)\big)\|_{L^2(\R\times\Theta)}
\prod_{j=2}^{3}\|\Delta_{N_j,L_j}(u)\|_{L^2(\R\times\Theta)}\leq
\\
\leq C_{\varepsilon}(N_2 N_3)^{\varepsilon}(L_0L_1L_2L_3)^{\beta}
\frac{N_1}{N_0}
\|\Delta_{N_0,L_0}(v)\|_{L^2(\R\times\Theta)}
\prod_{j=1}^{3}\|\Delta_{N_j,L_j}(u)\|_{L^2(\R\times\Theta)}.
\end{multline*}
Therefore, with $Q$ defined by (\ref{Q}), we can write
\begin{multline}\label{same}
E_1\leq
\sum_{L_0,L_1,L_2,L_3-{\rm dyadic }}
\\
\sum_{\stackrel{N_0\geq N_1\geq N_2\geq N_3}{N_0,N_1,N_2,N_3-{\rm dyadic }}}
L_{0}^{\beta-b'}(L_1L_2L_3)^{\beta-b}\Big(\frac{N_0}{N_1}\Big)^{\sigma-1}(N_2N_3)^{\varepsilon-\sigma_{1}}Q
\end{multline}
Let us take $\varepsilon>0$ such that $\varepsilon-\sigma_{1}<0$. 
Then as we did for the bound for $I_1$, we choose
$b'$ such that $\beta<b'<1/2$. We finally choose $b>1/2$ such that
$b+b'<1$. Using (\ref{mus1}) and
after summing geometric series in $L_0$, $L_1$, $L_2$, $L_3$, $N_2$, $N_3$, we can write that
$$
E_{1}\leq
C\|u\|^{2}_{X^{\sigma_{1},b}_{rad}(\R\times\Theta)}
\sum_{\stackrel{ N_0\geq N_1}{N_0,N_1-{\rm dyadic }}}
\Big(\frac{N_1}{N_0}\Big)^{1-\sigma}c(N_0)d(N_1),
$$
where $c(N_0)$ and $d(N_1)$ are defined by (\ref{cd}).
Therefore,
using (\ref{mus2}) and Lemma~\ref{precedent}, we arrive at the bound
$$
E_1\leq
C\|v\|_{X^{-\sigma,b'}_{rad}(\R\times\Theta)}
\|u\|^{2}_{X^{\sigma_{1},b}_{rad}(\R\times\Theta)}
\|u\|_{X^{\sigma,b}_{rad}(\R\times\Theta)}.
$$
Let us now turn to the bound for $E_2$. 
Using the formula
$$
\partial_{r}\Big(f(z(r))\Big)=(\partial_{r}z)\partial
f+(\partial_{r}\bar{z})\bar{\partial}f, 
$$
we can write
\begin{multline*}
\partial_{r}\Big(G_{1}(\Delta_{N_1}(u),\tilde{S}_{N_1}(u))\Big)=
\sum_{\stackrel{N_2\leq N_1}{N_2-{\rm dyadic }}}
\partial_{r}\Big(\Delta_{N_2}(u)\Big)\widetilde{G}_{11}^{N_2}(\Delta_{N_1}(u),\tilde{S}_{N_1}(u))
+
\\
\sum_{\stackrel{N_2\leq N_1}{N_2-{\rm dyadic }}}
\overline{\partial_{r}\Big(\Delta_{N_2}(u)\Big)}\widetilde{G}_{12}^{N_2}(\Delta_{N_1}(u),\tilde{S}_{N_1}(u)),
\end{multline*}
where thanks to the growth assumption on the nonlinearity, $\widetilde{G}_{1j}^{N_2}(z_1,z_2)$, $j=1,2$ 
satisfy
\begin{equation}\label{kino}
\sum_{k=1}^{2}
\Big(|\partial_{z_k}\widetilde{G}_{1j}^{N_2}(z_1,z_2)|+|\bar{\partial}_{z_k}\widetilde{G}_{1j}^{N_2}(z_1,z_2)|\Big)
\leq C.
\end{equation}
By expanding $\widetilde{G}_{1j}^{N_2}(z_1,z_2)$, $j=1,2$ in a telescopic series and using
(\ref{kino}), we get the bound
\begin{multline*}
E_2\leq C
\sum_{L_0,L_1,L_2,L_3-{\rm dyadic }}
\\
\sum_{\stackrel{N_0\geq N_1\geq N_2, N_1\geq N_3}{N_0,N_1,N_2,N_3-{\rm dyadic }}}
\int_{\R\times\Theta}
|\widetilde{\Delta}_{N_0,L_0}(v)
\Delta_{N_1,L_1}(u)
\partial_{r}\big(\Delta_{N_2,L_2}(u)\big)\Delta_{N_3,L_3}(u)|.
\end{multline*}
Using Proposition~\ref{thm3bis},  Proposition~\ref{thm3triis}, the
Cauchy-Schwarz inequality, (\ref{jm}) and (\ref{jm2}), we have
that for every $\varepsilon>0$ there exist $\beta<1/2$ and $C_{\varepsilon}$
such that
\begin{multline*}
\int_{\R\times\Theta}
|\widetilde{\Delta}_{N_0,L_0}(v)
\partial_{r}\big(\Delta_{N_1,L_1}(u)\big)\Delta_{N_2,L_2}(u)\Delta_{N_3,L_3}(u)|\leq
\\
\leq\|\widetilde{\Delta}_{N_0,L_0}(v)\Delta_{N_3,L_3}(u)\|_{L^2(\R\times\Theta)}
\|
\partial_{r}\big(\Delta_{N_2,L_2}(u)\big)
\Delta_{N_1,L_1}(u)\|_{L^2(\R\times\Theta)}
\leq
\\
\leq C_{\varepsilon}(N_2 N_3)^{\varepsilon}(L_0L_1L_2L_3)^{\beta}
\|
\widetilde{\Delta}_{N_0,L_0}(v)
\|_{L^2(\R\times\Theta)}
\\
\times
\|\partial_{r}\big(\Delta_{N_2,L_3}(u)\big)\|_{L^2(\R\times\Theta)}
\|\Delta_{N_1,L_1}(u)\|_{L^2(\R\times\Theta)}\|\Delta_{N_3,L_3}(u)\|_{L^2(\R\times\Theta)}\leq
\\
\leq C_{\varepsilon}(N_2 N_3)^{\varepsilon}(L_0L_1L_2L_3)^{\beta}
\frac{N_2}{N_0}
\|\Delta_{N_0,L_0}(v)\|_{L^2(\R\times\Theta)}
\prod_{j=1}^{3}\|\Delta_{N_j,L_j}(u)\|_{L^2(\R\times\Theta)}.
\end{multline*}
Next, for $N_2\leq N_1$, we can write,
$$
\Big(\frac{N_0}{N_1}\Big)^{\sigma}
\frac{N_2}{N_0}
\leq
\Big(\frac{N_0}{N_1}\Big)^{\sigma-1}
$$
and therefore, with $Q$ defined by (\ref{Q}), we can write
$$
E_2\leq
\sum_{L_0,L_1,L_2,L_3-{\rm dyadic }}\,\sum_{\stackrel{N_0\geq N_1\geq N_2, N_1\geq N_3}{N_0,N_1,N_2,N_3-{\rm dyadic }}}
L_{0}^{\beta-b'}(L_1L_2L_3)^{\beta-b}\Big(\frac{N_0}{N_1}\Big)^{\sigma-1}(N_2N_3)^{\varepsilon-\sigma_{1}}Q
$$
But the right hand-side of the above inequality is exactly the same as the 
right hand-side of of (\ref{same}). Therefore
$$
E_2\leq
C\|v\|_{X^{-\sigma,b'}_{rad}(\R\times\Theta)}
\|u\|^{2}_{X^{\sigma_{1},b}_{rad}(\R\times\Theta)}
\|u\|_{X^{\sigma,b}_{rad}(\R\times\Theta)}.
$$
This completes the proof of (\ref{main1}).
In the proof of (\ref{main1}), we analysed the expression 
$
\|F(u)\|_{X^{\sigma,-b'}_{rad}}.
$
The argument is based on successive expansions of $F(u)$ in telescopic series
and thus it works equally well if we replace $F(u)$ by $u\,G(v,w)$ where $G(z_1,z_2)$
satisfies the growth assumption
\begin{equation}\label{rastbis}
\big|
\partial^{k_1}_{z_1}\bar{\partial}^{k_2}_{z_1}
\partial^{l_1}_{z_2}\bar{\partial}^{l_2}_{z_2}
G(z_1,z_2)
\big|\leq C_{k_1,k_2,l_1,l_2}(1+|z_1|+|z_2|)^{\max(0,\alpha-k_1-k_2-l_1-l_2)}
\,.
\end{equation}
But this is exactly the situation that occurs in the analysis of
(\ref{main2}).
Indeed, one can write
$$
F(u)-F(v)=(u-v)G_{1}(u,v)+(\overline{u}-\overline{v})G_{2}(u,v)
$$
with $G_{j}(z_1,z_2)$, $j=1,2$ satisfying (\ref{rastbis}).
Since the analysis is very similar to the proof of (\ref{main1}), we shall
only outline the estimate for $(u-v)G_{1}(u,v)$. Again, we can suppose that
$F(u)$ is vanishing at order $3$ at $u=0$. Let us set
$$
w_1=u-v,\quad w_2=u,\quad w_3=v.
$$
One needs to bound
$$
\Big|
\int_{\R\times\Theta}
w_1G_{1}(w_2,w_3)\overline{w_4}
\Big|
$$
by
$$
C(1+\|w_2\|_{X^{\sigma,b}_{rad}(\R\times\Theta)}+\|w_3\|_{X^{\sigma,b}_{rad}(\R\times\Theta)})^2
\|w_1\|_{X^{\sigma,b}_{rad}(\R\times\Theta)}\|w_4\|_{X^{-\sigma,b'}_{rad}(\R\times\Theta)}\,.
$$
Next, we expand
$$
w_{1}=\sum_{N_1-{\rm dyadic}}\Delta_{N_1}(w_1),\quad w_{4}=\sum_{N_0-{\rm dyadic}}\Delta_{N_0}(w_4)
$$
and
$$
G_{1}(w_2,w_3)=\sum_{N_2-{\rm dyadic}}\Big(G_{1}(\tilde{S}_{2N_2}(w_2),\tilde{S}_{2N_2}(w_3))-G_{1}(\tilde{S}_{N_2}(w_2),\tilde{S}_{N_2}(w_3))\Big).
$$
Thus, modulo complex conjugations irrelevant in this discussion, one has to
evaluate quantities of type
\begin{multline}\label{derm}
\sum_{N_0,N_1,N_2-{\rm dyadic}}
\Big|
\int_{\R\times\Theta}\overline{\Delta_{N_0}(w_4)}\Delta_{N_1}(w_1)\Delta_{N_2}(w_j)
\\
H_{j}^{N_2}(\Delta_{N_2}(w_2),\tilde{S}_{N_2}(w_2),\Delta_{N_2}(w_3),\tilde{S}_{N_2}(w_3))
\Big|,\quad j=2,3,
\end{multline}
where $H_{j}^{N_2}(z_1,z_2,z_3,z_4)$ are smooth functions satisfying growth
restrictions at infinity coming from (\ref{rast}). In the analysis of
(\ref{derm}), we distinguish two cases for $N_0$, $N_1$, $N_2$ in the sum
defining (\ref{derm}). The first case is when $N_{0}\leq \max(N_1,N_2)$, In
this case, we expand once more $H_{j}^{N_2}$ which introduces a sum over $N_{3}-{\rm
dyadic}$, $N_{3}\leq N_2$ of terms $\Delta_{N_3}(w_{j})$ (or complex conjugate)
times a bounded function (thanks to the sub cubic nature of the nonlinearity).
The analysis is then exactly the same as for that of $I_1$ in the proof of
(\ref{main1}). If $N_{0}\geq \max(N_1,N_2)$, then we integrate by parts by the aid of
$\Delta_{N_0}(w_4)$ and analysis is the same as in the bound for 
$I_2$ in the proof of (\ref{main1}). This completes the proof of Proposition~\ref{main}.
\end{proof}
Let us now consider the integral equation corresponding to the problem
(\ref{1})-(\ref{2})
\begin{equation}\label{Duhamel}
u(t)=e^{it\Delta}u_0+i\int_{0}^{t}e^{i(t-\tau)\Delta}F(u(\tau))d\tau\,.
\end{equation}
With Proposition~\ref{main} in hand, we can deduce
the following estimates for the terms in the right hand-side of
(\ref{Duhamel}).
\begin{proposition}\label{duhbis}
Let $0<\sigma_1\leq \sigma<1/2$.
Then there exist two positive numbers $b,b'$ such
that $b+b'<1$, $b'<1/2<b$, there
exists $C>0$ such that for every $T\in]0,1]$, 
every $u,v\in X^{\sigma,b}_{rad}([-T,T]\times\Theta)$, every $u_0\in
H^{\sigma}_{rad}(\Theta)$, 
\begin{equation}\label{linear}
\big\|e^{it\Delta}u_0 \big\|_{X^{\sigma,b}_{rad}([-T,T]\times\Theta)}\leq
C\|u_0\|_{H^{\sigma}_{rad}(\Theta)}\, ,
\end{equation}
\begin{multline}\label{non1}
\Big\|
\int_{0}^{t}e^{i(t-\tau)\Delta}F(u(\tau))d\tau
\Big\|_{X^{\sigma,b}_{rad}([-T,T]\times\Theta)}
\leq
\\
\leq
CT^{1-b-b'}\Big(1+\|u\|^{2}_{X^{\sigma_{1},b}_{rad}([-T,T]\times\Theta)}\Big)
\|u\|_{X^{\sigma,b}_{rad}([-T,T]\times\Theta)}
\end{multline}
and
\begin{multline}\label{non2}
\Big\|
\int_{0}^{t}e^{i(t-\tau)\Delta}(F(u(\tau))-F(v(\tau)))d\tau
\Big\|_{X^{\sigma,b}_{rad}([-T,T]\times\Theta)}
\leq
\\
\leq CT^{1-b-b'}\Big(1+\|u\|^{2}_{X^{\sigma,b}_{rad}([-T,T]\times\Theta)}+\|v\|^{2}_{X^{\sigma,b}_{rad}([-T,T]\times\Theta)}\Big)
\|u-v\|_{X^{\sigma,b}_{rad}([-T,T]\times\Theta)}\,.
\end{multline}
\end{proposition}
\begin{proof}
Let $\psi\in C_{0}^{\infty}(\R)$ such that $\psi\equiv 1$ on $[-1,1]$.
Then, using (\ref{libre}), we can write
$$
\big\|e^{it\Delta}u_0 \big\|_{X^{\sigma,b}_{rad}([-T,T]\times\Theta)}\leq
\big\|\psi(t)e^{it\Delta}u_0 \big\|_{X^{\sigma,b}_{rad}(\R\times\Theta)}
=\|\psi\|_{H^b(\R)}\|u_0\|_{H^{\sigma}_{rad}(\Theta)}
$$
which proves (\ref{linear}). 
Let us remark that if 
$\tilde{u}\in X^{\sigma,b}_{rad}(\R\times\Theta)$
realises the $X^{\sigma,b}_{rad}([-T,T]\times\Theta)$ norm of $u$ then the
same $\tilde{u}$ realises all $X^{\sigma',b}_{rad}([-T,T]\times\Theta)$, $\sigma'<\sigma$ norms of $u$.
With this remark in hand, now the proofs of (\ref{non1}) and (\ref{non2}) follow from (\ref{main1}) and
(\ref{main2}) respectively, (\ref{libre}) and the inequality
\begin{equation}\label{ginibre}
\|\psi(t/T)\int_{0}^{t}f(\tau)d\tau\|_{H^{b}(\R)}\leq CT^{1-b-b'}\|f\|_{H^{-b'}(\R)}\,.
\end{equation}
We refer to \cite{Gi} for the proof of (\ref{ginibre}).
This completes the proof of Proposition~\ref{duhbis}.
\end{proof}
\section{Basic local well-posedness results for NLS and the truncated NLS}
Recall that we are interested in constructing solutions of the initial value problem
\begin{equation}\label{pak}
iu_t+\Delta u +F(u)=0, \quad u|_{t=0}=u_0 \,.
\end{equation}
We will approximate the solutions of (\ref{pak}) by the solutions of the ODE
\begin{equation}\label{pakbis}
iu_t+\Delta u +S_{N}(F(u))=0, \quad u|_{t=0}=u_{0}\in E_{N},
\end{equation}
for $N\gg 1$ (for the definition of the projector $S_N$, see Section~3, (\ref{formula}) above).
Equation (\ref{pakbis}) can be seen as a Hamiltonian ODE for $u=S_{N}(u)$.
More precisely, if 
$$
u=S_{N}(u)=\sum_{n=1}^{N}c_n\, e_{n,s},
$$
then the Hamiltonian of the ODE (\ref{pakbis}) is given by
$$
H(u,\overline{u})=
\sum_{n=1}^{N} z_n^{2-2s}\, |c_n|^2 - \int_{0}^1
V\Big(
\sum_{m= 1}^{N}c_m\, e_{m,s}(r)\Big)rdr\, .
$$
Multiplying (\ref{pakbis}) by $\overline{u}$ and integrating over $\Theta$ yields that
the $L^2$ norm is still a conserved quantity for (\ref{pakbis}). Therefore,
the Cauchy-Lipschitz theorem for ODE's implies the existence of global
dynamics for (\ref{pakbis}) for every $u_0\in L^2(\Theta)$. 
The $L^2$ conservation provides the bound
$$
\sum_{n=1}^{N}z_{n}^{-2s}|c_{n}(t)|^{2}\leq C
$$
uniformly in $t$. 
However, bounds on the quantities 
$$
\sum_{n=1}^{N}n^{\sigma}n^{-2s}|c_{n}(t)|^{2},\quad \sigma>0
$$
for long times are non trivial and go beyond the scope of the basic
Cauchy-Lipschitz theorem.
We next state the basic local well-posedness result for (\ref{pak}), which
unfortunately misses the $L^2$ theory.
\begin{proposition}\label{lwp}
Let us fix $0<\sigma_{1}\leq\sigma<1/2$.
Then there exist $b>1/2$, $\beta>0$, $C>0$, $\tilde{C}>0$, $c\in]0,1]$ such that for every 
$A>0$ if we set $T=c(1+A)^{-\beta}$ then for every $u_0\in
H^{\sigma_{1}}_{rad}(\Theta)$ satisfying $\|u_0\|_{H^{\sigma_{1}}}\leq A$ there exists a unique
solution of (\ref{pak}) in $X^{\sigma_{1},b}_{rad}([-T,T]\times \Theta)$. Moreover
$$
\|u\|_{L^{\infty}([-T,T];H^{\sigma_{1}}(\Theta))}
\leq 
C\|u\|_{X^{\sigma_{1},b}_{rad}([-T,T]\times \Theta)}
\leq \tilde{C}\|u_0\|_{H^{\sigma_{1}}(\Theta)}\, .
$$ 
If in addition $u_0\in H^{\sigma}_{rad}(\Theta)$ then
$$
\|u\|_{L^{\infty}([-T,T];H^{\sigma}(\Theta))}\leq 
C\|u\|_{X^{\sigma,b}_{rad}([-T,T]\times \Theta)}
\leq \tilde{C}\|u_0\|_{H^{\sigma}(\Theta)}\, .
$$ 
Finally if $u$ and $v$ are two solutions with data $u_0$, $v_0$ respectively,
satisfying 
$$
\|u_0\|_{H^{\sigma_{1}}}\leq A,\quad \|v_0\|_{H^{\sigma_{1}}}\leq A
$$
then 
$$
\|u-v\|_{L^{\infty}([-T,T];H^{\sigma_{1}}(\Theta))}\leq C\|u_0-v_0\|_{H^{\sigma_{1}}(\Theta)}\, .
$$ 
If in addition $u_0,v_0\in H^{\sigma}_{rad}(\Theta)$ then
$$
\|u-v\|_{L^{\infty}([-T,T];H^{\sigma}(\Theta))}\leq C\|u_0-v_0\|_{H^{\sigma}(\Theta)}\, .
$$
\end{proposition}
\begin{proof}
It is a direct application of  Proposition~\ref{duhbis} and the contraction
mapping principle to the map $\Phi_{u_0}(u)$ defined by the right hand-side of
(\ref{Duhamel}).
Indeed, for shortness, let us denote by $X^{\sigma}_{T}$ the Bourgain space 
$X^{\sigma,b}_{rad}([-T,T]\times \Theta)$, where $b$ is fixed in
Proposition~\ref{duhbis}. Then there exists $\theta>0$ ($\theta\equiv 1-b-b'$
with $b'$ fixed by Proposition~\ref{duhbis}) such that 
\begin{equation}\label{contr1}
\|\Phi_{u_0}(u)\|_{X^{\sigma_1}_{T}}\leq 
C\|u_0\|_{H^{\sigma_1}_{rad}(\Theta)}+CT^{\theta}(1+\|u\|_{X_{T}^{\sigma_{1}}}^{2})\|u\|_{X_{T}^{\sigma_1}}
\end{equation}
and
\begin{equation}\label{contr2}
\|\Phi_{u_0}(u)-\Phi_{u_0}(v)\|_{X^{\sigma_1}_{T}}\leq 
CT^{\theta}\|u-v\|_{X_{T}^{\sigma_1}}(1+\|u\|_{X_{T}^{\sigma_1}}^{2}+\|v\|_{X_{T}^{\sigma_1}}^{2})\,.
\end{equation}
Using (\ref{contr1}), we obtain that for every $u_0$ such that
$\|u_0\|_{H^{\sigma_{1}}}\leq A$ if we take 
$$T\sim(1+A)^{-\theta/2}$$ 
then the map $\Phi_{u_0}$ sends the ball of radius $2C\|u_0\|_{H^{\sigma_{1}}}$
of $X^{\sigma_1}_{T}$ into the same ball.
Thanks to (\ref{contr2}), with the same restriction on $T$ 
the map $\Phi_{u_0}$ is also a contraction on the ball of radius $2C\|u_0\|_{H^{\sigma_{1}}}$
of $X^{\sigma_1}_{T}$. 
The fixed point of this contraction is the needed local solution of
(\ref{pak}).
Proposition~\ref{duhbis} also yields the bound
\begin{equation*}
\|\Phi_{u_0}(u)\|_{X^{\sigma}_{T}}\leq 
C\|u_0\|_{H^{\sigma}_{rad}(\Theta)}+CT^{\theta}(1+\|u\|_{X_{T}^{\sigma_{1}}}^{2})\|u\|_{X_{T}^{\sigma}}
\end{equation*}
We obtain thus the propagation of higher regularity with the same restrictions on $T$.
Using Proposition~\ref{duhbis}, we get the bound
\begin{equation}\label{contr2pak}
\|\Phi_{u_0}(u)-\Phi_{u_0}(v)\|_{X^{\sigma}_{T}}\leq 
CT^{\theta}\|u-v\|_{X_{T}^{\sigma}}(1+\|u\|_{X_{T}^{\sigma}}^{2}+\|v\|_{X_{T}^{\sigma}}^{2})\,.
\end{equation}
Applying Proposition~\ref{duhbis}, (\ref{contr2}) and (\ref{contr2pak}) to the difference of two solutions
yields the quoted Lipschitz bound.
This completes the proof of Proposition~\ref{lwp}.
\end{proof}
\begin{remarque}\label{granitza}
If we are interested to prove propagation of higher Sobolev regularity,
with our methods we only can treat the domains of the powers of the Dirichlet Laplacian. For
example  we may expect to get that $H^1_0$ regularity is propagated by the flow. Similar
results for the Neumann Laplacian do not follow from our analysis. As
mention in the introduction, we do not pursue this here since  the measure
$\rho$ "lives" on functions of Sobolev regularity $H^{1/2-}$.
\end{remarque}
We state the analog of Proposition~\ref{lwp} for (\ref{pakbis}).
\begin{proposition}\label{lwpbis}
Let us fix $0<\sigma_{1}\leq\sigma<1/2$.
Then there exist $b>1/2$, $\beta>0$, $C>0$, $\tilde{C}>0$, $c\in]0,1]$ such that for every $A>0$ if we set
$T= c(1+A)^{-\beta}$ then for every $N\geq 1$, every 
$u_0\in H^{\sigma_{1}}_{rad}(\Theta)$ satisfying $\|u_0\|_{H^{\sigma_{1}}}\leq A$ there exists a unique
solution $u=S_{N}(u)$ of (\ref{pakbis}) in $X^{\sigma_{1},b}_{rad}([-T,T]\times \Theta)$. Moreover
\begin{equation*}
\|u\|_{L^{\infty}([-T,T];H^{\sigma_{1}}(\Theta))}
\leq 
C\|u\|_{X^{\sigma_{1},b}_{rad}([-T,T]\times \Theta)}
\leq \tilde{C}\|u_0\|_{H^{\sigma_{1}}(\Theta)}\, .
\end{equation*}
If in addition $u_0\in H^{\sigma}_{rad}(\Theta)$ then
\begin{equation*}
\|u\|_{L^{\infty}([-T,T];H^{\sigma}(\Theta))}
\leq 
C\|u\|_{X^{\sigma,b}_{rad}([-T,T]\times \Theta)}
\leq \tilde{C}\|u_0\|_{H^{\sigma}(\Theta)}\, .
\end{equation*}
Finally if $u$ and $v$ are two solutions with data $u_0$, $v_0$ respectively,
satisfying 
$$
\|u_0\|_{H^{\sigma_{1}}}\leq A,\quad \|v_0\|_{H^{\sigma_{1}}}\leq A
$$
then 
$$
\|u-v\|_{L^{\infty}([-T,T];H^{\sigma_{1}}(\Theta))}\leq C\|u_0-v_0\|_{H^{\sigma_{1}}(\Theta)}\, .
$$ 
If in addition $u_0,v_0\in H^{\sigma}_{rad}(\Theta)$ then
$$
\|u-v\|_{L^{\infty}([-T,T];H^{\sigma}(\Theta))}\leq 
C
\|u_0-v_0\|_{H^{\sigma}(\Theta)}\, .
$$
\end{proposition}
\begin{proof}
The only new point compared to Proposition~\ref{lwp} is to observe that $S_N$ is bounded,
uniformly in $N$ on the Bourgain spaces
$X^{\sigma,b}_{rad}([-T,T]\times\Theta)$, namely 
for every $u\in X^{\sigma,b}_{rad}([-T,T]\times\Theta)$,
$$
\|S_{N}(u)\|_{X^{\sigma,b}_{rad}([-T,T]\times\Theta)}\leq\|u\|_{X^{\sigma,b}_{rad}([-T,T]\times\Theta)}\,,
$$
a bound which is direct consequence of the definitions of
$X^{\sigma,b}_{rad}([-T,T]\times\Theta)$ and the projector $S_{N}$.
\end{proof}
\begin{remarque}
The main point in Proposition~\ref{lwpbis} is the uniformness of the bounds
with respect to $N$.
\end{remarque}
\section{Improved bounds for the truncated NLS }
In this sections, we improve the result of Proposition~\ref{lwpbis}. More
precisely, we show bounds on the $H^{\sigma}$ norm of the solutions of
(\ref{pakbis}), uniform in $N$ for initial data of ``large $\rho_{N}$ measure''.
Let us denote by $\Phi_{N}(t)$ the smooth flow map of (\ref{pakbis}) 
which is defined globally thanks to the $L^2$ conservation for (\ref{pakbis}).
The next statement results from an application of Liouville's theorem to (\ref{pakbis}).
\begin{proposition}\label{liouville}
The measure $\rho_{N}$ defined in Section~3 is invariant under the flow of the (\ref{pakbis}).
\end{proposition}
\begin{proof}
Set $c(t)\equiv(c_{n}(t))_{1\leq n\leq N}$, where
$$
u(t)=\sum_{n=1}^{N}c_{n}(t)e_{n,s}\,.
$$
In the coordinates $c_{n}$, the equation (\ref{pakbis}) can be written as
\begin{equation}\label{4bis}
iz_{n}^{-s}\dot{c_n}(t)-z_n^2\,z_{n}^{-s} c_n(t)+
\int_{\Theta}S_{N}(F(u(t)))\overline{e_{n}}=0,\quad 1\leq n\leq N.
\end{equation}
Next, equation (\ref{4bis}) can be written in a Hamiltonian format as follows
$$
\partial_{t}c_n=-iz_{n}^{2s}\frac{\partial H}{\partial \overline{c_n}},
\quad 
\partial_{t}\overline{c_n}=iz_{n}^{2s}\frac{\partial H}{\partial c_n},\quad 1\leq n\leq N,
$$
with
$$
H(c,\overline{c})=\sum_{n= 1}^{N} z_n^{2-2s}\, |c_n|^2 - \int_{0}^1
V\Big(
\sum_{m=1}^{N}c_m\, e_{m,s}(r)\Big)rdr\, .
$$
Since
$$
\sum_{n=1}^{N}
\Big(
\frac{\partial}{\partial c_n}\big(-iz_{n}^{2s}\frac{\partial H}{\partial \overline{c_n}}\big)
+
\frac{\partial}{\partial\overline{c_n}}
\big(iz_{n}^{2s}\frac{\partial H}{\partial c_n}\big)
\Big)=0,
$$
we can apply the Liouville theorem for divergence free vector fields to
conclude that the measure $dcd\overline{c}$ is invariant under the flow of
(\ref{pakbis}).
Recall that we denote by $\Phi_{N}(t)$, the flow of (\ref{pakbis}) and that the quantities $H(c,\overline{c})$ and 
$$
\|c\|^{2}\equiv\sum_{n=1}^{N}|c_n|^{2}
$$
are conserved under $\Phi_{N}(t)$.
Let $A$ be a Borel set of $E_{N}$. 
Recall that we denote by $\chi$ the characteristic function of the interval
$[0,R]$, $R>0$. Then 
$$
\rho_{N}(A)=\kappa_{N}\int_{A}e^{-H(c,\overline{c})}\chi(\|c\|)dcd\overline{c},
$$
where
$$
\kappa_{N}=\pi^{-N}\Big(\prod_{1\leq n\leq N}z_{n}^{2-2s}\Big).
$$
In addition
\begin{equation}\label{jacobi}
\rho_{N}\Big(\Phi_{N}(t)(A)\Big)=\kappa_{N}\int_{\Phi_{N}(t)(A)}e^{-H(c,\overline{c})}\chi(\|c\|)dcd\overline{c}.
\end{equation}
We can write
$$
\Phi_{N}(t)(A)=\big\{
(c,\overline{c})\,:\,(c,\overline{c})=\Phi_{N}(t)(b,\overline{b}),\quad
(b,\overline{b})\in A
\big\}.
$$
Let us perform the change of variables
$
(c,\overline{c})=\Phi_{N}(t)(b,\overline{b})
$
in the right hand-side of (\ref{jacobi}). Since $dcd\overline{c}$ is invariant
under $\Phi_{N}(t)$ the Jacobian of this variable change is one. Next by the
conservation laws
$$
H(\Phi_{N}(t)(b,\overline{b}))=H(b,\overline{b}),\quad
\|\Phi_{N}(t)(b)\|=\|b\|.
$$
Therefore
$$
\rho_{N}\Big(\Phi_{N}(t)(A)\Big)=\kappa_{N }\int_{A}e^{-H(b,\overline{b})}\chi(\|b\|)dbd\overline{b}=\rho_{N}(A).
$$
This completes the proof of Proposition~\ref{liouville}.
\end{proof}
Next, we state a bound for the solutions of (\ref{pakbis}) which gives a
control, independent of $N$ on norms which are stronger
then $L^2$ but weaker then $H^1$. 
\begin{proposition}\label{longtime}
For every integer $i\geq 1$, $\sigma\in[s,1/2[$ there exists a set
$$
\Sigma_{N,\sigma}^{i}\subset E_{N}
$$
such that
\begin{equation}\label{poslednopak}
\rho_{N}(E_{N}\backslash \Sigma_{N,\sigma}^{i})\leq 2^{-i},
\end{equation}
and for $u_0\in \Sigma_{N,\sigma}^{i}$ one has the bound
\begin{equation}\label{log}
\|\Phi_{N}(t)(u_0)\|_{H^{\sigma}}\leq C_{\sigma}(i+\log(1+|t|))^{\frac{1}{2}}\,.
\end{equation}
Moreover, for $N_1\leq N_2$ we have the inclusion $\Sigma_{N_1,\sigma}^{i}\subset \Sigma^{i}_{N_2,\sigma}$.
\end{proposition}
\begin{proof}
We will consider only the positive values of $t$, the analysis for $t<0$ being
the same.
For $\sigma\in [s,1/2[$, $i,j$ integers $\geq 1$, we set
$$
B_{N,\sigma}^{i,j}(D_{\sigma})=\Big\{u\in E_{N}\,:\,\|u\|_{H^{\sigma}(\Theta)}\leq D_{\sigma}(i+j)^{\frac{1}{2}},\quad\|u\|_{L^2(\Theta)}\leq R\Big\},
$$
where the number $D_{\sigma}\gg 1$ (independent of $i,j,N$) will be fixed later. 
Thanks to Proposition~\ref{lwpbis}, there exists $\tau\in ]0,1]$, $\tau \sim D_{\sigma}^{-\beta}(i+j)^{-\beta/2}$
for some $\beta>0$ and such that for $t\in[0,\tau]$,
\begin{equation}\label{preser}
\Phi_{N}(t)\big(B_{N,\sigma}^{i,j}(D_{\sigma})\big)\subset B_{N,\sigma}^{i,j}(CD_{\sigma})\, ,
\end{equation}
where $B_{N,\sigma}^{i,j}(CD_{\sigma})$ is defined similarly to
$B_{N,\sigma}^{i,j}(D_{\sigma})$ simply replacing $D_{\sigma}$ by
$CD_{\sigma}$ in the $H^{\sigma}$ bound for $u$.
Next, we set
$$
\Sigma_{N,\sigma}^{i,j}(D_{\sigma})=\bigcap_{k=0}^{[2^{j}/\tau]}\Phi_{N}(-k\tau)(B_{N,\sigma}^{i,j}(D_{\sigma}))\, ,
$$
where $[2^{j}/\tau]$ stays for the integer part of $2^{j}/\tau$. Using Proposition~\ref{liouville}, we can write
\begin{eqnarray*}
\rho_{N}(E_{N}\backslash\Sigma_{N,\sigma}^{i,j}(D_{\sigma}))
& = &
\rho_{N}\Big(\bigcup_{k=0} ^{[2^{j}/\tau]}(E_{N}\backslash\Phi_{N}(-k\tau)(B_{N,\sigma}^{i,j}(D_{\sigma})))\Big) 
\\
& \leq &
([2^{j}/\tau]+1)\rho_{N}(E_N\backslash B_{N,\sigma}^{i,j}(D_{\sigma}))
\\
& \leq &
C2^{j}D_{\sigma}^{\beta}(i+j)^{\beta/2}\rho_{N}(E_N\backslash B_{N,\sigma}^{i,j}(D_{\sigma}))\,.
\end{eqnarray*}
Let us now observe that 
\begin{eqnarray*}
\rho_{N}(E_N\backslash B_{N,\sigma}^{i,j}(D_{\sigma})) & = &
\rho
\Big(u\in H^{s}_{rad}(\Theta)\,:\,\|S_{N}(u)\|_{H^{\sigma}(\Theta)}> D_{\sigma}(i+j)^{\frac{1}{2}} 
\Big)
\\
& \leq &
\rho\Big(u\in H^{s}_{rad}(\Theta)\,:\,\|u\|_{H^{\sigma}(\Theta)}> D_{\sigma}(i+j)^{\frac{1}{2}} 
\Big)\,.
\end{eqnarray*}
Therefore, using Proposition~\ref{high-tris}, we can write
\begin{equation}\label{zvez}
\rho_{N}(E_{N}\backslash\Sigma_{N,\sigma}^{i,j}(D_{\sigma}))\leq
C_{\sigma}2^{j}D_{\sigma}^{\beta}(i+j)^{\beta/2}e^{-cD_{\sigma}^2(i+j)}\leq 2^{-(i+j)},
\end{equation}
provided $D_{\sigma}\gg 1$, depending on $\sigma$ but independent of $i,j,N$.
Thanks to (\ref{preser}), we obtain that for
$u_0\in\Sigma_{N,\sigma}^{i,j}(D_\sigma)$, the solution $u$ of (\ref{pakbis}) with data $u_0$ satisfies
\begin{equation}\label{jjj1}
\|u(t)\|_{H^{\sigma}(\Theta)}\leq CD_{\sigma}(i+j)^{\frac{1}{2}},\quad 0\leq t\leq 2^{j}\,.
\end{equation}
Next, we set
$$
\Sigma_{N,\sigma}^{i}=\bigcap_{j\geq 1}\Sigma_{N,\sigma}^{i,j}(D_{\sigma})\,.
$$
Thanks (\ref{zvez}),
\begin{equation}\label{posledno}
\rho_{N}(E_{N}\backslash \Sigma_{N,\sigma}^{i})\leq 2^{-i}\,.
\end{equation}
Next, using (\ref{jjj1}), we get (\ref{log}).
Observe that for $N_1\leq N_2$, we have the inclusion 
$B_{N_1,\sigma}^{i,j}(D_{\sigma})\subset B_{N_2,\sigma}^{i,j}(D_{\sigma})$ which implies that
$\Sigma_{N_1,\sigma}^{i,j}(D_{\sigma})\subset \Sigma_{N_2,\sigma}^{i,j}(D_{\sigma})$.
This in turn implies that for $N_1\leq N_2$, $\Sigma_{N_1,\sigma}^{i}\subset\Sigma_{N_2,\sigma}^{i}$.
This completes the proof of Proposition~\ref{longtime}.
\end{proof}
Next, we prove an invariance property of the sets $\Sigma^{i}_{N,\sigma}$
constructed in  Proposition~\ref{longtime}.
\begin{proposition}\label{longtimepak}
For every $\sigma\in ]s,1/2[$ every $\sigma_1\in [s,\sigma[$ every $t\in\R$
every integer $i\geq 1$ there exists $i_{1}$ such that for every $N\geq 1$, if $u_0\in
\Sigma^{i}_{N,\sigma}$ then one has
$$
\Phi_{N}(t)(u_0)\in \Sigma^{i+i_1}_{N,\sigma_1}.
$$
\end{proposition}
\begin{proof}
Again, we can suppose that $t>0$.
Set $u(t)\equiv\Phi_{N}(t)(u_0)$. 
If $u_{0}\in \Sigma_{N,\sigma}^{i}$ then for every integer
$j\geq 1$, we have the bound
$$
\|\Phi_{N}(t_1)(u_0)\|_{H^{\sigma}}\leq
C_{\sigma}(i+j)^{\frac{1}{2}},\quad 0\leq t_1\leq 2^j.
$$
Let $j_0\in\N$ (depending on $t$) be such that for every $j\geq 1$, $2^{j}+t\leq 2^{j+j_0}$.
Therefore, we have that
$$
\|\Phi_{N}(t_1)(u(t))\|_{H^{\sigma}}
=
\|\Phi_{N}(t+t_1)(u_0)\|_{H^{\sigma}}
\leq
C_{\sigma}(i+j+j_0)^{\frac{1}{2}},\quad 0\leq t_1\leq 2^j.
$$
The crucial observation is that thanks to the $L^2$ conservation law,
interpolating between the last bound and the $L^2$ conservation provides the
existence of $\theta\in]0,1[$ (depending on $\sigma$ and $\sigma_1$) such that
$$
\|\Phi_{N}(t_1)(u(t))\|_{H^{\sigma_1}}
\leq
c\Big[C_{\sigma}(i+j+j_0)\Big]^{\frac{\theta}{2}},\quad 0\leq t_1\leq 2^j.
$$
Next, we observe that since $\theta<1$, for $j_0\gg 1$,
$$
c\Big[C_{\sigma}(i+j+j_0)\Big]^{\frac{\theta}{2}}\leq
D_{\sigma_1}(i+j+j_0)^{\frac{1}{2}}\,.
$$
Thus
$$
\|\Phi_{N}(t_1)(u(t))\|_{H^{\sigma_1}}
\leq
D_{\sigma_1}(i+j+j_0)^{\frac{1}{2}}, \quad 0\leq t_1\leq 2^j.
$$
We can now conclude that $u(t)\in
\Sigma_{N,\sigma_1}^{i+j_0,j}(D_{\sigma})$ for every $j\geq
1$. Therefore 
$$
u(t)\in \Sigma_{N,\sigma_1}^{i+j_0}\,.
$$
This completes the proof of Proposition~\ref{longtimepak}.
\end{proof}
\begin{remarque}
The number $i_1$ is the same for every $i$, i.e. it depends only on
$t,\sigma,\sigma_1$.
This fact is however not of importance for the sequel.
\end{remarque}
\section{Global existence for NLS on a set of full $\rho$ measure }
The goal of this section is to compare the
flows of (\ref{pak}) and (\ref{pakbis}) on a set of full $\rho$ measure.
For an integer $i\geq 1$ and $\sigma\in[s,1/2[$, we set
$$
\Sigma_{\sigma}^{i}\equiv\bigcup_{N\geq 1}\Sigma_{N,\sigma}^{i}.
$$
where $\Sigma_{N,\sigma}^{i}$ are defined in Proposition~\ref{longtime}.
Let us denote by $\overline{\Sigma_{\sigma}^{i}}$ the closure of $\Sigma_{\sigma}^{i}$ in
$H^{\sigma}_{rad}(\Theta)$. Thus $\Sigma_{\sigma}^{i}$ is a closed set of $H^{\sigma}_{rad}(\Theta)$.
Then thanks to Lemma~\ref{lim2} and Proposition~\ref{longtime}, we can write
\begin{equation}\label{perm}
\rho(\overline{\Sigma_{\sigma}^{i}})\geq\limsup_{N\rightarrow \infty}\rho_{N}(\Sigma_{N,\sigma}^{i})
\geq 
\limsup_{N\rightarrow \infty}\big(\rho_{N}(E_{N})-2^{-i}\big)
=
\rho\big(H^s_{rad}(\Theta)\big)-2^{-i}.
\end{equation}
Next, we set
$$
\Sigma_{\sigma}\equiv\bigcup_{i\geq 1}\overline{\Sigma_{\sigma}^{i}}\,.
$$
In view of (\ref{perm}), we obtain that $\Sigma_{\sigma}$ is of full $\rho$
measure.
\\

Let $ l=(l_{j})_{j\in\N}$ be a increasing sequence of real numbers
such that $l_0=s$, $l_{j}<1/2$ and 
$$
\lim_{j\rightarrow\infty}l_{j}=1/2.
$$
Then, we set
\begin{equation}\label{intersection}
\Sigma=\bigcap_{\sigma\in l}\Sigma_{\sigma}
\end{equation}
The set $\Sigma$ is of full $\rho$ measure since every
$\Sigma_{\sigma}$ is of full $\rho$ measure and the
intersection in (\ref{intersection}) is countable.
The set $\Sigma$ is the statistical ensemble for the problem (\ref{pak}) and
the solutions of (\ref{pak}) with data in $\Sigma$ are globally defined. 
We have the following statement.
\begin{proposition}\label{compare}
For every $u_0\in\Sigma$, the local solution of (\ref{pak})
given by Proposition~\ref{lwp} is globally defined.
Moreover for every $t\in\R$, if we denote by $\Phi(t)$ the flow map of
(\ref{pak}) acting on $\Sigma$ then $\Phi(t)(\Sigma)=\Sigma$.
\end{proposition}
\begin{proof}
Let us fix $u_0\in\overline{\Sigma_{\sigma}^{i}}$, $\sigma\in l$, $\sigma_1\in]0,\sigma[$ and $T>0$.
Thus there exists a sequence $u_{0,k}\in\Sigma^{i}_{N_k,\sigma}$, where $N_k$
is tending to infinity, such that
$u_{0,k}$ converges to $u_0$ in $H^{\sigma}(\Theta)$. 
Thanks to Proposition~\ref{longtime}
\begin{equation}\label{ant}
\|\Phi_{N_k}(t)(u_{0,k})\|_{H^{\sigma}}\leq C_{\sigma}(i+\log(1+|t|))^{\frac{1}{2}}\,.
\end{equation}
Set 
$$
u_{N_k}(t)\equiv  \Phi_{N_k}(t)(u_{0,k})\,.
$$
Thanks to (\ref{ant}), there exists $\Lambda>1$, independent of $N_k$, such that
\begin{equation}\label{david1}
\|u_{N_k}(t)\|_{H^{\sigma}}\leq \Lambda,\quad |t|\leq T.
\end{equation}
Let us observe that (\ref{david1}), applied for $t=0$ implies that
$\|u_0\|_{H^{\sigma}}\leq\Lambda$ (after passing to the limit $N_k\rightarrow\infty$).
Let $\tau>0$ be the local existence time for (\ref{pak}), provided by
Proposition~\ref{lwp} for $A=\Lambda+1$. Recall that we can assume $\tau\sim\Lambda^{-\beta}$
for some $\beta>0$. Denote by $u(t)$ the solution of (\ref{pak}) with data
$u_0$ on the time interval $[-\tau,\tau]$. Set 
$$
v\equiv u-u_{N_k}.
$$
Then $v$ solves the equation
\begin{equation}\label{eqnv}
iv_t+\Delta v +F(u)-S_{N_k}(F(u_{N_k}))=0, \quad v|_{t=0}=u_0-u_{0,k} \, .
\end{equation}
Next we write
$$
F(u)-S_{N_k}(F(u_{N_k}))=S_{N_k}\big(F(u)-F(u_{N_k})\big)+(1-S_{N_k})F(u).
$$
Observe that the map $1-S_{N}$ sends $H^{\sigma}_{rad}(\Theta)$ to
$H^{\sigma_{1}}_{rad}(\Theta)$ with norm $\leq CN^{\sigma_{1}-\sigma}$.
Similarly, for $I\subset\R$ an interval, the map
$1-S_{N}$ sends $X^{\sigma,b}_{rad}(I\times\Theta)$ to
$X^{\sigma_{1},b}_{rad}(I\times\Theta)$ with norm $\leq CN^{\sigma_{1}-\sigma}$.
Moreover $S_{N}$ acts as a bounded operator (with norm $\leq 1$) on the Bourgain spaces
$X^{\sigma,b}_{rad}$. Therefore, using Proposition~\ref{duhbis}, we can write the
Duhamel formula associated to (\ref{eqnv}) and we obtain that there exists
$b>1/2$ and $\theta>0$ (depending only on $\sigma$, $\sigma_1$) such that one has the bound
\begin{eqnarray*}
\|v\|_{X^{\sigma_1,b}_{rad}([-\tau,\tau]\times\Theta)}
& \leq &
C\|u_0-u_{0,k}\|_{H^{\sigma_1}(\Theta)}
\\
&  & + C\tau^{\theta}\|v\|_{X^{\sigma_1,b}_{rad}([-\tau,\tau]\times\Theta)}
\big(
1+\|u\|_{X^{\sigma_1,b}_{rad}([-\tau,\tau]\times\Theta)}^{2}+
\\
& &
+\|u_{N_k}\|_{X^{\sigma_1,b}_{rad}([-\tau,\tau]\times\Theta)}^{2}
\big)
\\
& & + C\tau^{\theta}N_k^{\sigma_1-\sigma}\|u\|_{X^{\sigma,b}_{rad}([-\tau,\tau]\times\Theta)}
\big(
1+\|u\|_{X^{\sigma_1,b}_{rad}([-\tau,\tau]\times\Theta)}^{2}
\big).
\end{eqnarray*}
Using Proposition~\ref{lwp} and Proposition~\ref{lwpbis}, we get
\begin{eqnarray*}
\|v\|_{X^{\sigma_1,b}_{rad}([-\tau,\tau]\times\Theta)}
& \leq &
C\|u_0-u_{0,k}\|_{H^{\sigma_1}(\Theta)}
\\
& &
+
C\tau^{\theta}\|v\|_{X^{\sigma_1,b}_{rad}([-\tau,\tau]\times\Theta)}
(1+C\|u_{0}\|_{H^{\sigma_1}(\Theta)}^{2}+C\|u_{0,k}\|_{H^{\sigma_1}(\Theta)}^{2})
\\
& & +
C\tau^{\theta}N_k^{\sigma_1-\sigma}\|u_0\|_{H^{\sigma}(\Theta)}(1+C\|u_{0}\|_{H^{\sigma_1}(\Theta)}^{2})
\\
& \leq &
C\|u_0-u_{0,k}\|_{H^{\sigma_1}(\Theta)}+
C\tau^{\theta}\Lambda^{2}N_k^{\sigma_1-\sigma}\|u_0\|_{H^{\sigma}(\Theta)}
\\
& &
+C\tau^{\theta}\Lambda^{2}\|v\|_{X^{\sigma_1,b}_{rad}([-\tau,\tau]\times\Theta)}.
\end{eqnarray*}
Therefore, assuming in addition that $\tau\sim \Lambda^{-\theta/2}$, we obtain
$$
\|v\|_{X^{\sigma_{1},b}_{rad}([-\tau,\tau]\times\Theta)}\leq 
C\|u_0-u_{0,k}\|_{H^{\sigma_1}(\Theta)}+
CN_k^{\sigma_1-\sigma}\|u_0\|_{H^{\sigma}(\Theta)},\quad \tau\sim\Lambda^{-\beta},
$$
for some fixed positive real number $\beta$ and where the constant $C$ is independent of $N_k$.
Since $b>1/2$, the last inequality implies
\begin{equation}\label{vvvv}
\|v(t)\|_{H^{\sigma_1}(\Theta)}\leq 
C\|u_0-u_{0,k}\|_{H^{\sigma_1}(\Theta)}+
CN_k^{\sigma_1-\sigma}\|u_0\|_{H^{\sigma}(\Theta)},
\end{equation}
where $|t|\leq \tau\sim\Lambda^{-\beta}$, $\beta>0$.
By taking $N_k\gg 1$ and using the triangle inequality, we get
\begin{equation}\label{david2}
\|u(t)\|_{H^{\sigma_1}(\Theta)}\leq\Lambda+1,\quad |t|\leq\tau.
\end{equation}
The key quantity in this discussion is
$$
\|v(t)\|_{H^{\sigma_1}(\Theta)}+N_{k}^{\sigma_1-\sigma}\|u(t)\|_{H^{\sigma}(\Theta)}\,.
$$
We can iterate the argument for obtaining (\ref{vvvv}) on $[\tau,2\tau]$ thanks to the definition of
$\tau$ and the bounds (\ref{david1}) and (\ref{david2}). 
We obtain
\begin{equation*}
\|v(t)\|_{H^{\sigma_1}(\Theta)}\leq C\|v(\tau)\|_{H^{\sigma_1}(\Theta)}+CN_k^{\sigma_1-\sigma}\|u(\tau)\|_{H^{\sigma}(\Theta)},
\end{equation*}
where $t\in[\tau,2\tau]$ and  $\tau \sim\Lambda^{-\beta}$. Moreover, by  taking $N_k\gg 1$,
\begin{equation*}
\|u(t)\|_{H^{\sigma_{1}}(\Theta)}\leq\Lambda+1,\quad \tau \leq t\leq 2\tau.
\end{equation*}
Then, we can continue
by covering the interval $[-T,T]$ with intervals of size $\tau$, which yields
the existence of $u(t)$ on $[-T,T]$. 
Moreover $v$ satisfies the bound 
\begin{equation*}
\|v(t)\|_{H^{\sigma_{1}}(\Theta)}\leq 
C^{1+T}\Big(
N_k^{\sigma_1-\sigma}\|u_0\|_{H^{\sigma}(\Theta)}
+
\|u_0-u_{0,k}\|_{H^{\sigma_1}(\Theta)}
\Big),
\quad |t|\leq T.
\end{equation*}
Therefore by taking $N_k\gg 1$ (depending in particular on $T$), we
obtain that for every $\varepsilon>0$ there exists $N_{0}$ such that for
$N_k\geq N_0$ one has the inequality 
\begin{equation*}
\sup_{|t|\leq T}\|u(t)-\Phi_{N_k}(t)(u_{0,k})\|_{H^{\sigma_1}(\Theta)}<\varepsilon\,.
\end{equation*}
Hence for every $t\in[-T,T]$,
\begin{equation}\label{oliv}
\lim_{k\rightarrow\infty}
\|u(t)-\Phi_{N_k}(t)(u_{0,k})\|_{H^{\sigma_1}(\Theta)}=0\,.
\end{equation}
Since $T>0$ was chosen arbitrary, we obtain that for every
$u_0\in\overline{\Sigma_{\sigma}^{i}}$ the local solution of (\ref{pak}) is globally
defined. Since $i$ and $\sigma\in l$ are also arbitrary, we obtain that for every
$u_0\in\Sigma$, the the local solution of (\ref{pak}) is globally
defined. Let us denote by $\Phi(t)$ the flow of (\ref{pak}) acting on $\Sigma$.
Let us show the inclusion
\begin{equation}\label{inclu}
\Phi(t)(\Sigma)\subset\Sigma.
\end{equation}
Fix $u_0\in\Sigma$. It suffices to show that for every $\sigma_1\in l$, 
we have 
$$
\Phi(t)(u_0)\in \Sigma_{\sigma_1}\,.
$$
Let us take $\sigma\in ]\sigma_1,1/2[$, $\sigma\in l$. Since $u_0\in\Sigma$, we have that
$u_0\in \Sigma_{\sigma}$. 
Therefore there exists $i$ such that 
$
u_0\in \overline{\Sigma^{i}_{\sigma}}.
$
Let again $u_{0,k}\in\Sigma^{i}_{N_k,\sigma}$ be a sequence which tends to
$u_0$ in $H^{\sigma}(\Theta)$.
Thanks to Proposition~\ref{longtimepak} there exists $i_1$ such that
$$
\Phi_{N_k}(t)(u_{0,k})\in\Sigma^{i+i_1}_{N_k,\sigma_1}\,.
$$
Therefore using (\ref{oliv}), we obtain that
$$
\Phi(t)(u_0)\in \overline{\Sigma^{i+i_1}_{\sigma_1}}.
$$
Thus $\Phi(t)(u_0)\in\Sigma_{\sigma_1}$ which proves (\ref{inclu}).
Moreover the flow $\Phi(t)$ is reversible which
implies that 
$
\Phi(t)(\Sigma)=\Sigma.
$
Indeed, if $u\in\Sigma$ and $t\in\R$, we set $u_0\equiv\Phi(-t)u\in\Sigma$
(which is well-defined thanks to the previous analysis) and thus
$u=\Phi(t)u_{0}$,
i.e. $\Sigma\subset\Phi(t)(\Sigma)$.
This completes the proof of Proposition~\ref{compare}.
\end{proof}
We complete this section by getting a continuity property of $\Phi(t)$.
\begin{proposition}\label{gauss}
Let $u\in\Sigma$ and $u_n\in\Sigma$ be a sequence such that $u_n\rightarrow u$
in $H^s(\Theta)$. Then for every $t\in\R$,
$
\Phi(t)(u_n)\rightarrow \Phi(t)(u)
$
in $H^s(\Theta)$. In particular, for every $A$, a closed set in
$H^s_{rad}(\Theta)$ one has
$$
\Phi(t)(A\cap \Sigma)=\overline{\Phi(t)(A\cap \Sigma)}\cap\Sigma,
$$
where $\overline{\Phi(t)(A\cap \Sigma)}$ denotes the closure in
$H^s_{rad}(\Theta)$ of $\Phi(t)(A\cap \Sigma)$.
\end{proposition}
\begin{proof}
Since $u\in\Sigma$ there
exists $\Lambda\geq 1$ such that
$$
\sup_{|\tau|\leq |t|}\|\Phi(\tau)(u)\|_{H^s(\Theta)}\leq \Lambda.
$$
Let us denote by $\tau_0$ the local existence time in Proposition~\ref{lwp},
associated to $A=2\Lambda$. Then, by the continuity of the flow
$\Phi(\tau_0)(u_n)\rightarrow \Phi(\tau_0)(u)$ in $H^s(\Theta)$. Next, we cover
the interval $[0,t]$ by intervals of size $\tau_0$ and we apply the continuity
of the flow established in  Proposition~\ref{lwp} at each step. 
Therefore, we obtain that $\Phi(t)(u_n)\rightarrow \Phi(t)(u)$ in $H^s(\Theta)$.
Since $\Phi(t)(\Sigma)\subset\Sigma$, it is clear that
\begin{equation}\label{elen}
\Phi(t)(A\cap \Sigma)\subset\overline{\Phi(t)(A\cap \Sigma)}\cap\Sigma.
\end{equation}
Next, let us fix $u\in \overline{\Phi(t)(A\cap \Sigma)}\cap\Sigma$.
Thus there exists $v_n\in A\cap\Sigma$ such that $u_n\equiv\Phi(t)(v_n)$
converges to $u$ in $H^s(\Theta)$. Since $v_n\in\Sigma$ and $u\in\Sigma$, we
obtain that $u_n\in\Sigma$ and $\Phi(-t)(u)\in\Sigma$. Therefore, using the
continuity property we have just established, we obtain that
$v_{n}=\Phi(-t)(u_n)$ converges to $\Phi(-t)(u)$ in $H^s(\Theta)$. Since the
set $A$ is assumed closed, we obtain that $\Phi(-t)(u)\in A$. Thus
$u\in\Phi(t)(A\cap\Sigma)$ which gives the opposite to (\ref{elen}) inclusion.
This completes the proof of Proposition~\ref{gauss}.
\end{proof}
\section{Invariance of the measure $\rho$ }
In this section, we complete the proof of Theorem~\ref{glavna}.
Recall that we denote by $\Phi(t)$, $t\in\R$ the flow of (\ref{pak}) acting on $\Sigma$,
defined in (\ref{intersection}). 
Thanks to the continuity properties of $\Phi(t)$ displayed by Proposition~\ref{gauss}, we have that if
$A\subset\Sigma$ is a $\rho$ measurable set then so is $\Phi(t)(A)$.
Let us observe that thanks to the reversibility of the flow
$\Phi(t)$, it suffices to prove that for every $t\in\R$ and every $\rho$
measurable set $A\subset\Sigma$ one has the inequality
\begin{equation}\label{vili}
\rho\big(\Phi(t)(A)\big)\geq \rho(A).
\end{equation}
Let us show that it suffice to prove (\ref{vili}) only for closed sets of $H^s_{rad}(\Theta)$.
Indeed, by the regularity of the bounded Borel measures for every $\rho$ measurable
set $A\subset\Sigma$, 
we can find a sequence of closed sets $F_{n}\subset A$ such that 
$$
\rho(A)=\lim_{n\rightarrow\infty }\rho(F_{n})\,.
$$ 
Hence if we can prove (\ref{vili}) for the sets $F_n$ then we can write 
$$
\rho(A)=
\lim_{n\rightarrow\infty}\rho(F_n)\leq
\limsup_{n\rightarrow\infty}\rho\big(\Phi(t)(F_{n})\big)\leq \rho\big(\Phi(t)(A)\big).
$$
Therefore, it suffices to prove (\ref{vili}) for closed sets of
$H^s_{rad}(\Theta)$ which are included in $\Sigma$.
\\

Fix $\sigma\in]s,1/2[$, $\sigma\in l$.
Let us next show that it suffices to prove (\ref{vili}) for subsets of
$\Sigma$ which are compacts of $H^s_{rad}(\Theta)$ which are bounded in $H^{\sigma}_{rad}(\Theta)$.
Indeed, using Lemma~\ref{lim3}, we
can write that for every closed in $H^{s}_{rad}(\Theta)$ set $A\subset\Sigma$, one has
$$
\rho(A)=\lim_{R\rightarrow\infty}\rho(A\cap K_{R}),
$$
where $K_{R}$ is the closed ball of radius $R$ in $H^{\sigma}_{rad}(\Theta)$,
$\sigma\in ]s,1/2[$. Thus $A\cap K_{R}$ is a compact in $H^{s}_{rad}(\Theta)$ 
and if we can prove (\ref{vili}) for compacts  which are bounded in
$H^{\sigma}_{rad}(\Theta)$ then
$$
\rho(A)
\leq
\limsup_{R\rightarrow\infty}
\rho\Big(
\Phi(t)(A\cap K_{R})
\Big)\leq \rho(\Phi(t)(A)).
$$
Thus, it suffices to prove (\ref{vili}) for subsets of $\Sigma$
which are compacts in $H^{s}_{rad}(\Theta)$ and bounded in $H^{\sigma}_{rad}(\Theta)$.
\\

Let us now fix $t\in\R$ and $K\subset\Sigma$, a bounded set of
$H^{\sigma}_{rad}(\Theta)$ which is a compact in $H^{s}_{rad}(\Theta)$. 
Then we have the following lemma.
\begin{lemme}\label{yaki}
There exists a  ball ${\mathcal B}$, centered at the origin of 
$H^{s}_{rad}(\Theta)$ containing all $\Phi(\tau)(K)$, $|\tau|\leq |t|$.
\end{lemme}
\begin{proof}
The sets $\Phi(\tau)(K)$ are contained in a ball of $H^{s}_{rad}(\Theta)$ for
$|\tau|$ small enough, given by Proposition~\ref{lwp}.
We then argue by contradiction by supposing that there exists $T$ and a
sequence $u_n\in K$ such that
\begin{equation}\label{yakibis}
\lim_{n\rightarrow \infty}\|\Phi(T)(u_n)\|_{H^s(\Theta)}=\infty\,.
\end{equation}
Since $K$ is a compact, there exists a subsequence still denoted by $u_n$ and
$u\in K$ such that $u_n \rightarrow u$ in $H^{s}_{rad}(\Theta)$. 
Since $u\in\Sigma$, we can apply Proposition~\ref{gauss} and we
obtain that 
$\Phi(T)(u_n)\rightarrow \Phi(T)(u)$ in $H^s(\Theta)$ which contradicts
(\ref{yakibis}).
This completes the proof of Lemma~\ref{yaki}.
\end{proof}
Let us denote by $R_1$ the radius of ${\mathcal B}$. Set
$$
\tau_1\equiv c(1+R_1)^{-M}\, ,
$$
where $0<c\ll 1$ and $M\gg 1$ are two parameters to be fixed later.
A first restriction on $c$ and $M$ is to chose them so that $\tau_1$ is
smaller than the time existence provided by
Propositions~\ref{lwp},\ref{lwpbis} associated to $A=R_1$ (and $\sigma_1=s$).
It is then sufficient to prove that
\begin{equation}\label{restriction}
\rho(K)\leq\rho\big(\Phi(\tau)(K)\big),\quad |\tau|\leq \tau_1\,.
\end{equation}
Indeed, once (\ref{restriction}) is established, it suffices to cover $[0,t]$
by intervals of size $\sim\tau_1$ and to apply (\ref{restriction}) at each
step. Such an iteration is possible since at each step the image under
$\Phi(\tau)$ of the corresponding set remains in ${\mathcal B}$ and is
included in $\Sigma$.
\\

Let us now prove (\ref{restriction}). Fix $ \varepsilon>0$.
Denote by $B_{\varepsilon}$ the open ball centered at the origin and of radius
$\varepsilon$ of $H^{s}_{rad}(\Theta)$. 
Recall that we denote by $\Phi_{N}(t)$, $t\in\R$ the flow of (\ref{pakbis}).
Then using Proposition~\ref{lwpbis}, we infer that there exists $c>0$ such that
\begin{equation}\label{kenji}
\Phi_{N}(\tau)\Big(\big(K+B_{\varepsilon}\big)\cap E_{N}\Big)\subset
\Phi_{N}(\tau)(S_{N}(K))+B_{c\varepsilon},\quad N\gg 1.
\end{equation}
We now make appeal to the following lemma.
\begin{lemme}\label{kyoto}
For $N\gg 1$ one has the inclusion
$$
\Phi_{N}(\tau)(S_{N}(K))+B_{c\varepsilon}\subset
\Phi(\tau)(K)+\overline{B_{2c\varepsilon}}\,.
$$
\end{lemme}
\begin{proof}
The argument is similar to the proof of Proposition~\ref{compare}.
For $u_0\in K$, we denote by $u$ the solution of (\ref{pak}) with data $u_0$
and by $u_{N}$ the solution of (\ref{pakbis}) with data $S_{N}(u_0)$. Next, we
set $v\equiv u-u_N$. Then $v$ is a solution of
\begin{equation}\label{eqnvpak}
iv_t+\Delta v +F(u)-S_{N}(F(u_{N}))=0, \quad v|_{t=0}=(1-S_{N})u_0\,.
\end{equation}
By writing
$$
F(u)-S_{N}(F(u_{N}))=S_{N}\big(F(u)-F(u_{N})\big)+(1-S_{N})F(u)
$$
and using Proposition~\ref{duhbis}, we obtain that there exists
$b>1/2$ and $\theta>0$ such that one has
\begin{eqnarray*}
\|v\|_{X^{s,b}_{rad}([-\tau,\tau]\times\Theta)}
& \leq &
CN^{s-\sigma}\|u_0\|_{H^{\sigma}(\Theta)}
\\
&  & + C\tau^{\theta}\|v\|_{X^{s,b}_{rad}([-\tau,\tau]\times\Theta)}
\big(
1+\|u\|_{X^{s,b}_{rad}([-\tau,\tau]\times\Theta)}^{2}
\\
& &
+\|u_{N}\|_{X^{s,b}_{rad}([-\tau,\tau]\times\Theta)}^{2}
\big)
\\
& & + C\tau^{\theta}N^{s-\sigma}\|u\|_{X^{\sigma,b}_{rad}([-\tau,\tau]\times\Theta)}
\big(
1+\|u\|_{X^{s,b}_{rad}([-\tau,\tau]\times\Theta)}^{2}
\big).
\end{eqnarray*}
Using Proposition~\ref{lwp} and Proposition~\ref{lwpbis}, we get
\begin{eqnarray*}
\|v\|_{X^{s,b}_{rad}([-\tau,\tau]\times\Theta)}
& \leq &
CN^{s-\sigma}\|u_0\|_{H^{\sigma}(\Theta)}
\\
& &
+
C\tau^{\theta}\|v\|_{X^{s,b}_{rad}([-\tau,\tau]\times\Theta)}
(1+C\|u_{0}\|_{H^{s}(\Theta)}^{2})
\\
& & +
C\tau^{\theta}N^{s-\sigma}\|u_0\|_{H^{\sigma}(\Theta)}(1+C\|u_{0}\|_{H^{s}(\Theta)}^{2})\,.
\end{eqnarray*}
Coming back to the definition of $\tau_1$, by taking $c\ll 1$ and $M\gg 1$, we infer that
$$
\|v\|_{X^{s,b}_{rad}([-\tau,\tau]\times\Theta)}\leq 
CN^{s-\sigma}\|u_0\|_{H^{\sigma}(\Theta)}.
$$
Using that $u_0$ is in a bounded set of $H^{\sigma}_{rad}(\Theta)$ and since $b>1/2$, the last inequality implies
\begin{equation*}
\|v(t)\|_{H^{s}(\Theta)}\leq CN^{s-\sigma}\|u_0\|_{H^{\sigma}(\Theta)}\leq \tilde{C}N^{s-\sigma},\quad |t|\leq \tau.
\end{equation*}
This completes the proof of Lemma~\ref{kyoto}.
\end{proof}
Using 
(\ref{kenji}), Lemma~\ref{kyoto},
Lemma~\ref{lim2} and Proposition~\ref{liouville}, we can write
\begin{eqnarray*}
\rho\Big(\Phi(\tau)(K)+\overline{B_{2c\varepsilon}}\Big)& \geq &\limsup_{N\rightarrow \infty}
\rho_{N}\Big(\big(\Phi(\tau)(K)+\overline{B_{2c\varepsilon}}\big)\cap E_{N}\Big)
\\
& \geq &
\liminf_{N\rightarrow \infty}\rho_{N}\Big(\Phi_{N}(\tau)\big(\big(K+B_{\varepsilon}\big)\cap E_{N}\big)\Big)
\\
& = &
\liminf_{N\rightarrow \infty}\rho_{N}\Big(\big(K+B_{\varepsilon}\big)\cap E_{N}\Big)
\\
& \geq &
\rho\big(K+B_{\varepsilon}\big)\geq \rho(K).
\end{eqnarray*}
By letting $\varepsilon\rightarrow 0$, we obtain that $\rho(\Phi(\tau)(K))\geq \rho(K)$.
This completes the proof of (\ref{restriction}) which in turn completes the
proof of (\ref{vili}).
\\

This completes the proof of Theorem~\ref{glavna}.
\qed
\begin{remarque}
Let us notice that in the proof of Theorem~\ref{glavna}, we did not make appeal
to the conservation laws of (\ref{pak}). We only used the conservation laws of
(\ref{pakbis}) and thus the propagation of higher Sobolev regularity for (\ref{pak}) was not needed. 
\end{remarque} 
\section{Final remarks}
The result of Theorem~\ref{glavna} is obtained under the assumption
$\alpha<2$. Let us recall that if $\alpha=2$ with $F(u)=|u|^{2}u$ 
then one can construct initial data for (\ref{pak}) such that the local
solutions constructed in Proposition~\ref{lwp} develop singularities in finite time
(see \cite{Kav,BGT}). 
Observe that the data giving blow-up solutions in \cite{Kav} has to be
sufficiently smooth (at least $H^1$) in order to give sense of the quantities
involved in the well-known viriel identity. 
But one can show that for $\varphi_{\omega}$ defined by (\ref{5}) we have that
$\|\varphi_{\omega}\|_{H^1(\Theta)}$
is infinity almost surely.
It would be interesting to decide whether the obstruction to make work the
proof of Theorem~\ref{glavna} is related to a blow up phenomenon, i.e. can one
prove a blow up of the solutions of (\ref{pak}) with data on a set $A$ such
that $\rho(A)>0$ ?  A related and probably simpler question is whether one
can construct a blow up solution of NLS with data which is in $H^{s}$, $s<1$
but not in $H^1$ ?
\\

If we suppose the defocusing assumption $V(z)\leq 0$ then there is no problem
with the integrability of $f(u)$ and the $L^2$ cut-off is not needed. 
\\

Let us notice that the restriction $\alpha<2$ is too strong for the
well-posedness analysis of (\ref{1}) with data in ${\mathcal X}$.
Indeed this analysis seems to hold true for $\alpha<4$. Here is a rough
explanation. Essentially speaking, in order to
make work the nonlinear estimates with data of Sobolev regularity $<1/2$,
after $k\in\N$ expansions of the nonlinearity, for 
$$
N_2\geq N_3\geq \cdots\geq N_k
$$
one should control the expression
\begin{equation}\label{sob}
(N_2 N_3)^{\varepsilon}N_4\cdots N_k
\end{equation}
by
$$
C(N_2 N_3\cdots N_k)^{\sigma}
$$
for some $\sigma<1/2$. This leads to the restriction $k-3<\frac{1}{2}(k-1)$,
i.e. $k<5$ which corresponds to $\alpha<4$. 
In (\ref{sob}) the factor $N_4\cdots N_k$ appear from Sobolev embeddings which
in $2d$ costs $\frac{d}{2}=1$ derivatives (see \cite{BGTens} for a similar discussion).
However for $\alpha\geq 2$, the Sobolev inequality is no longer available to
give sense of $\int_{\Theta}V(u)$ for $u\in {\mathcal X}$. 
On the other hand one only needs to show that $\int_{\Theta}V(u)$ is finite
$\mu$ almost surely. This seems to be tractable by some Gaussian estimates and
the bounds of Lemma~\ref{Lp}. We plan to pursue this issue elsewhere.
\\

The measure $\mu$ constructed in Theorem~\ref{glavna} is obtained for
functions on the disc of radius $r=1$. Similar measures can be constructed for
any finite radius $r$ and the limiting behaviour of these measures as
$r\rightarrow\infty$ seems to be an interesting problem.
\\

One can also ask the question about ergodicity properties of the measure
$\rho$, i.e. the existence of ``non trivial'' $\rho$ measurable sets invariant
under the flow.
\\

Let us finally mention an extension of Theorem~\ref{glavna}. 
One can construct invariant measures leaving on
functions invariant by the rotations of the disc (see \cite{BGT}). In this
case, in the polar coordinates $(r,\varphi)$
on $\Theta$, the measure ``lives'' on the set of functions

\begin{equation}\label{higher}
e^{in\varphi}\,
\sum_{k\geq 1}\frac{g_{k}(\omega)}{z_{nk}}
\frac{J_{nk}(z_{nk}r)}{\|J_{nk}(z_{nk}\cdot)\|_{L^2(\Theta)}},
\end{equation}
where $J_{n}$, $n\geq 0$, $n\in\Z$ is the Bessel function of order $n$ and
$z_{nk}$ are its zeros (Theorem~\ref{glavna} corresponds to $n=0$). 
In (\ref{higher}), $g_{k}(\omega)$ is again a sequence of normalized
i.i.d. complex random variables.

\end{document}